**Курапов С. В.**

**Давидовский М. В.**

# АЛГОРИТМИЧЕСКИЕ МЕТОДЫ КОНЕЧНЫХ ДИСКРЕТНЫХ СТРУКТУР

# ТОПОЛОГИЧЕСКИЙ РИСУНОК ГРАФА

## часть 4

**(на правах рукописи)**

**Курапов С. В.**

**Давидовский М. В.**

В работе рассматриваются математические модели для создания топологического рисунка несепарабельного непланарного графа, основанные на методах теории вращения вершин Г. Рингеля. Индуцированная система циклов порождает топологический рисунок определенной толщины. Представлен метод определения местоположения мнимых вершин путём нахождения пересечения соединений на плоскости. В качестве основы используется топологический рисунок максимально плоского суграфа.

Для научных работников, студентов и аспирантов, специализирующихся на применении методов прикладной математики.

# ОГЛАВЛЕНИЕ

# Введение

Данная работа является дальнейшим развитием алгоритмических методов построения топологического рисунка графа, основанных на теории вращения вершин Г. Рингеля и Д. Янгса. Необходимые сведения о построении и описании топологического рисунка графа можно найти в следующих источниках информации: Курапов С. В., Давидовский М. В. Топологический рисунок графа (часть 1)[1], Курапов С. В., Давидовский М. В. Топологический рисунок графа (часть 2)[2], Курапов С. В., Давидовский М. В. Топологический рисунок графа (часть 3)[3].

Целью данной работы является построение и описание топологического рисунка непланарного несепарабельного графа. Составной частью описания топологического рисунка непланарного графа является понятие *толщины графа*, как наименьшего числа плоских подграфов, на которые можно разложить рёбра графа $G$. Мы будем исходить из положения, что для описания топологического рисунка неплоского графа нужно не только определить толщину, но и осуществить построение топологических рисунков плоских суграфов его составляющих.

Основой для построения топологического рисунка непланарного графа является система независимых изометрических циклов описывающая максимально плоский суграф графа $G$. В системе независимых циклов определяются сопряженные циклы. В результате образуется граф циклов $Y$, который ставит в соответствие циклу графа $G$ вершину графа циклов $Y$, а сопряженным ребрам графа $G$ ставит в соответствие ребра графа циклов $Y$. Учет связности вершин графа циклов $Y$ с вершинами графа $G$ порождает смешанный граф циклов $Q$.

Выделение из графа максимально плоского суграфа разделяет множество ребер графа $G$ на два подмножества: множество ребер принадлежащих максимально плоскому суграфу и множество ребер исключенных из графа в процессе планаризации.

Топологический рисунок максимально плоского суграфа определяет взаимное расположение вершин графа $G$, приводящее к исключению минимального количества ребер. Смешанный граф циклов $Q$, построенный на максимально плоском суграфе, позволяет определить маршруты для проведения исключенных соединений, тем самым задавая расположение соединений в топологическом пространстве. Процесс проведения соединения определяется местоположением мнимых вершин как топологическое пересечение соединения с ребрами максимально плоского суграфа. Подмножество циклов порожденное проведением непересекающихся соединений из подмножества исключенных в процессе планаризации, определяет топологический рисунок непланарного графа.

---

[1] https://arxiv.org/abs/2407.21564
[2] https://arxiv.org/abs/2407.21578
[3] https://arxiv.org/abs/2506.10936

Определение подмножества непересекающихся соединений для выбранной плоской части графа осуществляется методами векторной алгебры пересечений, разработанной группой математиков под руководством Раппопорта Л. И. Методы векторной алгебры пересечений позволяют определять пересечение соединений как пересечение проекций соединений на координатно-базисную систему векторов.

# Глава 16. ПЕРЕСЕЧЕНИЯ В НЕСЕПАРАБЕЛЬНОМ ГРАФЕ

## 16.1. Понятие толщины графа

Рассмотрим несепарабельные графы.

**Определение 16.1.** Связный неориентированный граф без петель и кратных рёбер, без мостов и точек сочленения, без вершин с валентностью меньшей трёх называется *несепарабельным графом*.

**Определение 16.2.** *Толщина* графа $G$ – это наименьшее число плоских подграфов, на которые можно разложить рёбра графа $G$.

То есть, если существует набор $k$ плоских графов, имеющих одинаковый набор вершин, объединение которых даёт граф $G$, то толщина графа $G$ не больше $k$. Таким образом, плоский граф имеет толщину $k=1$.

Данное определение не учитывает взаимное расположение вершин в пространстве и не рассматривает условие планарности подграфов.

Рассмотрим несепарабельный граф $K_7$.

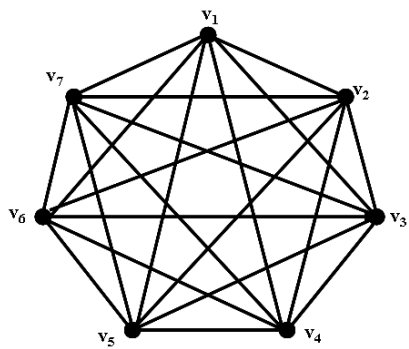

Рис. 16.1. Граф $K_7$.

Количество вершин графа = 7.
Количество рёбер графа = 21.
Количество единичных циклов = 35.

Смежность графа $K_7$:

вершина v1: v2 v3 v4 v5 v6 v7
вершина v2: v1 v3 v4 v5 v6 v7
вершина v3: v1 v2 v4 v5 v6 v7
вершина v4: v1 v2 v3 v5 v6 v7
вершина v5: v1 v2 v3 v4 v6 v7
вершина v6: v1 v2 v3 v4 v5 v7
вершина v7: v1 v2 v3 v4 v5 v6

Инцидентность графа $K_7$:

ребро $e_1$: $(v_1,v_2)$ или $(v_2,v_1)$;   ребро $e_2$: $(v_1,v_3)$ или $(v_3,v_1)$;
ребро $e_3$: $(v_1,v_4)$ или $(v_4,v_1)$;   ребро $e_4$: $(v_1,v_5)$ или $(v_5,v_1)$;
ребро $e_5$: $(v_1,v_6)$ или $(v_6,v_1)$;   ребро $e_6$: $(v_1,v_7)$ или $(v_7,v_1)$;
ребро $e_7$: $(v_2,v_3)$ или $(v_3,v_2)$;   ребро $e_8$: $(v_2,v_4)$ или $(v_4,v_2)$;

ребро $e_9$: $(v_2,v_5)$ или $(v_5,v_2)$;   ребро $e_{10}$: $(v_2,v_6)$ или $(v_6,v_2)$;
ребро $e_{11}$: $(v_2,v_7)$ или $(v_7,v_2)$;   ребро $e_{12}$: $(v_3,v_4)$ или $(v_4,v_3)$;
ребро $e_{13}$: $(v_3,v_5)$ или $(v_5,v_3)$;   ребро $e_{14}$: $(v_3.v_6)$ или $(v_6,v_3)$;
ребро $e_{15}$: $(v_3,v_7)$ или $(v_7,v_3)$;   ребро $e_{16}$: $(v_4.v_5)$ или $(v_5,v_4)$;
ребро $e_{17}$: $(v_4,v_6)$ или $(v_6,v_4)$;   ребро $e_{18}$: $(v_4,v_7)$ или $(v_7,v_4)$;
ребро $e_{19}$: $(v_5,v_6)$ или $(v_6,v_5)$;   ребро $e_{20}$: $(v_5.v_7)$ или $(v_7,v_5)$;
ребро $e_{21}$: $(v_6,v_7)$ или $(v_7.v_6)$.

Множество изометрических циклов графа:

цикл $c_1 = \{e_1,e_2,e_7\} \leftrightarrow \{v_1,v_2,v_3\}$;
цикл $c_2 = \{e_1,e_3,e_8\} \leftrightarrow \{v_1,v_2,v_4\}$;
цикл $c_3 = \{e_1,e_4,e_9\} \leftrightarrow \{v_1,v_2,v_5\}$;
цикл $c_4 = \{e_1,e_5,e_{10}\} \leftrightarrow \{v_1,v_2,v_6\}$;
цикл $c_5 = \{e_1,e_6,e_{11}\} \leftrightarrow \{v_1,v_2,v_7\}$;
цикл $c_6 = \{e_2,e_3,e_{12}\} \leftrightarrow \{v_1,v_3,v_4\}$;
цикл $c_7 = \{e_2,e_4,e_{13}\} \leftrightarrow \{v_1,v_3,v_5\}$;
цикл $c_8 = \{e_2,e_5,e_{14}\} \leftrightarrow \{v_1,v_3,v_6\}$;
цикл $c_9 = \{e_2,e_6,e_{15}\} \leftrightarrow \{v_1,v_3,v_7\}$;
цикл $c_{10} = \{e_3,e_4,e_{16}\} \leftrightarrow \{v_1,v_4,v_5\}$;
цикл $c_{11} = \{e_3,e_5,e_{17}\} \leftrightarrow \{v_1,v_4,v_6\}$;
цикл $c_{12} = \{e_3,e_6,e_{18}\} \leftrightarrow \{v_1,v_4,v_7\}$;
цикл $c_{13} = \{e_4,e_5,e_{19}\} \leftrightarrow \{v_1,v_5,v_6\}$;
цикл $c_{14} = \{e_4,e_6,e_{20}\} \leftrightarrow \{v_1,v_5,v_7\}$;
цикл $c_{15} = \{e_5,e_6,e_{21}\} \leftrightarrow \{v_1,v_6,v_7\}$;
цикл $c_{16} = \{e_7,e_8,e_{12}\} \leftrightarrow \{v_2,v_3,v_4\}$;
цикл $c_{17} = \{e_7,e_9,e_{13}\} \leftrightarrow \{v_2,v_3,v_5\}$;
цикл $c_{18} = \{e_7,e_{10},e_{14}\} \leftrightarrow \{v_2,v_3,v_6\}$;
цикл $c_{19} = \{e_7,e_{11},e_{15}\} \leftrightarrow \{v_2,v_3,v_7\}$;
цикл $c_{20} = \{e_8,e_9,e_{16}\} \leftrightarrow \{v_2,v_4,v_5\}$;
цикл $c_{21} = \{e_8,e_{10},e_{17}\} \leftrightarrow \{v_2,v_4,v_6\}$;
цикл $c_{22} = \{e_8,e_{11},e_{18}\} \leftrightarrow \{v_2,v_4,v_7\}$;
цикл $c_{23} = \{e_9,e_{10},e_{19}\} \leftrightarrow \{v_2,v_5,v_6\}$;
цикл $c_{24} = \{e_9,e_{11},e_{20}\} \leftrightarrow \{v_2,v_5,v_7\}$;
цикл $c_{25} = \{e_{10},e_{11},e_{21}\} \leftrightarrow \{v_2,v_6,v_7\}$;
цикл $c_{26} = \{e_{12},e_{13},e_{16}\} \leftrightarrow \{v_3,v_4,v_5\}$;
цикл $c_{27} = \{e_{12},e_{14},e_{17}\} \leftrightarrow \{v_3,v_4,v_6\}$;
цикл $c_{28} = \{e_{12},e_{15},e_{18}\} \leftrightarrow \{v_3,v_4,v_7\}$;
цикл $c_{29} = \{e_{13},e_{14},e_{19}\} \leftrightarrow \{v_3,v_5,v_6\}$;
цикл $c_{30} = \{e_{13},e_{15},e_{20}\} \leftrightarrow \{v_3,v_5,v_7\}$;
цикл $c_{31} = \{e_{14},e_{15},e_{21}\} \leftrightarrow \{v_3,v_6,v_7\}$;
цикл $c_{32} = \{e_{16},e_{17},e_{19}\} \leftrightarrow \{v_4,v_5,v_6\}$;
цикл $c_{33} = \{e_{16},e_{18},e_{20}\} \leftrightarrow \{v_4,v_5,v_7\}$;
цикл $c_{34} = \{e_{17},e_{18},e_{21}\} \leftrightarrow \{v_4,v_6,v_7\}$;
цикл $c_{35} = \{e_{19},e_{20},e_{21}\} \leftrightarrow \{v_5,v_6,v_7\}$.

Выделим систему изометрических циклов для построения максимально плоского суграфа.

цикл $c_1 = \{e_1,e_2,e_7\} \leftrightarrow \{v_1,v_2,v_3\}$;
цикл $c_5 = \{e_1,e_6,e_{11}\} \leftrightarrow \{v_1,v_2,v_7\}$;
цикл $c_{13} = \{e_4,e_5,e_{19}\} \leftrightarrow \{v_1,v_5,v_6\}$;
цикл $c_{15} = \{e_5,e_6,e_{21}\} \leftrightarrow \{v_1,v_6,v_7\}$;

цикл   $c_{19} = \{e_7, e_{11}, e_{15}\} \leftrightarrow \{v_2, v_3, v_7\}$;
цикл   $c_{26} = \{e_{12}, e_{13}, e_{16}\} \leftrightarrow \{v_3, v_4, v_5\}$;
цикл   $c_{28} = \{e_{12}, e_{15}, e_{18}\} \leftrightarrow \{v_3, v_4, v_7\}$;
цикл   $c_{33} = \{e_{16}, e_{18}, e_{20}\} \leftrightarrow \{v_4, v_5, v_7\}$;
цикл   $c_{35} = \{e_{19}, e_{20}, e_{21}\} \leftrightarrow \{v_5, v_6, v_7\}$.
обод $c_0$ = цикл $c_7 = \{e_2, e_4, e_{13}\} \leftrightarrow \{v_1, v_3, v_5\}$.

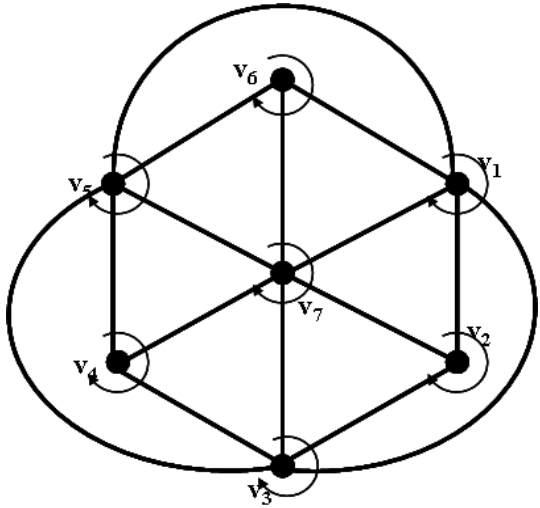 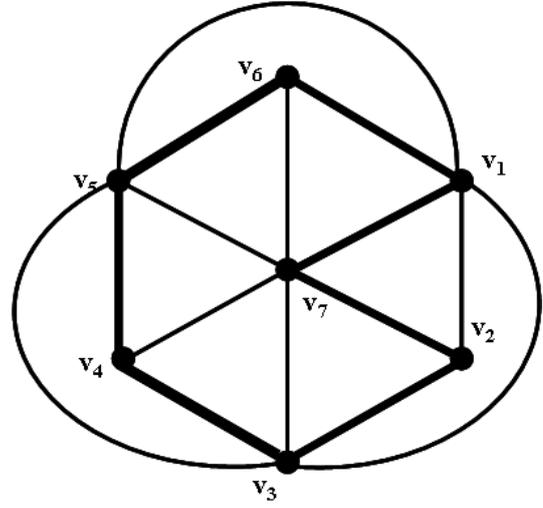

Рис. 16.2. Топологический рисунок максимально плоского суграфа.

Рис. 16.3. Гамильтонов цикл в максимально плоском суграфе.

Топологический рисунок расположения соединений, не вошедших в топологический рисунок максимально плоского суграфа, представлен на рис. 16.4. Введем исключенные в процессе планаризации соединения в топологический рисунок максимально плоского суграфа. Такие исключенные соединения показаны на рис. 16.5 красным цветом. Не трудно увидеть, что отличие в топологических рисунках заключается в различном взаимном расположении вершин.

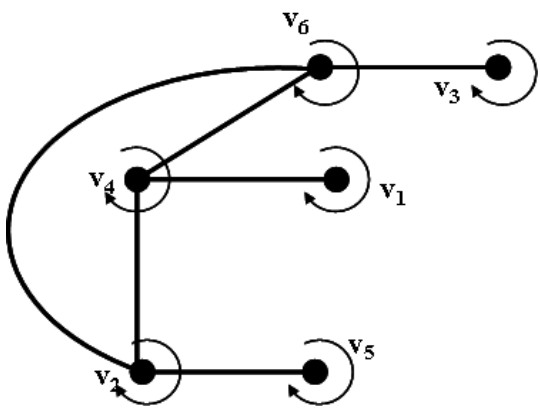 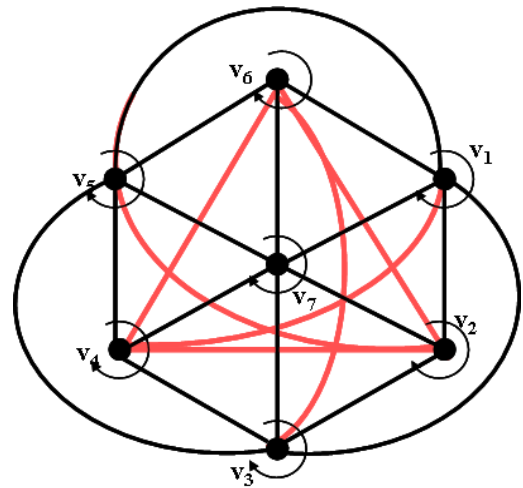

Рис. 16.4. Топологический рисунок исключенных соединений.

Рис. 16.5. Топологический рисунок плоского суграфа с исключенными соединениями.

Таким образом, Определение 16.2 характеризует свойства планарных составляющих графа $G$, не учитывая их связное относительное расположение элементов в пространстве.

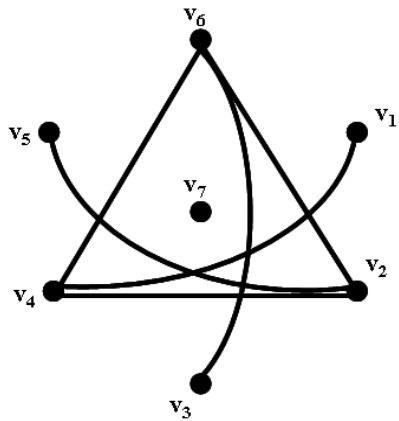

Рис. 16.6. Рисунок для ребер исключенных в процессе планаризации.

### 16.2. Гамильтонов цикл

Выделим в максимально плоском суграфе гамильтонов цикл (рис. 16.3). Система изометрических циклов, описывающая гамильтонов цикл, имеет вид:

цикл $c_{15} = \{e_5, e_6, e_{21}\} \leftrightarrow \{v_1, v_6, v_7\}$;
цикл $c_{19} = \{e_7, e_{11}, e_{15}\} \leftrightarrow \{v_2, v_3, v_7\}$;
цикл $c_{28} = \{e_{12}, e_{15}, e_{18}\} \leftrightarrow \{v_3, v_4, v_7\}$;
цикл $c_{33} = \{e_{16}, e_{18}, e_{20}\} \leftrightarrow \{v_4, v_5, v_7\}$;
цикл $c_{35} = \{e_{19}, e_{20}, e_{21}\} \leftrightarrow \{v_5, v_6, v_7\}$;
обод $c_0 = \{e_5, e_6, e_7, c_{11}.e_{12}, e_{16}, e_{19}\} \leftrightarrow \{v_1, v_2, v_3, v_4, v_5, v_6, v_7\}$.

Поставим в соответствие гамильтонову циклу координатно-базисную систему (КБС). Тогда для ребер, не вошедших в топологический рисунок плоской части графа, можно определить проекции на КБС (рис. 16.7).

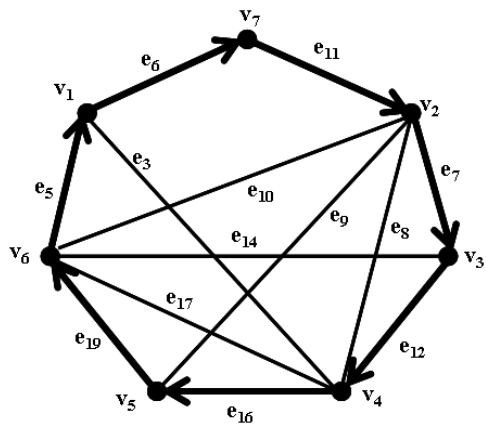

Рис. 16.7. Соединения и координатно-базисная система.

пр $e_3 = \{e_5, e_{16}, e_{19}\} = (v_4, v_5) + (v_5, v_6) + (v_6, v_1)$;
пр $e_8 = \{e_7, e_{12}\} = (v_2, v_3) + (v_3, v_4)$;
пр $e_9 = \{e_7, e_{12}, e_{16}\} = (v_2, v_3) + (v_3, v_4) + (v_4, v_5)$;
пр $e_{10} = \{e_5, e_6, e_{11}\} = (v_6, v_1) + (v_1, v_7) + (v_7, v_2)$;

пр $e_{14}$ = {$e_{12}, e_{16}, e_{19}$} = ($v_3, v_4$)+($v_4, v_5$)+($v_5, v_6$);
пр $e_{17}$ = {$e_{16}, e_{19}$} = ($v_4, v_5$)+($v_5, v_6$).

Определим пересечение соединений, как пересечение их проекций:

пр $e_3$ ∩ пр $e_8$ = {$e_5, e_{16}, e_{19}$} ∩ {$e_7, e_{12}$} = ∅;
пр $e_3$ ∩ пр $e_9$ = {$e_5, e_{16}, e_{19}$} ∩ {$e_7, e_{12}, e_{16}$} = $e_{16}$;
пр $e_3$ ∩ пр $e_{10}$ = {$e_5, e_{16}, e_{19}$} ∩ {$e_5, e_6, e_{11}$} = $e_5$;
пр $e_3$ ∩ пр $e_{14}$ = {$e_5, e_{16}, e_{19}$} ∩ {$e_{12}, e_{16}, e_{19}$} = {$e_{16}, e_{19}$};
пр $e_3$ ∩ пр $e_{17}$ = {$e_5, e_{16}, e_{19}$} ∩ {$e_{16}, e_{19}$} = пр $e_{17}$ ⊂ пр $e_3$ = ∅;
пр $e_8$ ∩ пр $e_9$ = {$e_7, e_{12}$} ∩ {$e_7, e_{12}, e_{16}$} = пр $e_8$ ⊂ пр $e_9$ = ∅;
пр $e_8$ ∩ пр $e_{10}$ = {$e_7, e_{12}$} ∩ {$e_5, e_6, e_{11}$} = ∅;
пр $e_8$ ∩ пр $e_{14}$ = {$e_7, e_{12}$} ∩ {$e_{12}, e_{16}, e_{19}$} = $e_{12}$;
пр $e_8$ ∩ пр $e_{17}$ = {$e_7, e_{12}$} ∩ {$e_{16}, e_{19}$} = ∅;
пр $e_9$ ∩ пр $e_{10}$ = {$e_7, e_{12}, e_{16}$} ∩ {$e_5, e_6, e_{11}$} = ∅;
пр $e_9$ ∩ пр $e_{14}$ = {$e_7, e_{12}, e_{16}$} ∩ {$e_{12}, e_{16}, e_{19}$} = {$e_{12}, e_{16}$};
пр $e_9$ ∩ пр $e_{17}$ = {$e_7, e_{12}, e_{16}$} ∩ {$e_{16}, e_{19}$} = $e_{16}$;
пр $e_{10}$ ∩ пр $e_{14}$ = {$e_5, e_6, e_{11}$} ∩ {$e_{12}, e_{16}, e_{19}$} = ∅;
пр $e_{10}$ ∩ пр $e_{17}$ = {$e_5, e_6, e_{11}$} ∩ {$e_{16}, e_{19}$} = ∅;
пр $e_{14}$ ∩ пр $e_{17}$ = {$e_{12}, e_{16}, e_{19}$} ∩ {$e_{16}, e_{19}$} = пр $e_{17}$ ⊂ пр $e_{14}$ = ∅.

Соединение e3 пересекается с 3-мя другими соединениями.
Соединение e8 пересекается с 1-им другим соединением.
Соединение e9 пересекается с 3-мя другими соединениями.
Соединение e10 пересекается с 1-им другим соединением.
Соединение e14 пересекается с 3-мя другими соединениями.
Соединение e17 пересекается с 1-им другим соединением.

Исключим ребро e3, как имеющее максимальное число пересечений и тогда изменится количество пересечений.

Соединение $e_8$ пересекается с 1-им другим соединением.
Соединение $e_9$ пересекается с 2-мя другими соединениями.
Соединение $e_{10}$ не имеет пересечений.
Соединение $e_{14}$ пересекается с 2-мя другими соединениями.
Соединение $e_{17}$ пересекается с 1-им другим соединением.

Исключим ребро $e_{14}$, как имеющее максимальное число пересечений, и тогда изменится количество пересечений.

Соединение $e_8$ не имеет пересечений.
Соединение e9 пересекается с 1-им другим соединением.
Соединение $e_{10}$ не имеет пересечений.
Соединение $e_{17}$ пересекается с 1-им другим соединением.

Исключим соединение $e_{17}$, как имеющее максимальное число пересечений, и тогда изменится количество пересечений.

Соединение $e_8$ = ($v_2, v_4$) не имеет пересечений.
Соединение $e_9$ = ($v_2, v_5$) не имеет пересечений.
Соединение $e_{10}$ = ($v_2, v_6$) не имеет пересечений.

С целью задания вращения вершин для описания плоского топологического рисунка зададим вращение изометрических циклов:

цикл $c_1$ = {$e_1, e_2, e_7$} ↔ <$v_1, v_3, v_2$> = ($v_1, v_3$)+($v_3, v_2$)+($v_2, v_1$);
цикл $c_5$ = {$e_1, e_6, e_{11}$} ↔ <$v_1, v_2, v_7$> = ($v_1, v_2$)+($v_2, v_7$)+($v_7, v_1$);

цикл $c_{13} = \{e_4,e_5,e_{19}\} \leftrightarrow \{v_1,v_6,v_5\} = (v_1,v_6)+(v_6,v_5)+(v_5,v_1)$;
цикл $c_{15} = \{e_5,e_6,e_{21}\} \leftrightarrow <v_1,v_7,v_6> = (v_1,v_7)+(v_7,v_6)+(v_6,v_1)$;
цикл $c_{19} = \{e_7,e_{11},e_{15}\} \leftrightarrow <v_2,v_3,v_7> = (v_2,v_3)+(v_3,v_7)+(v_7,v_2)$;
цикл $c_{26} = \{e_{12},e_{13},e_{16}\} \leftrightarrow <v_5,v_4,v_3> = (v_5,v_4)+(v_4,v_3)+(v_3,v_5)$;
цикл $c_{28} = \{e_{12},e_{15},e_{18}\} \leftrightarrow <v_3,v_4,v_7> = (v_3,v_4)+(v_4,v_7)+(v_7,v_3)$;
цикл $c_{33} = \{e_{16},e_{18},e_{20}\} \leftrightarrow <v_4,v_5,v_7> = (v_4,v_5)+(v_5,v_7)+(v_7,v_4)$;
цикл $c_{35} = \{e_{19},e_{20},e_{21}\} \leftrightarrow <v_5,v_6,v_7> = (v_5,v_6)+(v_6,v_7)+(v_7,v_5)$.
обод $c_0 =$ цикл $c_7 = \{e_2,e_4,e_{13}\} \leftrightarrow <v_5,v_3,v_1> = (v_5,v_3)+(v_3,v_1)+(v_1,v_5)$.

Поставим в соответствие каждому циклу в плоском топологическом рисунке вершину (рис. 16.8).

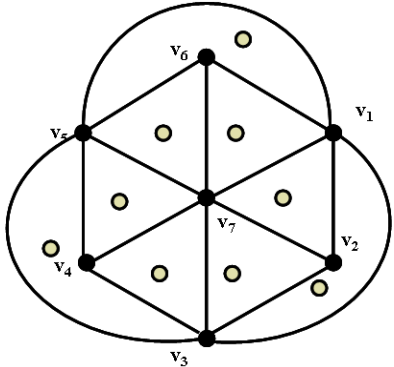 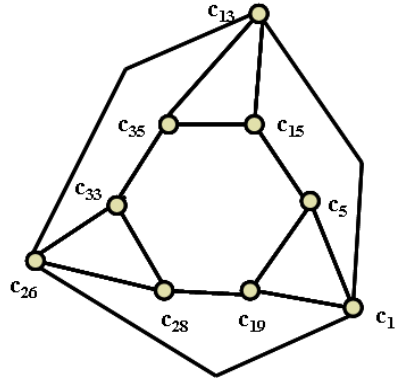 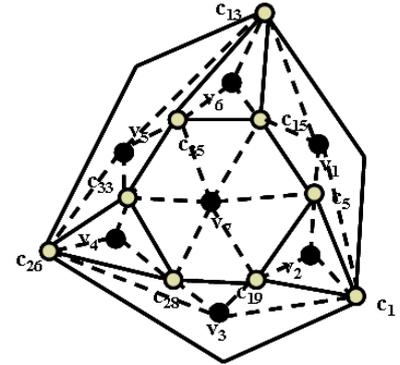

Рис. 16.8. Изометрические циклы и их вершины.    Рис. 16.9. Граф циклов.    Рис. 16.10. Смешанный граф циклов.

**Определение 16.3.** Два цикла, имеющие одно и только одно общее ребро, называются *сопряженными циклами*. Общее ребро называется *сопряженным ребром*.

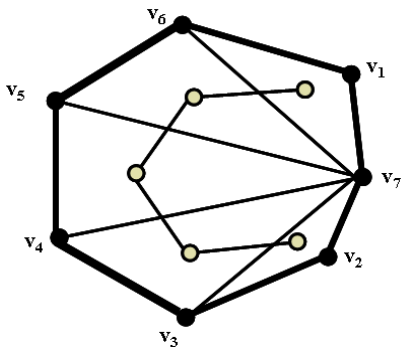 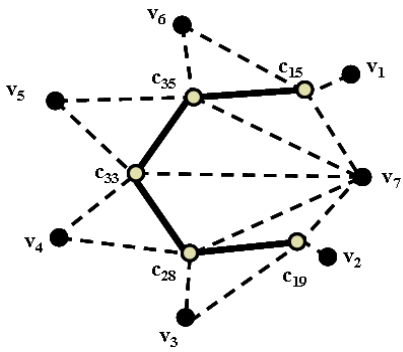 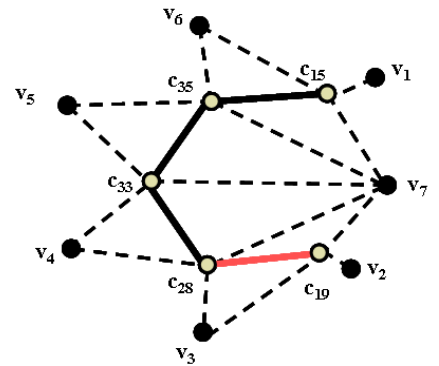

Рис. 16.11. Гамильтонов цикл и граф циклов.    Рис. 16.12. Смешанный граф циклов.    Рис. 16.13. Соединение циклов $c_{19}$ и $c_{28}$.

В плоском топологическом рисунке графа $K_7$ имеются циклы, содержащие ребро обода. В графе циклов соединим вершины таких циклов между собой ребрами. Также ребрами соединим вершины сопряженных циклов (рис. 16.9). Соединим вершины графа $K_7$ с вершинами циклов. Такой вид графа назовём *смешанным графом циклов* (рис. 16.10).

## 16.3. Смешанный граф циклов

Выделим гамильтонов цикл, состоящий из изометрических циклов $c_{19}, c_{28}, c_{33}, c_{35}, c_{15}$, и построим для него смешанный граф циклов (рис. 16.11-16.13).

Построим кратчайший маршрут в смешанном графе циклов между вершинами $v_2$ и $v_4$ (для ребра графа $e_8 = (v_2, v_4)$). В построенном маршруте будем рассматривать только ребра принадлежащие графу циклов. На рис. 16.13 данный маршрут состоит из одного ребра красного цвета.

Построенный маршрут соединяет вершины $v_2$ и $v_4$ и пересекает общее ребро циклов $c_{19}$ и $c_{28}$ ($\{e_7, e_{11}, e_{15}\} \cap \{e_{12}, e_{15}, e_{18}\} = e_{15} = (v_3, v_7)$). Обозначим это пересечение мнимой вершиной $v_8$. Рассмотрим векторную запись циклов $c_{19}$ и $c_{28}$.

$c_{19} = (v_2, v_3) + (v_3, v_7) + (v_7, v_2)$;
$c_{28} = (v_3, v_4) + (v_4, v_7) + (v_7, v_3)$.

Поместим вершину $v_8$ в ребро $e_{15}$. Векторную запись нужно начинать с концевых вершин:

$c_{19} = [(v_2, v_3) + (v_3, v_8)] + [(v_8, v_7) + (v_7, v_2)]$;
$c_{28} = [(v_4, v_7) + (v_7, v_8)] + [(v_8, v_3) + (v_3, v_4)]$.

Добавим в цикл $c_{19}$ ребро $(v_2, v_8)$ соединяющее вершину графа с мнимой вершиной. В цикл $c_{28}$ добавим ребро $(v_4, v_8)$ соединяющее вершину графа с мнимой вершиной.

$c_{19} = [(v_2, v_3) + (v_3, v_8)] + [(v_2.v_8) + (v_8, v_2)] + [(v_8, v_7) + (v_7, v_2)]$;
$c_{28} = [(v_4, v_7) + (v_7, v_8)] + [(v_4.v_8) + (v_8, v_4)] + [(v_8, v_3) + (v_3, v_4)]$.

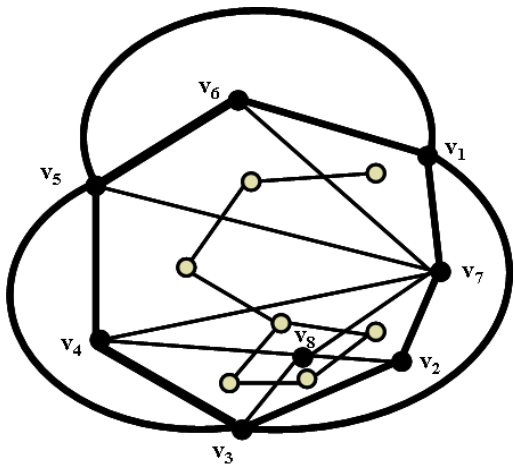 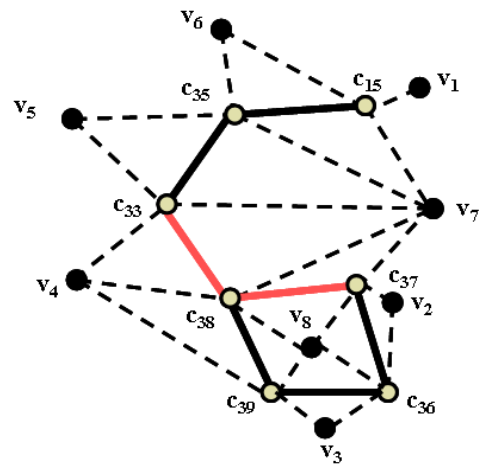

Рис. 16.14. Граф циклов.            Рис. 16.15. Смешанный граф циклов.

Выделим замкнутую последовательность ребер и сформируем новую систему циклов вместо циклов $c_{19}$ и $c_{28}$:

$c_{36} = (v_2, v_3) + (v_3, v_8) + (v_8.v_2)$;
$c_{37} = (v_2, v_8) + (v_8, v_7) + (v_7, v_2)$;
$c_{38} = (v_4, v_7) + (v_7, v_8) + (v_8.v_4)$;
$c_{39} = (v_4, v_8) + (v_8, v_3) + (v_3, v_4)$.

После операции введения в топологический рисунок мнимой вершины $v_8$ гамильтонов цикл состоит из циклов $c_{36}, c_{37}, c_{38}, c_{39}, c_{33}, c_{35}, c_{15}$. Построим смешанный граф циклов (рис. 16.15). Определим в смешанном графе минимальный маршрут для ребра $e_9 = (v_2, v_5)$. Минимальный маршрут состоит из двух ребер графа циклов $(c_{37}, c_{38})$ и $(c_{38}, c_{33})$. Новые мнимые вершины $v_9$ и $v_{10}$ будут расположены соответственно на пересечении циклов $c_{37} \cap c_{38}$ и $c_{38} \cap c_{33}$. Представим минимальный маршрут в виде кортежа вершин $<v_2, v_9, v_{10}, v_5>$.

Определим сопряженное неориентированное ребро как пересечение циклов $c_{37} \cap c_{38}$ с включением мнимой вершины $v_9$:

$c_{37} \cap c_{38} = [(v_2,v_8)+(v_8,v_7)+(v_7,v_2)] \cap [(v_4,v_7)+(v_7,v_8)+(v_8.v_4)] = (v_7,v_8)$.

Определим сопряженное неориентированное ребро как пересечение циклов $c_{38} \cap c_{33}$ с включением мнимой вершины $v_{10}$:

$c_{38} \cap c_{33} = [(v_4,v_7)+(v_7,v_8)+(v_8.v_4)] \cap [(v_4,v_5)+(v_5,v_7)+(v_7,v_4)] = (v_4,v_7)$.

Расположим циклы в последовательности:

$c_{37} = (v_2,v_8)+(v_8,v_7)+(v_7,v_2)$;
$c_{38} = (v_4,v_7)+(v_7,v_8)+(v_8.v_4)$;
$c_{33} = (v_5,v_7)+(v_7,v_4)+(v_4,v_5)$.

Вставляем в соответствующие ребра мнимые вершины:

$c_{37} = (v_2,v_8)+(v_8,v_9)+(v_9,v_7)+(v_7,v_2)$;
$c_{38} = (v_4,v_{10})+(v_{10},v_7)+(v_7,v_9)+(v_9,v_8)+(v_8.v_4)$;
$c_{33} = (v_5,v_7)+(v_7,v_{10})+(v_{10},v_4)+(v_4,v_5)$.

Упорядочим запись циклов оставляя на первом месте ориентированные ребра с вершинами из кортежа вершин минимального маршрута и отделяя ребра с мнимыми вершинами:

$c_{37} = [(v_2,v_8)+(v_8,v_9)]+[(v_9,v_7)+(v_7,v_2)]$;
$c_{38} = [(v_{10},v_7)+(v_7,v_9)]+[(v_9,v_8)+(v_8.v_4)+(v_4,v_{10})]$;
$c_{33} = [(v_5,v_7)+(v_7,v_{10})[+[(v_{10},v_4)+(v_4,v_5)]$.

Вставим отрезки соединения $(v_2,v_5)$:

$c_{37} = [(v_2,v_8)+(v_8,v_9)]+[(v_9,v_2)+(v_2,v_9)]+[(v_9,v_7)+(v_7,v_2)]$;
$c_{38} = [(v_{10},v_7)+(v_7,v_9)]+[(v_9,v_{10})+(v_{10},v_9)]+[(v_9,v_8)+(v_8.v_4)+(v_4,v_{10})]$;
$c_{33} = [(v_5,v_7)+(v_7,v_{10})]+[(v_{10},v_5)+(v_5,v_{10})]+[(v_{10},v_4)+(v_4,v_5)]$.

Выделяем замкнутую последовательность ориентированных ребер и формируем новую систему циклов вместо циклов $c_{37}, c_{38}, c_{33}$:

$c_{40} = (v_2,v_8)+(v_8,v_9)+(v_9,v_2)$;
$c_{41} = (v_2,v_9)+(v_9,v_7)+(v_7,v_2)$;
$c_{42} = (v_{10},v_7)+(v_7,v_9)+(v_9,v_{10})$;
$c_{43} = (v_{10},v_9)+(v_9,v_8)+(v_8.v_4)+(v_4,v_{10})$;
$c_{44} = (v_5,v_7)+(v_7,v_{10})+(v_{10},v_5)$
$c_{45} = (v_5,v_{10})+(v_{10},v_4)+(v_4,v_5)$.

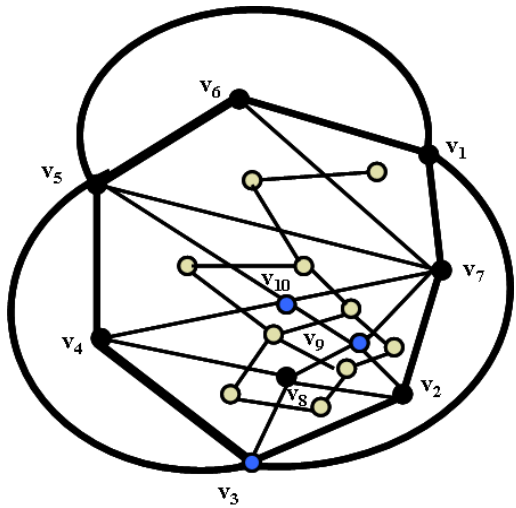
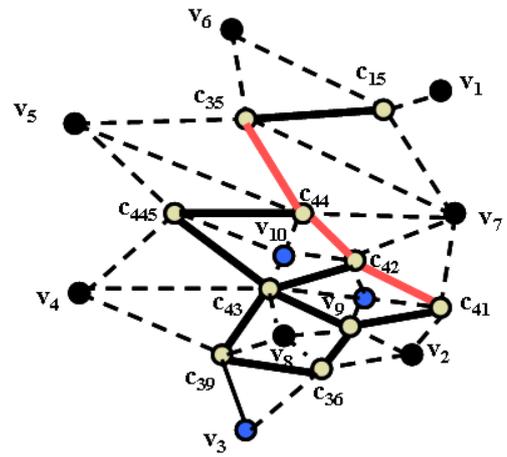

Рис. 16.16. Введение в топологический рисунок мнимых вершин $v_9, v_{10}$.

Рис. 16.17. Смешанный граф циклов.

В гамильтонов цикл входят следующие циклы:

$c_{15} = (v_1,v_7)+(v_7,v_6)+(v_6,v_1)$;
$c_{35} = (v_5,v_6)+(v_6,v_7)+(v_7,v_5)$;
$c_{36} = (v_2,v_3)+(v_3,v_8)+(v_8.v_2)$;
$c_{39} = (v_4,v_8)+(v_8,v_3)+(v_3,v_4)$;
$c_{40} = (v_2,v_8)+(v_8,v_9)+(v_9,v_2)$;
$c_{41} = (v_2,v_9)+(v_9,v_7)+(v_7,v_2)$;
$c_{42} = (v_{10},v_7)+(v_7,v_9)+(v_9,v_{10})$;
$c_{43} = (v_{10},v_9)+(v_9,v_8)+(v_8.v_4)+(v_4,v_{10})$;
$c_{44} = (v_5,v_7)+(v_7,v_{10})+(v_{10},v_5)$
$c_{45} = (v_5,v_{10})+(v_{10},v_4)+(v_4,v_5)$.

Определим минимальный маршрут для проведения ребра $(v_2,v_6)$ состоящий из последовательности вершин $c_{41}, c_{45}, c_{46}, c_{35}$ графа циклов. В сопряженное ребро циклов $c_{41}$ и $c_{42}$ будем вставлять вершину $v_{11}$, в сопряженное ребро циклов $c_{42}$ и $c_{44}$ – $v_{12}$, в сопряженное ребро циклов $c_{44}$ и $c_{35}$ – $v_{13}$. Минимальный маршрут запишем как кортеж вершин $<v_2, v_{11}, v_{12}, v_{13}, v_6>$.

Определим сопряженное неориентированное ребро как пересечение циклов $c_{41} \cap c_{45}$ с включением мнимой вершины $v_{11}$:

$c_{41} \cap c_{45} = [(v_2,v_9)+(v_9,v_7)+(v_7,v_2)] \cap [(v_9,v_{10})+(v_{10},v_7)+(v_7,v_9)] = (v_7,v_9)$.

Определим сопряженное неориентированное ребро как пересечение циклов $c_{45} \cap c_{46}$ с включением мнимой вершины $v_{12}$:

$c_{45} \cap c_{46} = [(v_9,v_{10})+(v_{10},v_7)+(v_7,v_9)] \cap [(v_5,v_7)+(v_7,v_{10})+(v_{10},v_5)] = (v_7,v_{10})$.

Определим сопряженное неориентированное ребро как пересечение циклов $c_{46} \cap c_{35}$ с включением мнимой вершины $v_{13}$:

$c_{46} \cap c_{35} = [(v_5,v_7)+(v_7,v_{10})+(v_{10},v_5)] \cap [(v_5,v_6)+(v_6,v_7)+(v_7,v_5)] = (v_5,v_7)$.

Рассмотрим следующую последовательность циклов:

$c_{41} = (v_2,v_9)+(v_9,v_7)+(v_7,v_2)$;
$c_{42} = (v_{10},v_7)+(v_7,v_9)+(v_9,v_{10})$;
$c_{44} = (v_5,v_7)+(v_7,v_{10})+(v_{10},v_5)$

$c_{35} = (v_5,v_6)+(v_6,v_7)+(v_7,v_5).$

Вставим мнимые вершины в соответствующие ребра:

$c_{41} = (v_2,v_9)+(v_9,v_{11})+(v_{11},v_7)+(v_7,v_2);$
$c_{42} = (v_{10},v_{12})+(v_{12},v_7)+(v_7,v_{11})+(v_{11},v_9)+(v_9,v_{10});$
$c_{44} = (v_5,v_{13})+(v_{13},v_7)+(v_7,v_{12})+(v_{12},v_{10})+(v_{10},v_5)$
$c_{35} = (v_5,v_6)+(v_6,v_7)+(v_7,v_{13})+(v_{13},v_5).$

Упорядочим запись циклов, оставляя на первом месте ориентированные ребра с вершинами из кортежа вершин минимального маршрута и отделяя ребра с мнимыми вершинами:

$c_{41} = [(v_2,v_9)+(v_9,v_{11})]+[(v_{11},v_7)+(v_7,v_2)];$
$c_{42} = [(v_{12},v_7)+(v_7,v_{11})]+[(v_{11},v_9)+(v_9,v_{10})+(v_{10},v_{12})];$
$c_{44} = [(v_{13},v_7)+(v_7,v_{12})]+[(v_{12},v_{10})+(v_{10},v_5)+(v_5,v_{13})];$
$c_{35} = [(v_6,v_7)+(v_7,v_{13})]+[(v_{13},v_5)+(v_5,v_6)].$

Вставим отрезки соединения $(v_2,v_6)$:

$c_{41} = [(v_2,v_9)+(v_9,v_{11})]+[(v_2,v_{11})+(v_{11},v_2)]+[(v_{11},v_7)+(v_7,v_2)];$
$c_{42} = [(v_{12},v_7)+(v_7,v_{11})]+[(v_{11},v_{12})+(v_{12},v_{11})]+[(v_{11},v_9)+(v_9,v_{10})+(v_{10},v_{12})];$
$c_{44} = [(v_{13},v_7)+(v_7,v_{12})]+[(v_{12},v_{13})+(v_{13},v_{12})]+[(v_{12},v_{10})+(v_{10},v_5)+(v_5,v_{13})];$
$c_{35} = [(v_6,v_7)+(v_7,v_{13})]+[(v_{13},v_6)+(v_6,v_{13})]+[(v_{13},v_5)+(v_5,v_6)].$

Выделяем замкнутую последовательность ориентированных ребер и формируем новую систему циклов вместо циклов $c_{37},c_{38},c_{33}$:

$c_{46} = (v_2,v_9)+(v_9,v_{11})+(v_2,v_{11});$
$c_{47} = (v_{11},v_2)+(v_{11},v_7)+(v_7,v_2);$
$c_{48} = (v_{12},v_7)+(v_7,v_{11})+(v_{11},v_{12});$
$c_{49} = (v_{12},v_{11})+(v_{11},v_9)+(v_9,v_{10})+(v_{10},v_{12});$
$c_{50} = (v_{13},v_7)+(v_7,v_{12})+(v_{12},v_{13});$
$c_{51} = (v_{13},v_{12})+(v_{12},v_{10})+(v_{10},v_5)+(v_5,v_{13});$
$c_{52} = (v_6,v_7)+(v_7,v_{13})+(v_{13},v_6);$
$c_{53} = (v_6,v_{13})+(v_{13},v_5)+(v_5,v_6).$

Подмножество циклов определяет гамильтонов цикл:

$c_{15} = (v_1,v_7)+(v_7,v_6)+(v_6,v_1);$
$c_{36} = (v_2,v_3)+(v_3,v_8)+(v_8.v_2);$
$c_{39} = (v_4,v_8)+(v_8,v_3)+(v_3,v_4);$
$c_{40} = (v_2,v_8)+(v_8,v_9)+(v_9,v_2);$
$c_{43} = (v_4,v_{10})+(v_{10},v_9)+(v_9,v_8)+(v_8.v_4);$
$c_{45} = (v_5,v_{10})+(v_{10},v_4)+(v_4,v_5);$
$c_{46} = (v_2,v_9)+(v_9,v_{11})+(v_2,v_{11});$
$c_{47} = (v_{11},v_2)+(v_{11},v_7)+(v_7,v_2);$
$c_{48} = (v_{12},v_7)+(v_7,v_{11})+(v_{11},v_{12});$
$c_{49} = (v_{12},v_{11})+(v_{11},v_9)+(v_9,v_{10})+(v_{10},v_{12});$
$c_{50} = (v_{13},v_7)+(v_7,v_{12})+(v_{12},v_{13});$
$c_{51} = (v_{13},v_{12})+(v_{12},v_{10})+(v_{10},v_5)+(v_5,v_{13});$
$c_{52} = (v_6,v_7)+(v_7,v_{13})+(v_{13},v_6);$
$c_{53} = (v_6,v_{13})+(v_{13},v_5)+(v_5,v_6).$

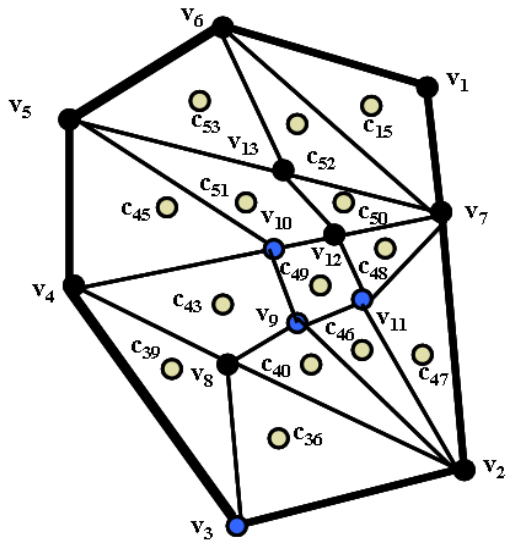 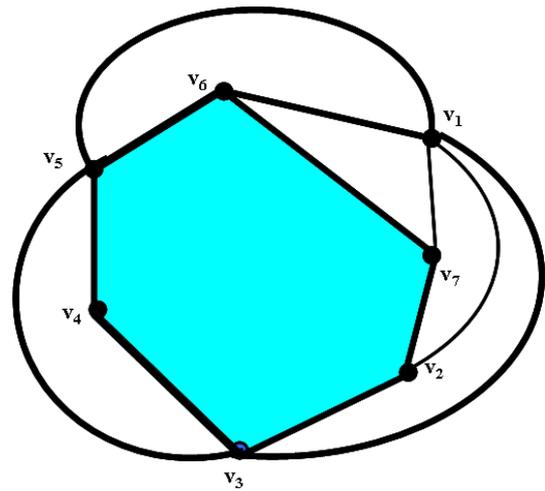

Рис. 16.18. Состав гамильтонова цикла после укладки ребра $(v_2,v_6)$.

Рис. 16.19. Удаление циклов с мнимыми вершинами.

Удалим в топологическом рисунке графа все циклы, имеющие в своем составе мнимые вершины (рис. 16.19). В результате удаления строится подмножество циклов с ободом $<v_6,v_1,v_3,v_2,v_7>$. Поставим в соответствие ободу подмножества циклов координатно-базисную систему векторов. Удаленное в процессе планаризации ребро $(v_3,v_6)$ не имеет пересечений в проекции на КБС (рис. 16.20).

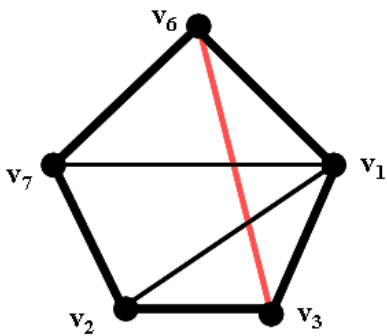 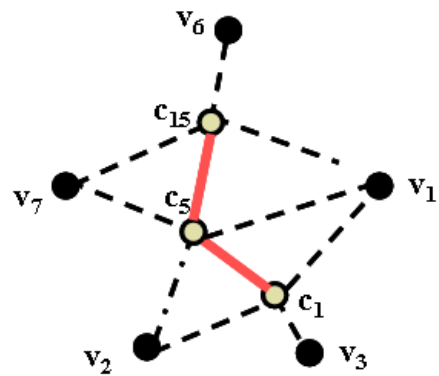

Рис. 16.20. Проведение соединения $(v_3,v_6)$.

Рис. 16.21. Смешанный граф циклов.

Выделим в смешанном графе минимальный маршрут, состоящий из отрезка $\{c_1,c_5\}$ с включением мнимой вершины $v_{14}$ и отрезка $\{c_5,c_{15}\}$ с включением мнимой вершины $v_{15}$. Запишем минимальный маршрут в виде кортежа вершин $<v_3,v_{14},v_{15},v_6>$.

Определим неориентированное сопряженное ребро как пересечение циклов $c_1 \cap c_5$ с включением мнимой вершины $v_{14}$:

$c_1 = (v_1,v_3)+(v_3,v_2)+(v_2,v_1)$;
$c_5 = (v_1,v_2)+(v_2,v_7)+(v_7,v_1)$.
$c_1 \cap c_5 = (v_2,v_1)$.

Определим неориентированное сопряженное ребро как пересечение циклов $c_5 \cap c_{15}$ с включением мнимой вершины $v_{15}$:

$c_5 = (v_1,v_2)+(v_2,v_7)+(v_7,v_1)$.
$c_{15} = (v_1,v_7)+(v_7,v_6)+(v_6,v_1)$.
$c_5\, c_{15} = (v_7,v_1)$.

Рассмотрим систему циклов:

$c_1 = (v_1,v_3)+(v_3,v_2)+(v_2,v_1)$;
$c_5 = (v_1,v_2)+(v_2,v_7)+(v_7,v_1)$.
$c_{15} = (v_1,v_7)+(v_7,v_6)+(v_6,v_1)$.

В соответствующие ребра вставляем мнимые вершины:

$c_1 = (v_1,v_3)+(v_3,v_2)+(v_2,v_{14})+(v_{14},v_1)$;
$c_5 = (v_1,v_{14})+(v_{14},v_2)+(v_2,v_7)+(v_7,v_{15})+(v_{15},v_1)$;
$c_{15} = (v_1,v_{15})+(v_{15},v_7)+(v_7,v_6)+(v_6,v_1)$.

Упорядочим запись циклов, оставляя на первом месте ориентированные ребра с вершинами из кортежа вершин минимального маршрута и отделяя ребра с мнимыми вершинами:

$c_1 = [(v_3,v_2)+(v_2,v_{14})]+[(v_{14},v_1)+(v_1,v_3)]$;
$c_5 = [(v_{14},v_2)+(v_2,v_7)+(v_7,v_{15})]+[(v_{15},v_1)+(v_1,v_{14})]$;
$c_{15} = [(v_6,v_1)+(v_1,v_{15})]+[(v_{15},v_7)+(v_7,v_6)]$.

Вставим отрезки соединения $(v_3,v_6)$:

$c_1 = [(v_3,v_2)+(v_2,v_{14})]+[(v_{14},v_3)+(v_3,v_{14})]+[(v_{14},v_1)+(v_1,v_3)]$;
$c_5 = [(v_{14},v_2)+(v_2,v_7)+(v_7,v_{15})]+[(v_{15},v_{14})+(v_{14},v_{15})]+[(v_{15},v_1)+(v_1,v_{14})]$;
$c_{15} = [(v_6,v_1)+(v_1,v_{15})]+[(v_{15},v_6)+(v_6,v_{15})]+[(v_{15},v_7)+(v_7,v_6)]$.

Выделим замкнутую последовательность ребер и сформируем новую систему циклов вместо циклов $c_1, c_5, c_{15}$:

$c_{54} = (v_3,v_2)+(v_2,v_{14})+(v_{14},v_3)$;
$c_{55} = (v_3,v_{14})+(v_{14},v_1)+(v_1,v_3)$;
$c_{56} = (v_{14},v_2)+(v_2,v_7)+(v_7,v_{15})+(v_{15},v_{14})$
$c_{57} = (v_{14},v_{15})+(v_{15},v_1)+(v_1,v_{14})$;
$c_{58} = (v_6,v_1)+(v_1,v_{15})+(v_{15},v_6)$;
$c_{59} = (v_6,v_{15})+(v_{15},v_7)+(v_7,v_6)$.

Формируем множество циклов, характеризующее топологический рисунок 2-го слоя (рис. 16.22):

$c_{13} = (v_6,v_5)+(v_5,v_1)+(v_1,v_6)$;
$c_{26} = (v_5,v_4)+(v_4,v_3)+(v_3,v_5)$;
$c_{36} = (v_2,v_3)+(v_3,v_8)+(v_8.v_2)$;
$c_{39} = (v_4,v_8)+(v_8,v_3)+(v_3,v_4)$;
$c_{40} = (v_2,v_8)+(v_8,v_9)+(v_9,v_2)$;
$c_{43} = (v_4,v_{10})+(v_{10},v_9)+(v_9,v_8)+(v_8.v_4)$;
$c_{45} = (v_5,v_{10})+(v_{10},v_4)+(v_4,v_5)$;
$c_{46} = (v_2,v_9)+(v_9,v_{11})+(v_2,v_{11})$;
$c_{47} = (v_{11},v_2)+(v_{11},v_7)+(v_7,v_2)$;
$c_{48} = (v_{12},v_7)+(v_7,v_{11})+(v_{11},v_{12})$;
$c_{49} = (v_{12},v_{11})+(v_{11},v_9)+(v_9,v_{10})+(v_{10},v_{12})$;

$c_{50} = (v_{13},v_7)+(v_7,v_{12})+(v_{12},v_{13})$;
$c_{51} = (v_{13},v_{12})+(v_{12},v_{10})+(v_{10},v_5)+(v_5,v_{13})$;
$c_{52} = (v_6,v_7)+(v_7,v_{13})+(v_{13},v_6)$;
$c_{53} = (v_6,v_{13})+(v_{13},v_5)+(v_5,v_6)$;
$c_{54} = (v_3,v_2)+(v_2,v_{14})+(v_{14},v_3)$;
$c_{55} = (v_3,v_{14})+(v_{14},v_1)+(v_1,v_3)$;
$c_{56} = (v_{14},v_2)+(v_2,v_7)+(v_7,v_{15})+(v_{15},v_{14})$;
$c_{57} = (v_{14},v_{15})+(v_{15},v_1)+(v_1,v_{14})$;
$c_{58} = (v_6,v_1)+(v_1,v_{15})+(v_{15},v_6)$;
$c_{59} = (v_6,v_{15})+(v_{15},v_7)+(v_7,v_6)$.

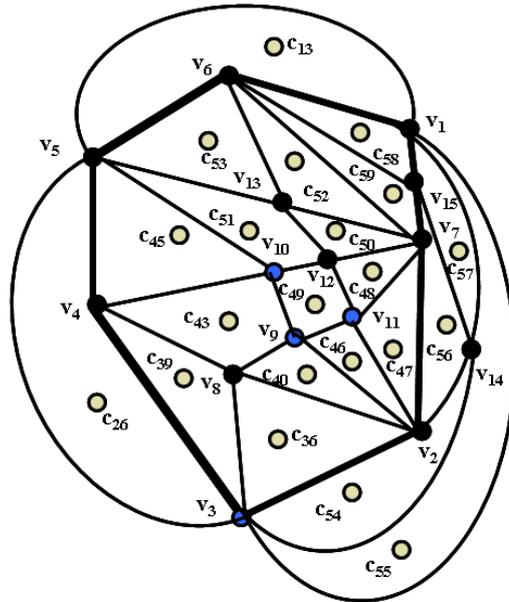

Рис. 16.22. Топологический рисунок 2-го слоя графа *G*.

### 16.4. Топологический рисунок следующего слоя

После построения топологического рисунка 2-го слоя остались два не уложенных соединения – $(v_1,v_4)$ и $(v_4,v_6)$. Для проведения оставшихся соединений снова выделим в плоском топологическом рисунке гамильтонов цикл. Поставим в соответствие гамильтонову циклу координатно-базисную систему векторов и убедимся, что соединения не пересекаются (рис. 16.23).

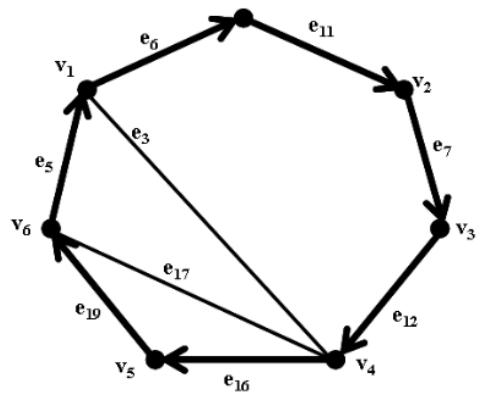

Рис. 16.23. Координатно-базисная система векторов.

Выделим минимальный маршрут в смешанном графе циклов (рис. 16.24). Минимальный маршрут состоит из отрезка ($c_{45}$,$c_{51}$) с введением мнимой вершины $v_{16}$ и отрезка ($c_{51}$,$c_{53}$) с введением мнимой вершины $v_{17}$. Запишем минимальный маршрут в виде кортежа вершин <$v_4$,$v_{16}$,$v_{17}$,$v_6$>.

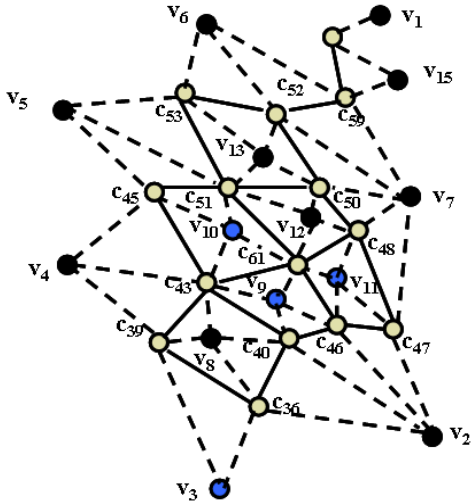 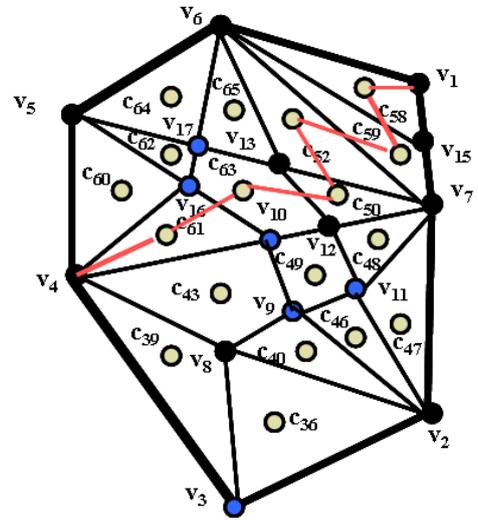

Рис. 16.24. Смешанный граф циклов.  Рис. 16.25. Проведение соединения ($v_4$,$v_1$).

Определим сопряженное неориентированное ребро как пересечение циклов $c_{41} \cap c_{45}$ с включением мнимой вершины $v_{16}$:

$c_{45} \cap c_{51}$ = [($v_5$,$v_{10}$)+($v_{10}$,$v_4$)+($v_4$,$v_5$)] $\cap$ [($v_{13}$,$v_{12}$)+($v_{12}$,$v_{10}$)+($v_{10}$,$v_5$)+($v_5$,$v_{13}$)] = ($v_5$,$v_{10}$);

Определим сопряженное неориентированное ребро как пересечение циклов $c_{45} \cap c_{46}$ с включением мнимой вершины $v_{17}$:

$c_{51} \cap c_{53}$ = [($v_{13}$,$v_{12}$)+($v_{12}$,$v_{10}$)+($v_{10}$,$v_5$)+($v_5$,$v_{13}$)] $\cap$ [($v_6$,$v_{13}$)+($v_{13}$,$v_5$)+($v_5$,$v_6$)] = ($v_5$,$v_{13}$).
$c_{45}$ = ($v_5$,$v_{10}$)+($v_{10}$,$v_4$)+($v_4$,$v_5$);
$c_{51}$ = ($v_{13}$,$v_{12}$)+($v_{12}$,$v_{10}$)+($v_{10}$,$v_5$)+($v_5$,$v_{13}$);
$c_{53}$ = ($v_6$,$v_{13}$)+($v_{13}$,$v_5$)+($v_5$,$v_6$).

Введем мнимые вершины в соответствующие ребра:

$c_{45}$ = ($v_5$,$v_{16}$)+($v_{16}$,$v_{10}$)+($v_{10}$,$v_4$)+($v_4$,$v_5$);
$c_{51}$ = ($v_{13}$,$v_{12}$)+($v_{12}$,$v_{10}$)+($v_{10}$,$v_{16}$)+($v_{16}$,$v_5$)+($v_5$,$v_{17}$)+($v_{17}$,$v_{13}$);
$c_{53}$ = ($v_6$,$v_{13}$)+($v_{13}$,$v_{17}$)+($v_{17}$,$v_5$)+($v_5$,$v_6$).

Упорядочим запись циклов, оставляя на первом месте ориентированные ребра с вершинами из кортежа вершин минимального маршрута и отделяя ребра с мнимыми вершинами:

$c_{45}$ = [($v_4$,$v_5$)+($v_5$,$v_{16}$)]+[($v_{16}$,$v_{10}$)+($v_{10}$,$v_4$)];
$c_{51}$ = [($v_{16}$,$v_5$)+($v_5$,$v_{17}$)]+[($v_{17}$,$v_{13}$)+($v_{13}$,$v_{12}$)+($v_{12}$,$v_{10}$)+($v_{10}$,$v_{16}$)];
$c_{53}$ = [($v_6$,$v_{13}$)+($v_{13}$,$v_{17}$)]+[($v_{17}$,$v_5$)+($v_5$,$v_6$)].

Вставим отрезки соединения ($v_4$,$v_6$):

$c_{45}$ = [($v_4$,$v_5$)+($v_5$,$v_{16}$)]+[($v_{16}$,$v_4$)+($v_4$,$v_{16}$)]+[($v_{16}$,$v_{10}$)+($v_{10}$,$v_4$)];
$c_{51}$ = [($v_{16}$,$v_5$)+($v_5$,$v_{17}$)]+[($v_{17}$,$v_{16}$)+($v_{16}$,$v_{17}$)]+[($v_{17}$,$v_{13}$)+($v_{13}$,$v_{12}$)+($v_{12}$,$v_{10}$)+($v_{10}$,$v_{16}$)];
$c_{53}$ = [($v_6$,$v_{13}$)+($v_{13}$,$v_{17}$)]+[($v_{17}$,$v_6$)+($v_6$,$v_{17}$)]+[($v_{17}$,$v_5$)+($v_5$,$v_6$)].

Выделим замкнутую последовательность рёбер и сформируем новую систему циклов вместо циклов $c_{45}, c_{51}, c_{53}$.

$c_{60} = (v_4,v_5)+(v_5,v_{16})+(v_{16},v_4);$
$c_{61} = (v_4,v_{16})+(v_{16},v_{10})+(v_{10},v_4);$
$c_{62} = (v_{16},v_5)+(v_5,v_{17})+(v_{17},v_{16});$
$c_{63} = (v_{16},v_{17})+(v_{17},v_{13})+(v_{13},v_{12})+(v_{12},v_{10})+(v_{10},v_{16});$
$c_{64} = (v_6,v_{13})+(v_{13},v_{17})+(v_{17},v_6);$
$c_{65} = (v_6,v_{17})+(v_{17},v_5)+(v_5,v_6).$

Строим смешанный граф циклов. Выделим минимальный маршрут в смешанном графе циклов. Этот маршрут состоит из отрезка $(c_{61},c_{63})$ с введением мнимой вершины $v_{18}$, отрезка $(c_{63},c_{50})$ с введением мнимой вершины $v_{19}$, отрезка $(c_{50},c_{52})$ с введением мнимой вершины $v_{20}$, отрезка $(c_{52},c_{59})$ с введением мнимой вершины $v_{21}$, отрезка $(c_{59},c_{58})$ с введением мнимой вершины $v_{22}$. Запишем минимальный маршрут в виде кортежа вершин $<v_4,v_{18},v_{19},v_{20},v_{21},v_{22},v_6>$.

Определим сопряженное неориентированное ребро как пересечение циклов $c_{61} \cap c_{63}$ с включением мнимой вершины $v_{18}$:

$c_{61} \cap c_{63} = [(v_4,v_{16})+(v_{16},v_{10})+(v_{10},v_4)] \cap [(v_{16},v_{17})+(v_{17},v_{13})+(v_{13},v_{12})+(v_{12},v_{10})+(v_{10},v_{16})] =$
$= (v_{16},v_{10}).$

Определим сопряженное неориентированное ребро как пересечение циклов $c_{61} \cap c_{63}$ с включением мнимой вершины $v_{19}$:

$c_{63} \cap c_{50} = [(v_{16},v_{17})+(v_{17},v_{13})+(v_{13},v_{12})+(v_{12},v_{10})+(v_{10},v_{16})] \cap [(v_{13},v_7)+(v_7,v_{12})+(v_{12},v_{13})] =$
$= (v_{12},v_{13}).$

Определим сопряженное неориентированное ребро как пересечение циклов $c_{61} \cap c_{63}$ с включением мнимой вершины $v_{20}$:

$c_{50} \cap c_{52} = [(v_{13},v_7)+(v_7,v_{12})+(v_{12},v_{13})] \cap [(v_6,v_7)+(v_7,v_{13})+(v_{13},v_6)] = (v_7,v_{13}).$

Определим сопряженное неориентированное ребро как пересечение циклов $c_{61} \cap c_{63}$ с включением мнимой вершины $v_{21}$:

$c_{52} \cap c_{59} = [(v_6,v_7)+(v_7,v_{13})+(v_{13},v_6)] \cap [(v_6,v_{15})+(v_{15},v_7)+(v_7,v_6)] = (v_7,v_6).$

Определим сопряженное неориентированное ребро как пересечение циклов $c_{61} \cap c_{63}$ с включением мнимой вершины $v_{22}$:

$c_{59} \cap c_{58} = [(v_6,v_{15})+(v_{15},v_7)+(v_7,v_6)] \cap [(v_6,v_1)+(v_1,v_{15})+(v_{15},v_6)] = (v_{15},v_6).$

Рассмотрим следующую последовательность циклов:

$c_{61} = (v_4,v_{16})+(v_{16},v_{10})+(v_{10},v_4);$
$c_{63} = (v_{16},v_{17})+(v_{17},v_{13})+(v_{13},v_{12})+(v_{12},v_{10})+(v_{10},v_{16});$
$c_{50} = (v_{13},v_7)+(v_7,v_{12})+(v_{12},v_{13});$
$c_{52} = (v_6,v_7)+(v_7,v_{13})+(v_{13},v_6);$
$c_{59} = (v_6,v_{15})+(v_{15},v_7)+(v_7,v_6);$
$c_{58} = (v_6,v_1)+(v_1,v_{15})+(v_{15},v_6).$

Введём мнимые вершины в соответствующие рёбра:

$c_{61} = (v_4,v_{16})+(v_{16},v_{18})+(v_{18},v_{10})+(v_{10},v_4);$

$c_{63} = (v_{16},v_{17})+(v_{17},v_{13})+(v_{13},v_{19})+(v_{19},v_{12})+(v_{12},v_{10})+(v_{10},v_{18})+(v_{18}.v_{16})$;
$c_{50} = (v_{13},v_{20})+(v_{20},v_{7})+(v_{7},v_{12})+(v_{12},v_{19})+(v_{19},v_{13})$;
$c_{52} = (v_{6},v_{21})+(v_{21},v_{7})+(v_{7},v_{20})+(v_{20},v_{13})+(v_{13},v_{6})$;
$c_{59} = (v_{6},v_{22})+(v_{22},v_{15})+(v_{15},v_{7})+(v_{7},v_{21})+(v_{21},v_{6})$;
$c_{58} = (v_{6},v_{1})+(v_{1},v_{15})+(v_{15},v_{22})+(v_{22},v_{6})$.

Упорядочим запись циклов, оставляя на первом месте ориентированные ребра с вершинами из кортежа вершин минимального маршрута и отделяя ребра с мнимыми вершинами:

$c_{61} = [(v_{4},v_{16})+(v_{16},v_{18})]+[(v_{18},v_{10})+(v_{10},v_{4})]$;
$c_{63} = [(v_{18}.v_{16})+(v_{16},v_{17})+(v_{17},v_{13})+(v_{13},v_{19})]+[(v_{19},v_{12})+(v_{12},v_{10})+(v_{10},v_{18})]$;
$c_{50} = [(v_{19},v_{13}) (v_{13},v_{20})]+[(v_{20},v_{7})+(v_{7},v_{12})+(v_{12},v_{19})]$;
$c_{52} = [(v_{20},v_{13})+(v_{13},v_{6})+(v_{6},v_{21})]+[(v_{21},v_{7})+(v_{7},v_{20})]$;
$c_{59} = [(v_{21},v_{6}) (v_{6},v_{22})]+[(v_{22},v_{15})+(v_{15},v_{7})+(v_{7},v_{21})]$;
$c_{58} = [(v_{1},v_{15})+(v_{15},v_{22})]+[(v_{22},v_{6})+(v_{6},v_{1})]$.

Вставим отрезки соединения $(v_{4},v_{1})$:

$c_{61} = [(v_{4},v_{16})+(v_{16},v_{18})]+ [(v_{18},v_{4})+(v_{4},v_{18})]+[(v_{18},v_{10})+(v_{10},v_{4})]$;
$c_{63} = [(v_{18}.v_{16})+(v_{16},v_{17})+(v_{17},v_{13})+(v_{13},v_{19})]+ [(v_{19},v_{18})+(v_{18},v_{19})]+[(v_{19},v_{12})+(v_{12},v_{10})+ +(v_{10},v_{18})]$;
$c_{50} = [(v_{19},v_{13}) (v_{13},v_{20})]+[(v_{20},v_{19})+(v_{19},v_{20})]+[(v_{20},v_{7})+(v_{7},v_{12})+(v_{12},v_{19})]$;
$c_{52} = [(v_{20},v_{13})+(v_{13},v_{6})+(v_{6},v_{21})]+ [(v_{21},v_{20})+(v_{20},v_{21})]+[(v_{21},v_{7})+(v_{7},v_{20})]$;
$c_{59} = [(v_{21},v_{6}) (v_{6},v_{22})]+[(v_{22},v_{21})+(v_{21},v_{22})]+[(v_{22},v_{15})+(v_{15},v_{7})+(v_{7},v_{21})]$;
$c_{58} = [(v_{1},v_{15})+(v_{15},v_{22})]+[(v_{22},v_{1})+(v_{1},v_{22})]+[(v_{22},v_{6})+(v_{6},v_{1})]$.

Выделим замкнутую последовательность ребер и сформируем новую систему циклов вместо циклов $c_{61}, c_{63}, c_{50}, c_{52}, c_{59}, c_{58}$:

$c_{66} = (v_{4},v_{16})+(v_{16},v_{18})+(v_{18},v_{4})$;
$c_{67} = (v_{4},v_{18})+(v_{18},v_{10})+(v_{10},v_{4})$;
$c_{68} = (v_{18}.v_{16})+(v_{16},v_{17})+(v_{17},v_{13})+(v_{13},v_{19})+(v_{19},v_{18})$;
$c_{69} = (v_{18},v_{19})+(v_{19},v_{12})+(v_{12},v_{10})+(v_{10},v_{18})$;
$c_{70} = (v_{19},v_{13})+(v_{13},v_{20})+(v_{20},v_{19})$;
$c_{71} = (v_{19},v_{20})+(v_{20},v_{7})+(v_{7},v_{12})+(v_{12},v_{19})$;
$c_{72} = (v_{20},v_{13})+(v_{13},v_{6})+(v_{6},v_{21})+(v_{21},v_{20})$;
$c_{73} =+(v_{20},v_{21})]+[(v_{21},v_{7})+(v_{7},v_{20})]$;
$c_{74} = [(v_{21},v_{6}) (v_{6},v_{22})]+[(v_{22},v_{21})+$
$c_{75} = (v_{21},v_{22})]+[(v_{22},v_{15})+(v_{15},v_{7})+(v_{7},v_{21})]$;
$c_{76} = (v_{1},v_{15})+(v_{15},v_{22})+(v_{22},v_{1})$;
$c_{77} = (v_{1},v_{22})+(v_{22},v_{6})+(v_{6},v_{1})$.

Формируем множество циклов, характеризующее топологический рисунок 3-го слоя (рис. 16.26):

$c_{13} = (v_{6},v_{5})+(v_{5},v_{1})+(v_{1},v_{6})$;
$c_{26} = (v_{5},v_{4})+(v_{4},v_{3})+(v_{3},v_{5})$;
$c_{36} = (v_{2},v_{3})+(v_{3},v_{8})+(v_{8}.v_{2})$;
$c_{39} = (v_{4},v_{8})+(v_{8},v_{3})+(v_{3},v_{4})$;
$c_{40} = (v_{2},v_{8})+(v_{8},v_{9})+(v_{9},v_{2})$;
$c_{43} = (v_{4},v_{10})+(v_{10},v_{9})+(v_{9},v_{8})+(v_{8}.v_{4})$;
$c_{46} = (v_{2},v_{9})+(v_{9},v_{11})+(v_{2},v_{11})$;
$c_{47} = (v_{11},v_{2})+(v_{11},v_{7})+(v_{7},v_{2})$;
$c_{48} = (v_{12},v_{7})+(v_{7},v_{11})+(v_{11},v_{12})$;

$c_{49} = (v_{12},v_{11})+(v_{11},v_9)+(v_9,v_{10})+(v_{10},v_{12})$;
$c_{54} = (v_3,v_2)+(v_2,v_{14})+(v_{14},v_3)$;
$c_{55} = (v_3,v_{14})+(v_{14},v_1)+(v_1,v_3)$;
$c_{56} = (v_{14},v_2)+(v_2,v_7)+(v_7,v_{15})+(v_{15},v_{14})$;
$c_{57} = (v_{14},v_{15})+(v_{15},v_1)+(v_1,v_{14})$;
$c_{60} = (v_4,v_5)+(v_5,v_{16})+(v_{16},v_4)$;
$c_{62} = (v_{16},v_5)+(v_5,v_{17})+(v_{17},v_{16})$;
$c_{64} = (v_6,v_{13})+(v_{13},v_{17})+(v_{17},v_6)$;
$c_{65} = (v_6,v_{17})+(v_{17},v_5)+(v_5,v_6)$.
$c_{66} = (v_4,v_{16})+(v_{16},v_{18})+(v_{18},v_4)$;
$c_{67} = (v_4,v_{18})+(v_{18},v_{10})+(v_{10},v_4)$;
$c_{68} = (v_{18}.v_{16})+(v_{16},v_{17})+(v_{17},v_{13})+(v_{13},v_{19})+(v_{19},v_{18})$;
$c_{69} = (v_{18},v_{19})+(v_{19},v_{12})+(v_{12},v_{10})+(v_{10},v_{18})$;
$c_{70} = (v_{19},v_{13})+(v_{13},v_{20})+(v_{20},v_{19})$;
$c_{71} = (v_{19},v_{20})+(v_{20},v_7)+(v_7,v_{12})+(v_{12},v_{19})$;
$c_{72} = (v_{20},v_{13})+(v_{13},v_6)+(v_6,v_{21})+(v_{21},v_{20})$;
$c_{73} =+(v_{20},v_{21})+(v_{21},v_7)+(v_7,v_{20})]$;
$c_{74} = (v_{21},v_6)\ (v_6,v_{22})+(v_{22},v_{21})$;
$c_{75} = (v_{21},v_{22})+(v_{22},v_{15})+(v_{15},v_7)+(v_7,v_{21})]$;
$c_{76} = (v_1,v_{15})+(v_{15},v_{22})+(v_{22},v_1)$;
$c_{77} = (v_1,v_{22})+(v_{22},v_6)+(v_6,v_1)$.

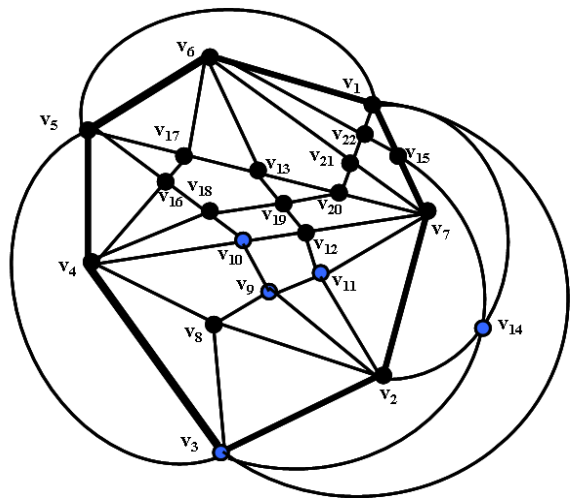

Рис. 16.26. Топологический рисунок графа с 3-мя слоями.

Рассмотрим последовательность расположения мнимых вершин, расположенных на соединениях удаленных в процессе планаризации:

$<v_1.v_4> = <v_{22},v_{21},v_{20},v_{19},v_{18}>$;
$<v_2.v_4> = <v_8>$;
$<v_2.v_5> = <v_9,v_{10},v_{18},v_{16}>$;
$<v_2.v_6> = <v_{11},v_{12},v_{19},v_{13}>$;
$<v_4.v_6> = <v_{16},v_{17}>$;
$<v_3.v_6> = <v_{14},v_{15},v_{22}>$.

**Комментарии**

В этой главе показано, что описание процесса построения топологического рисунка непланарного графа начинается с выделения максимально плоского суграфа графа.

Выделяется подмножество рёбер исключённых из графа в процессе планаризации. На примере графа К$_7$ показан процесс определения подмножества непересекающихся соединений внутри гамильтонова цикла. Процесс реализуется методами векторной алгебры пересечений.

Также в главе вводятся понятия *сопряжённого ребра* и *сопряжённых циклов*. С использованием этих понятий строится смешанный граф циклов для подмножества циклов, характеризующих гамильтонов цикл. Для проведения соединения осуществляется поиск минимального маршрута, состоящего из подмножества циклов. Проведение соединения разделяет выделенные циклы на соответствующие части путём введения мнимых вершин. В результате образуется новая система циклов, описывающая топологический рисунок суграфа с учётом проведения соединений.

Кроме того, в главе показано, что процесс проведения соединений не ограничивается построениями внутри гамильтонова цикла. В частности, процесс построения координатно-базисной системы векторов для проведения непересекающихся соединений может состоять из обода для любой независимой системы изометрических циклов. Главное, чтобы общая кольцевая сумма циклов была линейно независима.

При невозможности дальнейшего проведения соединений осуществляется построение топологического рисунка следующего слоя.

# Глава 17. ТОПОЛОГИЧЕСКИЙ РИСУНОК НЕПЛАНАРНОГО ГРАФА

## 17.1. Обод максимально плоского суграфа

За основу примем топологический рисунок 2-го слоя. Обод топологического рисунка 2-го слоя разбивает пространство $R^2$ на две области (части): внешнюю и внутреннею. Мы рассмотрели случай проведения пересекающихся соединений во внутренней области и перешли к проведению соединений в следующем слое. Рассмотрим далее проведение пересекающихся соединений во внешней области.

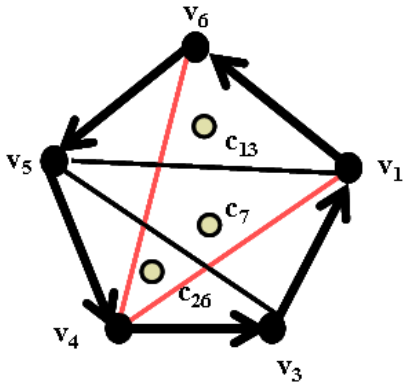
Рис. 17.1. Координатно-базисная система векторов для внешней области.

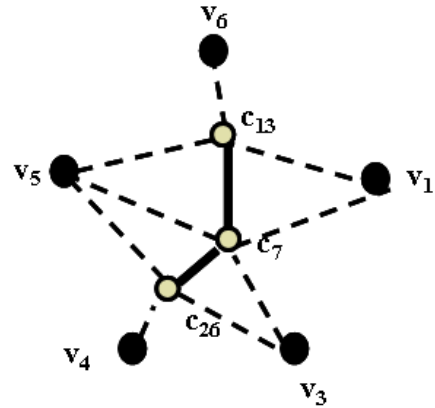
Рис. 17.2. Смешанный граф циклов.

Поставим в соответствие ободу 2-го слоя координатно-базисную систему векторов. Рассмотрев проекции убедимся, что соединения $(v_4,v_6)$ и $(v_4,v_1)$ не пересекаются (рис. 17.1). Построим смешанный граф циклов.

Выделим в смешанном графе циклов минимальный маршрут для соединения $(v_4,v_1)$. Минимальный маршрут состоит из отрезка $(c_{26},c_7)$ с включением мнимой вершины $v_{23}$.

Определим сопряженное неориентированное ребро как пересечение циклов $c_{26} \cap c_7$ с включением мнимой вершины $v_{23}$:

$c_{26} \cap c_7 = [(v_5,v_4)+(v_4,v_3)+(v_3,v_5)] \cap [(v_5,v_3)+(v_3,v_1)+(v_1,v_5)] = (v_3,v_5)$.

Рассмотрим следующую последовательность циклов:

$c_{26} = (v_5,v_4)+(v_4,v_3)+(v_3,v_5);$
$c_7 = (v_5,v_3)+(v_3,v_1)+(v_1,v_5).$

Вставим мнимую вершину в ориентированную часть ребра $(v_4,v_1)$:

$c_{26} = (v_5,v_4)+(v_4,v_3)+(v_3,v_{23})+(v_{23},v_5);$
$c_7 = (v_5,v_{23})+(v_{23},v_3)+(v_3,v_1)+(v_1,v_5).$

Упорядочим запись циклов, оставляя на первом месте ориентированные ребра с вершинами из кортежа вершин минимального маршрута и отделяя ребра с мнимыми вершинами:

$c_{26} = [(v_4,v_3)+(v_3,v_{23})]+[(v_{23},v_5)+(v_5,v_4)];$
$c_7 = [(v_1,v_5)+(v_5,v_{23})]+[(v_{23},v_3)+(v_3,v_1)].$

Вставим отрезки соединения $(v_4,v_1)$:

$c_{26} = [(v_4,v_3)+(v_3,v_{23})]+[(v_{23},v_4)+(v_4,v_{23})]+[(v_{23},v_5)+(v_5,v_4)]$;
$c_7 = [(v_1,v_5)+(v_5,v_{23})]+[(v_{23},v_1)+(v_1,v_{23})]+[(v_{23},v_3)+(v_3,v_1)]$.

Выделим замкнутую последовательность рёбер и сформируем новую систему циклов вместо циклов $c_7, c_{26}$:

$c_{78} = (v_4,v_3)+(v_3,v_{23})+(v_{23},v_4)$;
$c_{79} = (v_4,v_{23})+(v_{23},v_5)+(v_5,v_4)$;
$c_{80} = (v_1,v_5)+(v_5,v_{23})+(v_{23},v_1)$;
$c_{81} = (v_1,v_{23})+(v_{23},v_3)+(v_3,v_1)$.

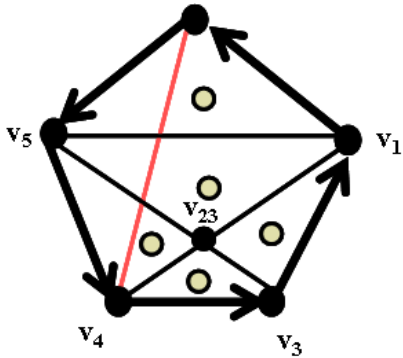 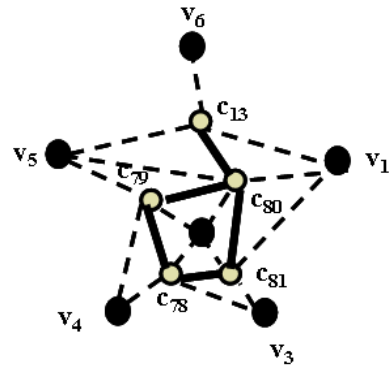

Рис. 17.3. Новое положение циклов.     Рис. 17.4. Смешанный граф циклов.

Построим смешанный граф циклов, для чего выделим минимальный маршрут в смешанном графе циклов. Минимальный маршрут состоит из отрезка $(c_{78},c_{79})$ с введением мнимой вершины $v_{24}$ и отрезка $(c_{79},c_{13})$ с введением мнимой вершины $v_{25}$. Запишем минимальный маршрут в виде кортежа вершин $<v_4,v_{24},v_{25},v_6>$.

Определим сопряжённое неориентированное ребро как пересечение циклов $c_{78} \cap c_{79}$ с включением мнимой вершины $v_{24}$:

$c_{79} \cap c_{80} = [(v_4,v_{23})+(v_{23},v_5)+(v_5,v_4)] \cap [(v_1,v_5)+(v_5,v_{23})+(v_{23},v_1)] = (v_5,v_{23})$.

Определим сопряжённое неориентированное ребро как пересечение циклов $c_{79} \cap c_{13}$ с включением мнимой вершины $v_{25}$:

$c_{80} \cap c_{13} = [(v_1,v_5)+(v_5,v_{23})+(v_{23},v_1)] \cap [(v_6,v_5)+(v_5,v_1)+(v_1,v_6)] = (v_1,v_5)$.

Рассмотрим следующую последовательность циклов:

$c_{79} = (v_4,v_{23})+(v_{23},v_5)+(v_5,v_4)$;
$c_{80} = (v_1,v_5)+(v_5,v_{23})+(v_{23},v_1)$;
$c_{13} = (v_6,v_5)+(v_5,v_1)+(v_1,v_6)$.

Введём мнимые вершины в соответствующие рёбра:

$c_{79} = (v_4,v_{23})+(v_{23},v_{24})+(v_{24},v_5)+(v_5,v_4)$;
$c_{80} = (v_1,v_{25})+(v_{25},v_5)+(v_5,v_{24})+(v_{24},v_{23})+(v_{23},v_1)$;
$c_{13} = (v_6,v_5)+(v_5,v_{25})+(v_{25},v_1)+(v_1,v_6)$.

Упорядочим запись циклов, оставляя на первом месте ориентированные рёбра с вершинами из кортежа вершин минимального маршрута и отделяя рёбра с мнимыми вершинами:

$c_{79} = [(v_4,v_{23})+(v_{23},v_{24})]+[(v_{24},v_5)+(v_5,v_4)]$;
$c_{80} = [(v_{24},v_{23})+(v_{23},v_1)+(v_1,v_{25})]+[(v_{25},v_5)+(v_5,v_{24})]$;
$c_{13} = ](v_6,v_5)+(v_5,v_{25})]+[(v_{25},v_1)+(v_1,v_6)]$.

Вставим отрезки соединения $(v_4,v_6)$:

$c_{79} = [(v_4,v_{23})+(v_{23},v_{24})]+[v_{24},v_4)+(v_4,v_{24})]+[(v_{24},v_5)+(v_5,v_4)]$;
$c_{80} = [(v_{24},v_{23})+(v_{23},v_1)+(v_1,v_{25})]+[(v_{25},v_{24})+(v_{24},v_{25})]+[(v_{25},v_5)+(v_5,v_{24})]$;
$c_{13} = [(v_6,v_5)+(v_5,v_{25})]+[v_{25},v_6)+(v_6,v_{25})]+[(v_{25},v_1)+(v_1,v_6)]$.

Выделим замкнутую последовательность ребер и сформируем новую систему циклов вместо циклов $c_{79}, c_{80}, c_{13}$:

$c_{82} = (v_4,v_{23})+(v_{23},v_{24})+v_{24},v_4)$;
$c_{83} = (v_4,v_{24})+(v_{24},v_5)+(v_5,v_4)$;
$c_{84} = (v_{24},v_{23})+(v_{23},v_1)+(v_1,v_{25})+(v_{25},v_{24})$;
$c_{85} = (v_{24},v_{25})+(v_{25},v_5)+(v_5,v_{24})$;
$c_{86} = (v_6,v_5)+(v_5,v_{25})+v_{25},v_6)$;
$c_{87} = (v_6,v_{25})+(v_{25},v_1)+(v_1,v_6)$.

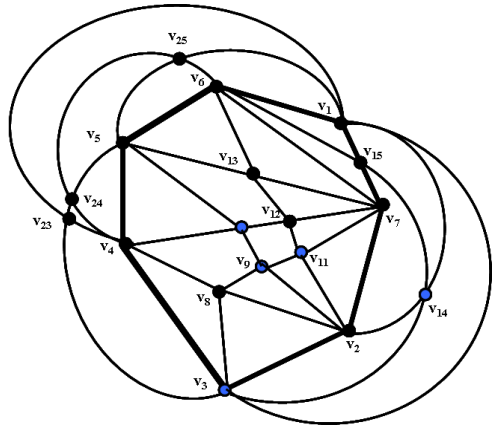

Рис. 17.5. Топологический рисунок графа с учетом пересечений.

Множество циклов, характеризующих топологический рисунок, имеет вид:

$c_{36} = (v_2,v_3)+(v_3,v_8)+(v_8.v_2)$;
$c_{39} = (v_4,v_8)+(v_8,v_3)+(v_3,v_4)$;
$c_{40} = (v_2,v_8)+(v_8,v_9)+(v_9,v_2)$;
$c_{43} = (v_4,v_{10})+(v_{10},v_9)+(v_9,v_8)+(v_8.v_4)$;
$c_{45} = (v_5,v_{10})+(v_{10},v_4)+(v_4,v_5)$;
$c_{46} = (v_2,v_9)+(v_9,v_{11})+(v_2,v_{11})$;
$c_{47} = (v_{11},v_2)+(v_{11},v_7)+(v_7,v_2)$;
$c_{48} = (v_{12},v_7)+(v_7,v_{11})+(v_{11},v_{12})$;
$c_{49} = (v_{12},v_{11})+(v_{11},v_9)+(v_9,v_{10})+(v_{10},v_{12})$;
$c_{50} = (v_{13},v_7)+(v_7,v_{12})+(v_{12},v_{13})$;
$c_{51} = (v_{13},v_{12})+(v_{12},v_{10})+(v_{10},v_5)+(v_5,v_{13})$;
$c_{52} = (v_6,v_7)+(v_7,v_{13})+(v_{13},v_6)$;
$c_{53} = (v_6,v_{13})+(v_{13},v_5)+(v_5,v_6)$;
$c_{54} = (v_3,v_2)+(v_2,v_{14})+(v_{14},v_3)$;
$c_{55} = (v_3,v_{14})+(v_{14},v_1)+(v_1,v_3)$;
$c_{56} = (v_{14},v_2)+(v_2,v_7)+(v_7,v_{15})+(v_{15},v_{14})$;
$c_{57} = (v_{14},v_{15})+(v_{15},v_1)+(v_1,v_{14})$;
$c_{58} = (v_6,v_1)+(v_1,v_{15})+(v_{15},v_6)$;
$c_{59} = (v_6,v_{15})+(v_{15},v_7)+(v_7,v_6)$.

$c_{78} = (v_4,v_3)+(v_3,v_{23})+(v_{23},v_4);$
$c_{81} = (v_1,v_{23})+(v_{23},v_3)+(v_3,v_1);$
$c_{82} = (v_4,v_{23})+(v_{23},v_{24})+v_{24},v_4);$
$c_{83} = (v_4,v_{24})+(v_{24},v_5)+(v_5,v_4);$
$c_{84} = (v_{24},v_{23})+(v_{23},v_1)+(v_1,v_{25})+(v_{25},v_{24});$
$c_{85} = (v_{24},v_{25})+(v_{25},v_5)+(v_5,v_{24});$
$c_{86} = (v_6,v_5)+(v_5,v_{25})+v_{25},v_6);$
$c_{87} = (v_6,v_{25})+(v_{25},v_1)+(v_1,v_6).$

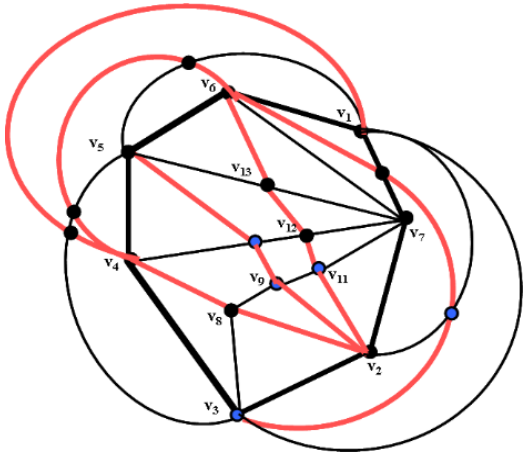

Рис. 17.6. Выделение соединений исключенных в процессе планаризации (красный цвет).

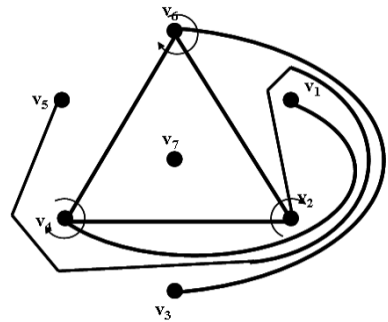

Рис. 17.7. Вращение вершин в плоском рисунке исключеных соединений.

Вращение вершин на рис 17.6 и рис. 17.7 для соединений, исключенных в процессе планаризации, совпадает.

Рассмотрим последовательность расположения мнимых вершин, расположенных на соединениях удаленных в процессе планаризации:

$<v_1.v_4> = <v_{23}>;$
$<v_2.v_4> = <v_8>;$
$<v_2.v_5> = <v_9,v_{10}>;$
$<v_2.v_6> = <v_{11},v_{12},v_{13}>;$
$<v_4.v_6> = <v_{24},v_{25}>;$
$<v_3.v_6> = <v_{14},v_{15}>.$

### Комментарии

В главе рассмотрены вопросы построения топологического рисунка непланарного графа с учетом проведения соединений в пространстве вне обода максимально плоского суграфа.

Таким образом показано, что может существовать множество видов изображения топологического рисунка непланарного графа разной толщины. Представлен вид топологического рисунка непланарного графа, построенный без учета части пространства вне обода максимально плоского суграфа. Другой вид топологического рисунка непланарного графа, характеризующий классическое понятие толщины графа, назовём *топологическим рисунком графа определенной толщины*.

# Глава 18. МАРШРУТЫ В СМЕШАННОМ ГРАФЕ ЦИКЛОВ

## 18.1. Полный граф К$_8$

Рассмотрим полный граф на 8 вершин К$_8$.

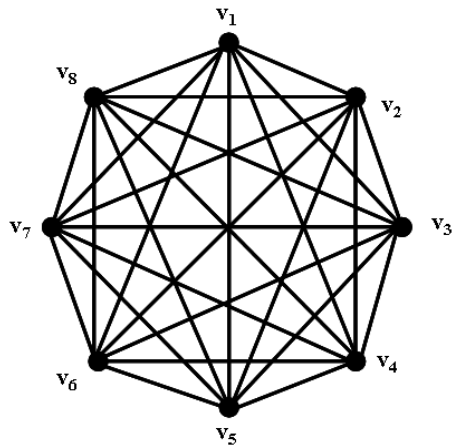

Рис. 18.1. Граф К$_8$.

Количество вершин графа = 8.

Количество ребер графа = 28.

Количество единичных циклов = 56.

Смежность графа К$_8$:

вершина $v_1$: $v_2$ $v_3$ $v_4$ $v_5$ $v_6$ $v_7$ $v_8$
вершина $v_2$: $v_1$ $v_3$ $v_4$ $v_5$ $v_6$ $v_7$ $v_8$
вершина $v_3$: $v_1$ $v_2$ $v_4$ $v_5$ $v_6$ $v_7$ $v_8$
вершина $v_4$: $v_1$ $v_2$ $v_3$ $v_5$ $v_6$ $v_7$ $v_8$
вершина $v_5$: $v_1$ $v_2$ $v_3$ $v_4$ $v_6$ $v_7$ $v_8$
вершина $v_6$: $v_1$ $v_2$ $v_3$ $v_4$ $v_5$ $v_7$ $v_8$
вершина $v_7$: $v_1$ $v_2$ $v_3$ $v_4$ $v_5$ $v_6$ $v_8$
вершина $v_8$: $v_1$ $v_2$ $v_3$ $v_4$ $v_5$ $v_6$ $v_7$

Инцидентность графа К$_8$:

ребро $e_1$: $(v_1,v_2)$ или $(v_2,v_1)$;
ребро $e_2$: $(v_1,v_3)$ или $(v_3,v_1)$;
ребро $e_3$: $(v_1,v_4)$ или $(v_4,v_1)$;
ребро $e_4$: $(v_1,v_5)$ или $(v_5,v_1)$;
ребро $e_5$: $(v_1,v_6)$ или $(v_6,v_1)$;
ребро $e_6$: $(v_1,v_7)$ или $(v_7,v_1)$;
ребро $e_7$: $(v_1,v_8)$ или $(v_8,v_1)$;
ребро $e_8$: $(v_2,v_3)$ или $(v_3,v_2)$;
ребро $e_9$: $(v_2,v_4)$ или $(v_4,v_2)$;
ребро $e_{10}$: $(v_2,v_5)$ или $(v_5,v_2)$;
ребро $e_{11}$: $(v_2,v_6)$ или $(v_6,v_2)$;
ребро $e_{12}$: $(v_2,v_7)$ или $(v_7,v_2)$;
ребро $e_{13}$: $(v_2,v_8)$ или $(v_8,v_2)$;
ребро $e_{14}$: $(v_3,v_4)$ или $(v_4,v_3)$;
ребро $e_{15}$: $(v_3,v_5)$ или $(v_5,v_3)$;
ребро $e_{16}$: $(v_3,v_6)$ или $(v_6,v_3)$;
ребро $e_{17}$: $(v_3,v_7)$ или $(v_7,v_3)$;
ребро $e_{18}$: $(v_3,v_8)$ или $(v_8,v_3)$;
ребро $e_{19}$: $(v_4,v_5)$ или $(v_5,v_4)$;
ребро $e_{20}$: $(v_4,v_6)$ или $(v_6,v_4)$;
ребро $e_{21}$: $(v_4,v_7)$ или $(v_7,v_4)$;
ребро $e_{22}$: $(v_4,v_8)$ или $(v_8,v_4)$;
ребро $e_{23}$: $(v_5,v_6)$ или $(v_6,v_5)$;
ребро $e_{24}$: $(v_5,v_7)$ или $(v_7,v_5)$;
ребро $e_{25}$: $(v_5,v_8)$ или $(v_8,v_5)$;
ребро $e_{26}$: $(v_6,v_7)$ или $(v_7,v_6)$;
ребро $e_{27}$: $(v_6,v_8)$ или $(v_8,v_6)$;
ребро $e_{28}$: $(v_7,v_8)$ или $(v_8,v_7)$.

Множество изометрических циклов графа:

цикл   $c_1 = \{e_1,e_2,e_8\} \leftrightarrow \{v_1,v_2,v_3\}$;
цикл   $c_2 = \{e_1,e_3,e_9\} \leftrightarrow \{v_1,v_2,v_4\}$;
цикл   $c_3 = \{e_1,e_4,e_{10}\} \leftrightarrow \{v_1,v_2,v_5\}$;
цикл   $c_4 = \{e_1,e_5,e_{11}\} \leftrightarrow \{v_1,v_2,v_6\}$;
цикл   $c_5 = \{e_1.e_6.e_{12}\} \leftrightarrow \{v_1,v_2,v_7\}$;
цикл   $c_6 = \{e_1.e_7.e_{13}\} \leftrightarrow \{v_1,v_2,v_8\}$;
цикл   $c_7 = \{e_2.e_3.e_{14}\} \leftrightarrow \{v_1,v_3,v_4\}$;
цикл   $c_8 = \{e_2.e_4.e_{15}\} \leftrightarrow \{v_1,v_3,v_5\}$;
цикл   $c_9 = \{e_2.e_5.e_{16}\} \leftrightarrow \{v_1,v_3,v_6\}$;
цикл  $c_{10} = \{e_2,e_6,e_{17}\} \leftrightarrow \{v_1,v_3,v_7\}$;
цикл  $c_{11} = \{e_2,e_7,e_{18}\} \leftrightarrow \{v_1.v_3.v_8\}$;
цикл  $c_{12} = \{e_3,e_4,e_{19}\} \leftrightarrow \{v_1,v_4,v_5\}$;
цикл  $c_{13} = \{e_3,e_5,e_{20}\} \leftrightarrow \{v_1,v_4,v_6\}$;
цикл  $c_{14} = \{e_3,e_6,e_{21}\} \leftrightarrow \{v_1,v_4,v_7\}$;
цикл  $c_{15} = \{e_3,e_7,e_{22}\} \leftrightarrow \{v_1,v_4,v_8\}$;
цикл  $c_{16} = \{e_4,e_5,e_{23}\} \leftrightarrow \{v_1,v_5,v_6\}$;
цикл  $c_{17} = \{e_4,e_6,e_{24}\} \leftrightarrow \{v_1,v_5,v_7\}$;
цикл  $c_{18} = \{e_4,e_7,e_{25}\} \leftrightarrow \{v_1,v_5,v_8\}$;
цикл  $c_{19} = \{e_5,e_6,e_{26}\} \leftrightarrow \{v_1,v_6,v_7\}$;
цикл  $c_{20} = \{e_5,e_7,e_{27}\} \leftrightarrow \{v_1,v_6,v_8\}$;
цикл  $c_{21} = \{e_6,e_7,e_{28}\} \leftrightarrow \{v_1,v_7,v_8\}$;
цикл  $c_{22} = \{e_8,e_9,e_{14}\} \leftrightarrow \{v_2,v_3,v_4\}$;
цикл  $c_{23} = \{e_8,e_{10},e_{15}\} \leftrightarrow \{v_2,v_3,v_5\}$;
цикл  $c_{24} = \{e_8,e_{11},e_{16}\} \leftrightarrow \{v_2,v_3,v_6\}$;
цикл  $c_{25} = \{e_8,e_{12},e_{17}\} \leftrightarrow \{v_2,v_3,v_7\}$;
цикл  $c_{26} = \{e_8,e_{13},e_{18}\} \leftrightarrow \{v_2,v_3,v_8\}$;
цикл  $c_{27} = \{e_9,e_{10},e_{19}\} \leftrightarrow \{v_2,v_4,v_5\}$;
цикл  $c_{28} = \{e_9,e_{11},e_{20}\} \leftrightarrow \{v_2,v_4,v_6\}$;
цикл  $c_{29} = \{e_9,e_{12},e_{21}\} \leftrightarrow \{v_2,v_4,v_7\}$;
цикл  $c_{30} = \{e_9,e_{13},e_{22}\} \leftrightarrow \{v_2,v_4,v_8\}$;
цикл  $c_{31} = \{e_{10},e_{11},e_{23}\} \leftrightarrow \{v_2,v_5,v_6\}$;
цикл  $c_{32} = \{e_{10},e_{12},e_{24}\} \leftrightarrow \{v_2,v_5,v_7\}$;
цикл  $c_{33} = \{e_{10},e_{13},e_{25}\} \leftrightarrow \{v_2,v_5,v_8\}$;
цикл  $c_{34} = \{e_{11},e_{12},e_{26}\} \leftrightarrow \{v_2,v_6,v_7\}$;
цикл  $c_{35} = \{e_{11},e_{13},e_{27}\} \leftrightarrow \{v_2,v_6,v_8\}$;
цикл  $c_{36} = \{e_{12},e_{13},e_{28}\} \leftrightarrow \{v_2,v_7,v_8\}$;
цикл  $c_{37} = \{e_{14},e_{15},e_{19}\} \leftrightarrow \{v_3,v_4,v_5\}$;
цикл  $c_{38} = \{e_{14},e_{16},e_{20}\} \leftrightarrow \{v_3,v_4,v_6\}$;
цикл  $c_{39} = \{e_{14},e_{17},e_{21}\} \leftrightarrow \{v_3,v_4,v_7\}$;
цикл  $c_{40} = \{e_{14},e_{18},e_{22}\} \leftrightarrow \{v_3,v_4,v_8\}$;
цикл  $c_{41} = \{e_{15},e_{16},e_{23}\} \leftrightarrow \{v_3,v_5,v_6\}$;
цикл  $c_{42} = \{e_{15},e_{17},e_{24}\} \leftrightarrow \{v_3,v_5,v_7\}$;
цикл  $c_{43} = \{e_{15},e_{18},e_{25}\} \leftrightarrow \{v_3,v_5,v_8\}$;
цикл  $c_{44} = \{e_{16},e_{17},e_{26}\} \leftrightarrow \{v_3,v_6,v_7\}$;
цикл  $c_{45} = \{e_{16},e_{18},e_{27}\} \leftrightarrow \{v_3,v_6,v_8\}$;
цикл  $c_{46} = \{e_{17},e_{18},e_{28}\} \leftrightarrow \{v_3,v_7,v_8\}$;
цикл  $c_{47} = \{e_{19},e_{20},e_{23}\} \leftrightarrow \{v_4,v_5,v_6\}$;
цикл  $c_{48} = \{e_{19},e_{21},e_{24}\} \leftrightarrow \{v_4,v_5,v_7\}$;
цикл  $c_{49} = \{e_{19},e_2,e_{25}\} \leftrightarrow \{v_4,v_5,v_8\}$;
цикл  $c_{50} = \{e_{20},e_{21},e_{26}\} \leftrightarrow \{v_4,v_6,v_7\}$;

цикл  $c_{51} = \{e_{20},e_{22},e_{27}\} \leftrightarrow \{v_4,v_6,v_8\}$;
цикл  $c_{52} = \{e_{21},e_{22},e_{28}\} \leftrightarrow \{v_4,v_7,v_8\}$;
цикл  $c_{53} = \{e_{23},e_{24},e_{26}\} \leftrightarrow \{v_5,v_6,v_7\}$;
цикл  $c_{54} = \{e_{23},e_{25},e_{27}\} \leftrightarrow \{v_5,v_6,v_8\}$;
цикл  $c_{55} = \{e_{24},e_{25},e_{28}\} \leftrightarrow \{v_5,v_7,v_8\}$;
цикл  $c_{56} = \{e_{26},e_{27},e_{28}\} \leftrightarrow \{v_6,v_7,v_8\}$.

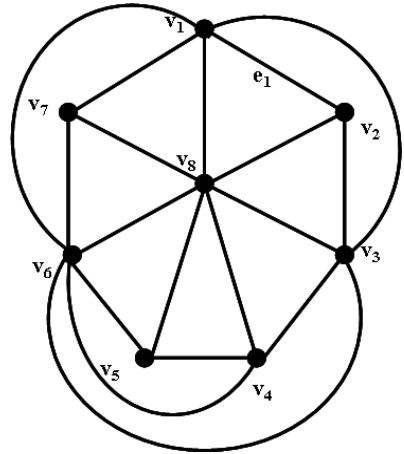

Рис. 18.2. Максимально плоский суграф

Множество циклов характеризующих максимально плоский суграф:

цикл  $c_1 = \{e_1,e_2,e_8\} \leftrightarrow \{v_1,v_2,v_3\}$;
цикл  $c_6 = \{e_1.e_7.e_{13}\} \leftrightarrow \{v_1,v_2,v_8\}$;
цикл  $c_{19} = \{e_5,e_6,e_{26}\} \leftrightarrow \{v_1,v_6,v_7\}$;
цикл  $c_{21} = \{e_6,e_7,e_{28}\} \leftrightarrow \{v_1,v_7,v_8\}$;
цикл  $c_{26} = \{e_8,e_{13},e_{18}\} \leftrightarrow \{v_2,v_3,v_8\}$;
цикл  $c_{38} = \{e_{14},e_{16},e_{20}\} \leftrightarrow \{v_3,v_4,v_6\}$;
цикл  $c_{40} = \{e_{14},e_{18},e_{22}\} \leftrightarrow \{v_3,v_4,v_8\}$;
цикл  $c_{47} = \{e_{19},e_{20},e_{23}\} \leftrightarrow \{v_4,v_5,v_6\}$;
цикл  $c_{49} = \{e_{19},e_2,e_{25}\} \leftrightarrow \{v_4,v_5,v_8\}$;
цикл  $c_{54} = \{e_{23},e_{25},e_{27}\} \leftrightarrow \{v_5,v_6,v_8\}$;
цикл  $c_{56} = \{e_{26},e_{27},e_{28}\} \leftrightarrow \{v_6,v_7,v_8\}$.

Зададим вращение циклов в векторном виде:

$c_1 = (v_1,v_2)+(v_2,v_3)+(v_3,v_1)$;
$c_6 = (v_2,v_1)+(v_1,v_8)+(v_8,v_2)$;
$c_{21} = (v_8,v_1)+(v_1,v_7)+(v_7,v_8)$;
$c_{19} = (v_7,v_1)+(v_1,v_6)+(v_6,v_7)$;
$c_{56} = (v_7,v_6)+(v_6,v_8)+(v_8,v_7)$;
$c_{54} = (v_8,v_6)+(v_6,v_5)+(v_5,v_8)$;
$c_{49} = (v_8,v_5)+(v_5,v_4)+(v_4,v_8)$;
$c_{47} = (v_4,v_5)+(v_5,v_6)+(v_6,v_4)$;
$c_{40} = (v_8,v_4)+(v_4,v_3)+(v_3,v_8)$;
$c_{26} = (v_8,v_3)+(v_3,v_2)+(v_2,v_8)$;
$c_{38} = (v_4,v_6)+(v_6,v_3)+(v_3,v_4)$.

Ребра исключенные из процесса планаризации:

ребро  $e_3$: $(v_1,v_4)$ или $(v_4,v_1)$;        ребро  $e_4$: $(v_1,v_5)$ или $(v_5,v_1)$;
ребро  $e_9$: $(v_2,v_4)$ или $(v_4,v_2)$;        ребро  $e_{10}$: $(v_2,v_5)$ или $(v_5,v_2)$;
ребро  $e_{11}$: $(v_2,v_6)$ или $(v_6,v_2)$;      ребро  $e_{12}$: $(v_2,v_7)$ или $(v_7,v_2)$;

ребро $e_{15}$: $(v_3,v_5)$ или $(v_5,v_3)$; ребро $e_{17}$: $(v_3,v_7)$ или $(v_7,v_3)$;
ребро $e_{21}$: $(v_4,v_7)$ или $(v_7,v_4)$; ребро $e_{24}$: $(v_5,v_7)$ или $(v_7,v_5)$;

Построим в максимально плоском суграфе гамильтонов цикл. Поставим в соответствие гамильтонову циклу координатно-базисную систему векторов. Спроектируем множество соединений исключенных в процессе планаризации (красный цвет) на координатно-базисную систему векторов (рис. 18.3). Выделим множество непересекающихся соединений (красный цвет) по проекциям векторов на КБС (рис. 18.4).

Выделим подмножество циклов, расположенных внутри гамильтонова цикла:

$c_6 = (v_2,v_1)+(v_1,v_8)+(v_8,v_2)$;
$c_{21} = (v_8,v_1)+(v_1,v_7)+(v_7,v_8)$;
$c_{56} = (v_7,v_6)+(v_6,v_8)+(v_8,v_7)$;
$c_{54} = (v_8,v_6)+(v_6,v_5)+(v_5,v_8)$;
$c_{49} = (v_8,v_5)+(v_5,v_4)+(v_4,v_8)$;
$c_{40} = (v_8,v_4)+(v_4,v_3)+(v_3,v_8)$.

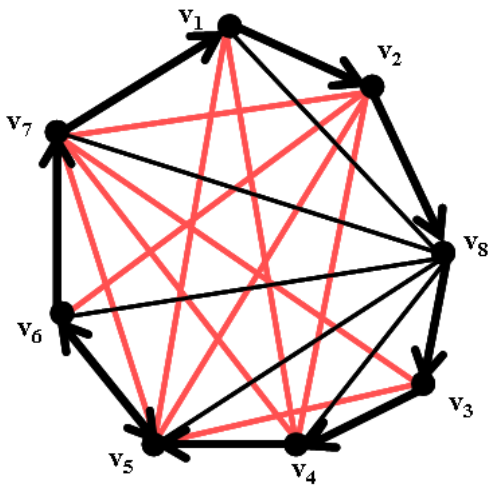 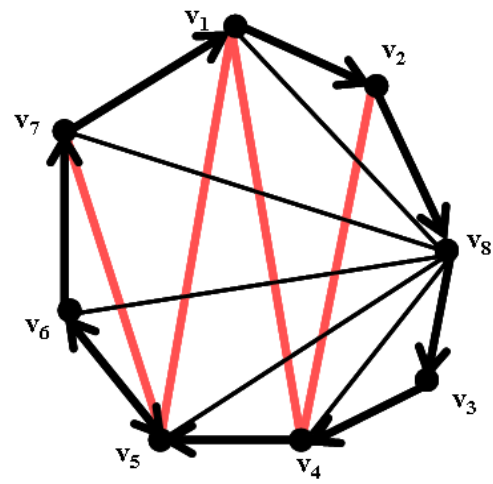

Рис. 18.3. Соединения и КБС.    Рис. 18.4. Выделенные соединения.

### 18.2. Соединение $e_{24} = (v_5,v_7)$

Построим смешанный граф циклов для множества циклов, образующих внутренность гамильтонова цикла. На рис. 18.5 представлен процесс выделения минимального маршрута, соединяющего циклы $c_{54}$ и $c_{56}$, для проведения соединения $(v_5,v_7)$. Минимальный маршрут состоит из отрезка $(c_{54},c_{56})$ с введением мнимой вершины $v_9$. Запишем минимальный маршрут в виде кортежа вершин $<v_5,v_9,v_7>$.

Определим сопряженное неориентированное ребро как пересечение циклов $c_{54} \cap c_{56}$ с включением мнимой вершины $v_9$:

$c_{54} \cap c_{56} = [(v_8,v_6)+(v_6,v_5)+(v_5,v_8)] \cap [(v_7,v_6)+(v_6,v_8)+(v_8,v_7)] = (v_6,v_8)$.

Рассмотрим следующую последовательность циклов:

$c_{54} = (v_8,v_6)+(v_6,v_5)+(v_5,v_8)$;
$c_{56} = (v_7,v_6)+(v_6,v_8)+(v_8,v_7)$.

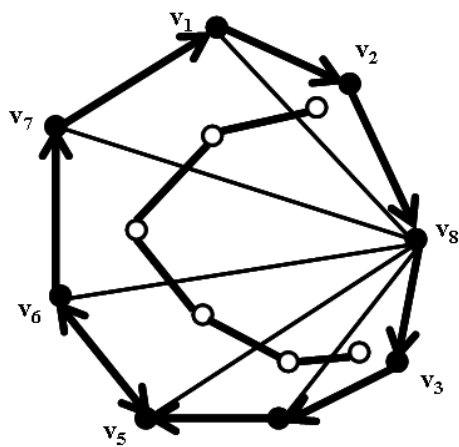
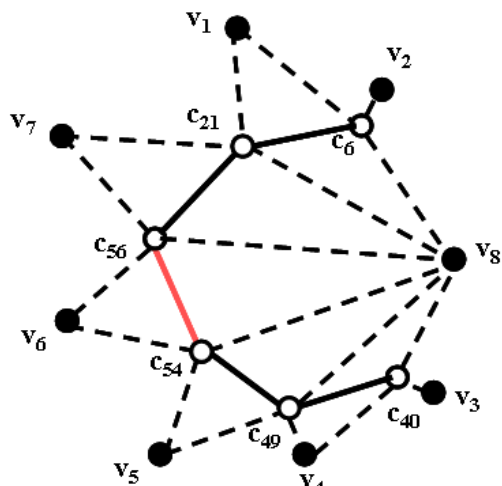

Рис. 18.5. Граф циклов.    Рис. 18.6. Смешанный граф циклов.

Введём мнимую вершину в соответствующее ребро:

$c_{54} = (v_8,v_9)+(v_9,v_6)+(v_6,v_5)+(v_5,v_8)$;
$c_{56} = (v_7,v_6)+(v_6,v_9)+(v_9,v_8)+(v_8,v_7)$.

Упорядочим запись циклов, оставляя на первом месте ориентированные рёбра с вершинами из кортежа вершин минимального маршрута и отделяя рёбра с мнимыми вершинами:

$c_{54} = [(v_5,v_8)\ (v_8,v_9)]+[(v_9,v_6)+(v_6,v_5)]$;
$c_{56} = [(v_7,v_6)+(v_6,v_9)]+[(v_9,v_8)+(v_8,v_7)]$.

Вставим отрезки соединения $(v_5,v_7)$:

$c_{54} = [(v_5,v_8)\ (v_8,v_9)]+[v_9,v_5)+(v_5,v_9)]+[(v_9,v_6)+(v_6,v_5)]$;
$c_{56} = [(v_7,v_6)+(v_6,v_9)]+[v_9,v_7)+(v_7,v_9)]+[(v_9,v_8)+(v_8,v_7)]$.

Выделим замкнутую последовательность рёбер и сформируем новую систему циклов вместо циклов $c_{54}, c_{56}$:

$c_{57} = (v_5,v_8)+(v_8,v_9)+(v_9,v_5)$;
$c_{58} = (v_5,v_9)+(v_9,v_6)+(v_6,v_5)$;
$c_{59} = (v_7,v_6)+(v_6,v_9)+(v_9,v_7)$;
$c_{60} = (v_7,v_9)+(v_9,v_8)+(v_8,v_7)$.

В результате образуется множество циклов, характеризующих топологический рисунок:

$c_6 = (v_2,v_1)+(v_1,v_8)+(v_8,v_2)$;
$c_{21} = (v_8,v_1)+(v_1,v_7)+(v_7,v_8)$;
$c_{49} = (v_8,v_5)+(v_5,v_4)+(v_4,v_8)$;
$c_{40} = (v_8,v_4)+(v_4,v_3)+(v_3,v_8)$;
$c_{57} = (v_5,v_8)+(v_8,v_9)+(v_9,v_5)$;
$c_{58} = (v_5,v_9)+(v_9,v_6)+(v_6,v_5)$;
$c_{59} = (v_7,v_6)+(v_6,v_9)+(v_9,v_7)$;
$c_{60} = (v_7,v_9)+(v_9,v_8)+(v_8,v_7)$;
обод $= (v_1,v_2)+(v_2,v_8)+(v_8,v_3)+(v_3,v_4)+(v_4,v_5)+(v_5,v_6)+(v_6,v_7)+(v_7,v_1)$.

В результате проведения соединения $(v_5,v_7)$ образуется новая система циклов и новый граф циклов (рис. 18.7). С учетом соединения вершин графа и вершин графа циклов строится смешанный граф циклов (рис. 18.8). При построении смешанного графа циклов учтем тот факт, что вновь введенное соединение (красный цвет) пересекает некоторые ребра графа циклов (рис. 18.7). Впредь будем исключать такие ребра из построения смешанного графа циклов (рис. 18.9).

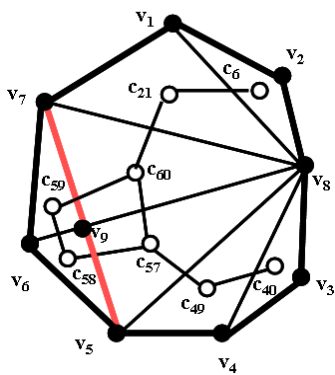
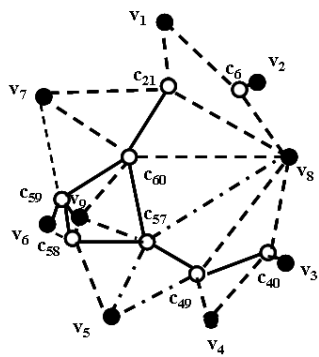
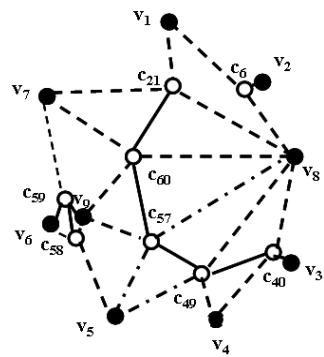

Рис. 18.7. Граф циклов        Рис. 18.8. Смешанный граф циклов.        Рис. 18.9. Смешанный граф циклов с учетом введения мнимой вершины.

### 18.3. Соединение $e_4 = (v_5, v_1)$

Строим минимальный маршрут в смешанном графе циклов для проведения соединения $(v_5, v_1)$. Запишем минимальный маршрут в виде кортежа $<v_5, c_{57}, c_{60}, c_{21}, v_1>$ с мнимыми вершинами $\{v_{10}, v_{11}\}$.

Определим сопряженное неориентированное ребро как пересечение циклов $c_{57} \cap c_{60}$ с включением мнимой вершины $v_{10}$:

$c_{57} \cap c_{60} = [(v_5,v_8)+(v_8,v_9)+(v_9,v_5)] \cap [(v_7,v_9)+(v_9,v_8)+(v_8,v_7)] = (v_8,v_9)$.

Определим сопряженное неориентированное ребро как пересечение циклов $c_{60} \cap c_{21}$ с включением мнимой вершины $v_{11}$:

$c_{60} \cap c_{21} = [(v_7,v_9)+(v_9,v_8)+(v_8,v_7)] \cap [(v_8,v_1)+(v_1,v_7)+(v_7,v_8)] = (v_7,v_8)$.

Рассмотрим следующую последовательность циклов:

$c_{57} = (v_5,v_8)+(v_8,v_9)+(v_9,v_5)$;
$c_{60} = (v_7,v_9)+(v_9,v_8)+(v_8,v_7)$;
$c_{21} = (v_8,v_1)+(v_1,v_7)+(v_7,v_8)$.

Введём мнимые вершины в соответствующие ребра:

$c_{57} = (v_5,v_8)+(v_8,v_{10})+(v_{10},v_9)+(v_9,v_5)$;
$c_{60} = (v_7,v_9)+(v_9,v_{10})+(v_{10},v_8)+(v_8,v_{11})+(v_{11},v_7)$;
$c_{21} = (v_8,v_1)+(v_1,v_7)+(v_7,v_{11})+(v_{11},v_8)$.

Упорядочим запись циклов, оставляя на первом месте ориентированные ребра с вершинами из кортежа вершин минимального маршрута и отделяя ребра с мнимыми вершинами:

$c_{57} = [(v_5,v_8)+(v_8,v_{10})]+[(v_{10},v_9)+(v_9,v_5)];$
$c_{60} = [(v_{10},v_8)+(v_8,v_{11})]+[(v_{11},v_7)+(v_7,v_9)+(v_9,v_{10})];$
$c_{21} = [(v_1,v_7)+(v_7,v_{11})]+[(v_{11},v_8)+(v_8,v_1)].$

Вставим отрезки соединения $(v_5,v_1)$:

$c_{57} = [(v_5,v_8)+(v_8,v_{10})]+[(v_{10},v_5)+(v_5,v_{10})]+[(v_{10},v_9)+(v_9,v_5)];$
$c_{60} = [(v_{10},v_8)+(v_8,v_{11})]+[(v_{11},v_{10})+(v_{10},v_{11})]+[(v_{11},v_7)+(v_7,v_9)+(v_9,v_{10})];$
$c_{21} = [(v_1,v_7)+(v_7,v_{11})]+[(v_{11},v_1)+(v_{11},v_1)]+[(v_{11},v_8)+(v_8,v_1)].$

Выделим замкнутую последовательность рёбер и сформируем новую систему циклов вместо циклов $\{c_{21},c_{60},c_{57}\}$:

$c_{61} = (v_5,v_8)+(v_8,v_{10})+(v_{10},v_5);$
$c_{62} = (v_5,v_{10})+(v_{10},v_9)+(v_9,v_5);$
$c_{63} = (v_{10},v_8)+(v_8,v_{11})+(v_{11},v_{10});$
$c_{64} = (v_{10},v_{11})+(v_{11},v_7)+(v_7,v_9)+(v_9,v_{10});$
$c_{65} = (v_1,v_7)+(v_7,v_{11})+(v_{11},v_1);$
$c_{66} = (v_{11},v_1)+(v_{11},v_8)+(v_8,v_1).$

Проведение соединения $(v_5,v_1)$ порождает новую систему циклов. В свою очередь, на основе новой системы циклов строится смешанный граф циклов (рис. 18.11).

В результате образуется множество циклов, характеризующих топологический рисунок:

$c_6 = (v_2,v_1)+(v_1,v_8)+(v_8,v_2);$
$c_{49} = (v_8,v_5)+(v_5,v_4)+(v_4,v_8);$
$c_{40} = (v_8,v_4)+(v_4,v_3)+(v_3,v_8);$
$c_{58} = (v_5,v_9)+(v_9,v_6)+(v_6,v_5);$
$c_{59} = (v_7,v_6)+D(v_6,v_9)+(v_9,v_7);$
$c_{61} = (v_5,v_8)+(v_8,v_{10})+(v_{10},v_5);$
$c_{62} = (v_5,v_{10})+(v_{10},v_9)+(v_9,v_5);$
$c_{63} = (v_{10},v_8)+(v_8,v_{11})+(v_{11},v_{10});$
$c_{64} = (v_{10},v_{11})+(v_{11},v_7)+(v_7,v_9)+(v_9,v_{10});$
$c_{65} = (v_1,v_7)+(v_7,v_{11})+(v_{11},v_1);$
$c_{66} = (v_{11},v_1)+(v_{11},v_8)+(v_8,v_1);$
обод $= (v_1,v_2)+(v_2,v_8)+(v_8,v_3)+(v_3,v_4)+(v_4,v_5)+(v_5,v_6)+(v_6,v_7)+(v_7,v_1).$

### 18.4. Соединение $e_3 = (v_4,v_1)$

Выделим минимальный маршрут для проведения соединения $(v_4,v_1)$. Запишем минимальный маршрут в виде кортежа $<v_4,c_{49},c_{61},c_{63},c_{66},v_1>$ с учётом введения мнимых вершин $\{v_{12},v_{13},v_{14}\}$.

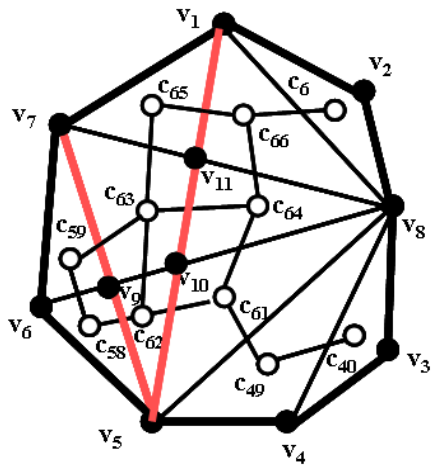 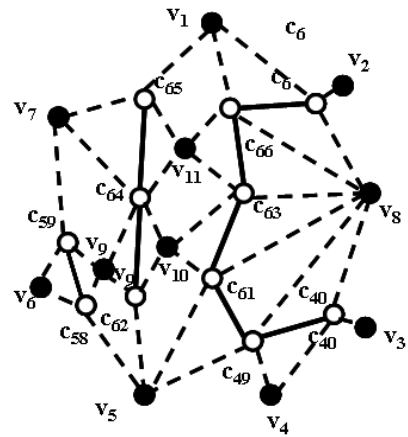

Рис. 18.10. Проведение соединения $(v_5, v_1)$.     Рис. 18.11. Смешанный граф циклов.

Определим сопряженное неориентированное ребро как пересечение циклов $c_{49} \cap c_{61}$ с включением мнимой вершины $v_{12}$:

$c_{49} \cap c_{61} = [(v_8,v_5)+(v_5,v_4)+(v_4,v_8)] \cap [(v_5,v_8)+(v_8,v_{10})+(v_{10},v_5)] = (v_8,v_5)$.

Определим сопряженное неориентированное ребро как пересечение циклов $c_{61} \cap c_{64}$ с включением мнимой вершины $v_{13}$:

$c_{61} \cap c_{63} = [(v_5,v_8)+(v_8,v_{10})+(v_{10},v_5)] \cap [(v_{10},v_8)+(v_8,v_{11})+(v_{11},v_{10})] = (v_{10},v_8)$.

Определим сопряженное неориентированное ребро как пересечение циклов $c_{63} \cap c_{66}$ с включением мнимой вершины $v_{14}$:

$c_{63} \cap c_{66} = [(v_{10},v_8)+(v_8,v_{11})+(v_{11},v_{10})] \cap [(v_1,v_{11})+(v_{11},v_8)+(v_8,v_1)] = (v_{11},v_8)$.

Рассмотрим следующую последовательность циклов:

$c_{49} = (v_8,v_5)+(v_5,v_4)+(v_4,v_8)$;
$c_{61} = (v_5,v_8)+(v_8,v_{10})+(v_{10},v_5)$;
$c_{63} = (v_{10},v_8)+(v_8,v_{11})+(v_{11},v_{10})$;
$c_{66} = (v_1,v_{11})+(v_{11},v_8)+(v_8,v_1)$.

Введём мнимые вершины в соответствующие рёбра:

$c_{49} = (v_8,v_{12})+(v_{12},v_5)+(v_5,v_4)+(v_4,v_8)$;
$c_{61} = (v_5,v_{12})+(v_{12},v_8)+(v_8,v_{13})+(v_{13},v_{10})+(v_{10},v_5)$;
$c_{63} = (v_{10},v_{13})+(v_{13},v_8)+(v_8,v_{14})+(v_{14},v_{11})+(v_{11},v_{10})$;
$c_{66} = (v_1,v_{11})+(v_{11},v_{14})+(v_{14},v_8)+(v_8,v_1)$.

Упорядочим запись циклов, оставляя на первом месте ориентированные рёбра с вершинами из кортежа вершин минимального маршрута и отделяя рёбра с мнимыми вершинами:

$c_{49} = [(v_4,v_8)+(v_8,v_{12})]+[(v_{12},v_5)+(v_5,v_4)]$;
$c_{61} = [(v_{12},v_8)+(v_8,v_{13})]+[(v_{13},v_{10})+(v_{10},v_5)+(v_5,v_{12})]$;
$c_{63} = [(v_{13},v_8)+(v_8,v_{14})]+[(v_{14},v_{11})+(v_{11},v_{10})+(v_{10},v_{13})]$;
$c_{66} = [(v_1,v_{11})+(v_{11},v_{14})]+[(v_{14},v_8)+(v_8,v_1)]$.

Вставим отрезки соединения $(v_4,v_1)$:

$c_{49} = [(v_4,v_8)+(v_8,v_{12})]+[(v_{12},v_4)+(v_4,v_{12})]+[(v_{12},v_5)+(v_5,v_4)]$;
$c_{61} = [(v_{12},v_8)+(v_8,v_{13})]+[(v_{13},v_{12})+(v_{12},v_{13})]+[(v_{13},v_{10})+(v_{10},v_5)+(v_5,v_{12})]$;

$c_{63} = [(v_{13},v_8)+(v_8,v_{14})]+[(v_{14},v_{13})+(v_{13},v_{14})]+[(v_{14},v_{11})+(v_{11},v_{10})+(v_{10},v_{13})]$;
$c_{66} = [(v_1,v_{11})+(v_{11},v_{14})]+[(v_{14},v_1)+(v_1,v_{14})]+[(v_{14},v_8)+(v_8,v_1)]$.

Выделим замкнутую последовательность ребер и сформируем новую систему циклов вместо циклов $\{c_{49},c_{61},c_{63},c_{66}\}$:

$c_{67} = (v_4,v_8)+(v_8,v_{12})+(v_{12},v_4)$;
$c_{68} = (v_4,v_{12})+(v_{12},v_5)+(v_5,v_4)$;
$c_{69} = (v_{12},v_8)+(v_8,v_{13})+(v_{13},v_{12})$;
$c_{70} = (v_{12},v_{13})+(v_{13},v_{10})+(v_{10},v_5)+(v_5,v_{12})$;
$c_{71} = (v_{13},v_8)+(v_8,v_{14})+(v_{14},v_{13})$;
$c_{72} = (v_{13},v_{14})+(v_{14},v_{11})+(v_{11},v_{10})+(v_{10},v_{13})$;
$c_{73} = (v_1,v_{11})+(v_{11},v_{14})+(v_{14},v_1)$;
$c_{74} = (v_1,v_{14})+(v_{14},v_8)+(v_8,v_1)$.

Проведение соединения $(v_4,v_1)$ порождает новую систему циклов. В свою очередь, на основе новой системы циклов строится смешанный граф циклов (рис. 18.13).

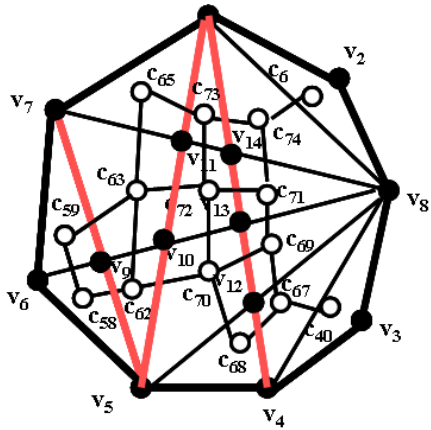 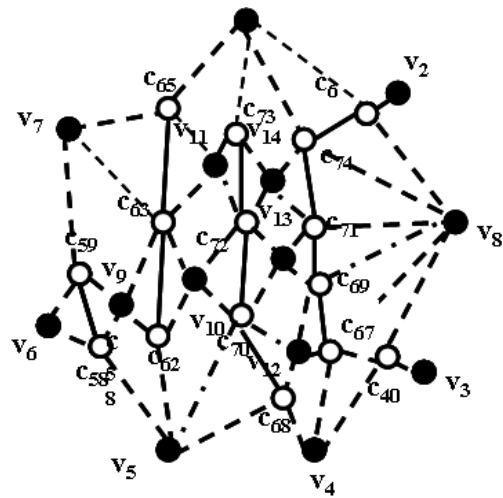

Рис. 18.12. Проведение соединения $(v_4,v_1)$.   рис. 18.13. Смешанный граф циклов.

В результате образуется множество циклов, характеризующих топологический рисунок:

$c_6 = (v_2,v_1)+(v_1,v_8)+(v_8,v_2)$;
$c_{40} = (v_8,v_4)+(v_4,v_3)+(v_3,v_8)$;
$c_{58} = (v_5,v_9)+(v_9,v_6)+(v_6,v_5)$;
$c_{59} = (v_7,v_6)+(v_6,v_9)+(v_9,v_7)$;
$c_{62} = (v_5,v_{10})+(v_{10},v_9)+(v_9,v_5)$;
$c_{64} = (v_{10},v_{11})+(v_{11},v_7)+(v_7,v_9)+(v_9,v_{10})$;
$c_{65} = (v_1,v_7)+(v_7,v_{11})+(v_{11},v_1)$;
$c_{67} = (v_4,v_8)+(v_8,v_{12})+(v_{12},v_4)$;
$c_{68} = (v_4,v_{12})+(v_{12},v_5)+(v_5,v_4)$;
$c_{69} = (v_{12},v_8)+(v_8,v_{13})+(v_{13},v_{12})$;
$c_{70} = (v_{12},v_{13})+(v_{13},v_{10})+(v_{10},v_5)+(v_5,v_{12})$;
$c_{71} = (v_{13},v_8)+(v_8,v_{14})+(v_{14},v_{13})$;
$c_{72} = (v_{13},v_{14})+(v_{14},v_{11})+(v_{11},v_{10})+(v_{10},v_{13})$;
$c_{73} = (v_1,v_{11})+(v_{11},v_{14})+(v_{14},v_1)$;
$c_{74} = (v_1,v_{14})+(v_{14},v_8)+(v_8,v_1)$;

обод = $(v_1,v_2)+(v_2,v_8)+(v_8,v_3)+(v_3,v_4)+(v_4,v_5)+(v_5,v_6)+(v_6,v_7)+(v_7,v_1)$.

### 18.5. Соединение $e_9 = (v_4, v_2)$

Выделим минимальный маршрут для проведения соединения $(v_4, v_2)$. Запишем минимальный маршрут в виде кортежа $<v_4, c_{67}, c_{69}, c_{71}, c_{74}, c_6, v_2>$ с учетом введения мнимых вершин $\{v_{15}, v_{16}, v_{17}, v_{18}\}$.

Определим сопряженное неориентированное ребро как пересечение циклов $c_{67} \cap c_{69}$ с включением мнимой вершины $v_{15}$:

$c_{69} \cap c_{71} = [(v_4,v_8)+(v_8,v_{12})+(v_{12},v_4)] \cap [(v_{12},v_8)+(v_8,v_{13})+(v_{13},v_{12})] = (v_8,v_{12})$.

Определим сопряженное неориентированное ребро как пересечение циклов $c_{69} \cap c_{71}$ с включением мнимой вершины $v_{16}$:

$c_{69} \cap c_{71} = [(v_{12},v_8)+(v_8,v_{13})+(v_{13},v_{12})] \cap [(v_{13},v_8)+(v_8,v_{14})+(v_{14},v_{13})] = (v_{13},v_8)$.

Определим сопряженное неориентированное ребро как пересечение циклов $c_{71} \cap c_{74}$ с включением мнимой вершины $v_{17}$:

$c_{71} \cap c_{74} = [(v_{13},v_8)+(v_8,v_{14})+(v_{14},v_{13})] \cap [(v_1,v_{14})+(v_{14},v_8)+(v_8,v_1)] = (v_{14},v_8)$.

Определим сопряженное неориентированное ребро как пересечение циклов $c_{74} \cap c_6$ с включением мнимой вершины $v_{18}$:

$c_{74} \cap c_6 = [(v_1,v_{14})+(v_{14},v_8)+(v_8,v_1)] \cap [(v_2,v_1)+(v_1,v_8)+(v_8,v_2)] = (v_1,v_8)$.

Рассмотрим следующую последовательность циклов:

$c_{67} = (v_4,v_8)+(v_8,v_{12})+(v_{12},v_4)$;
$c_{69} = (v_{12},v_8)+(v_8,v_{13})+(v_{13},v_{12})$;
$c_{71} = (v_{13},v_8)+(v_8,v_{14})+(v_{14},v_{13})$;
$c_{74} = (v_1,v_{14})+(v_{14},v_8)+(v_8,v_1)$;
$c_6 = (v_2,v_1)+(v_1,v_8)+(v_8,v_2)$.

Введём мнимые вершины в соответствующие рёбра:

$c_{67} = (v_4,v_8)+(v_8,v_{15})+(v_{15},v_{12})+(v_{12},v_4)$;
$c_{69} = (v_{12},v_{15})+(v_{15},v_8)+(v_8,v_{16})+(v_{16},v_{13})+(v_{13},v_{12})$;
$c_{71} = (v_{13},v_{16})+(v_{16},v_8)+(v_8,v_{17})+(v_{17},v_{14})+(v_{14},v_{13})$;
$c_{74} = (v_1,v_{14})+(v_{14},v_{17})+(v_{17},v_8)+(v_8,v_{18})+(v_{18},v_1)$;
$c_6 = (v_2,v_1)+(v_1,v_{18})+(v_{18},v_8)+(v_8,v_2)$.

Упорядочим запись циклов, оставляя на первом месте ориентированные рёбра с вершинами из кортежа вершин минимального маршрута и отделяя рёбра с мнимыми вершинами:

$c_{67} = [(v_4,v_8)+(v_8,v_{15})]+[(v_{15},v_{12})+(v_{12},v_4)]$;
$c_{69} = [(v_{15},v_8)+(v_8,v_{16})]+[(v_{16},v_{13})+(v_{13},v_{12})+(v_{12},v_{15})]$;
$c_{71} = [(v_{16},v_8)+(v_8,v_{17})]+[(v_{17},v_{14})+(v_{14},v_{13})+(v_{13},v_{16})]$;
$c_{74} = [(v_{17},v_8)+(v_8,v_{18})]+[(v_{18},v_1)+(v_1,v_{14})+(v_{14},v_{17})]$;
$c_6 = [(v_2,v_1)+(v_1,v_{18})]+[(v_{18},v_8)+(v_8,v_2)]$.

Вставим отрезки соединения $(v_4,v_2)$:

$c_{67} = [(v_4,v_8)+(v_8,v_{15})]+[(v_{15},v_4)+(v_4,v_{15})]+[(v_{15},v_{12})+(v_{12},v_4)]$;
$c_{69} = [(v_{15},v_8)+(v_8,v_{16})]+[(v_{16},v_{15})+(v_{15},v_{16})]+[(v_{16},v_{13})+(v_{13},v_{12})+(v_{12},v_{15})]$;
$c_{71} = [(v_{16},v_8)+(v_8,v_{17})]+[(v_{17},v_{16})+(v_{16},v_{17})]+[(v_{17},v_{14})+(v_{14},v_{13})+(v_{13},v_{16})]$;
$c_{74} = [(v_{17},v_8)+(v_8,v_{18})]+[(v_{18},v_{17})+(v_{17},v_{18})]+[(v_{18},v_1)+(v_1,v_{14})+(v_{14},v_{17})]$;

$c_6 = [(v_2,v_1)+(v_1,v_{18})]+[(v_{18},v_2)+(v_2,v_{18})]+[(v_{18},v_8)+(v_8,v_2)]$.

Выделим замкнутую последовательность рёбер и сформируем новую систему циклов вместо циклов $\{c_{67},c_{69},c_{71},c_{74},c_6\}$:

$c_{75} = (v_4,v_8)+(v_8,v_{15})+(v_{15},v_4)$;
$c_{76} = (v_4,v_{15})+(v_{15},v_{12})+(v_{12},v_4)$;
$c_{77} = (v_{15},v_8)+(v_8,v_{16})+(v_{16},v_{15})$;
$c_{78} = (v_{15},v_{16})+(v_{16},v_{13})+(v_{13},v_{12})+(v_{12},v_{15})$;
$c_{79} = (v_{16},v_8)+(v_8,v_{17})+(v_{17},v_{16})$;
$c_{80} = (v_{16},v_{17})+(v_{17},v_{14})+(v_{14},v_{13})+(v_{13},v_{16})$;
$c_{81} = (v_{17},v_8)+(v_8,v_{18})+(v_{18},v_{17})$;
$c_{82} = (v_{17},v_{18})+(v_{18},v_1)+(v_1,v_{14})+(v_{14},v_{17})$;
$c_{83} = (v_2,v_1)+(v_1,v_{18})+(v_{18},v_2)$;
$c_{84} = (v_2,v_{18})+(v_{18},v_8)+(v_8,v_2)$.

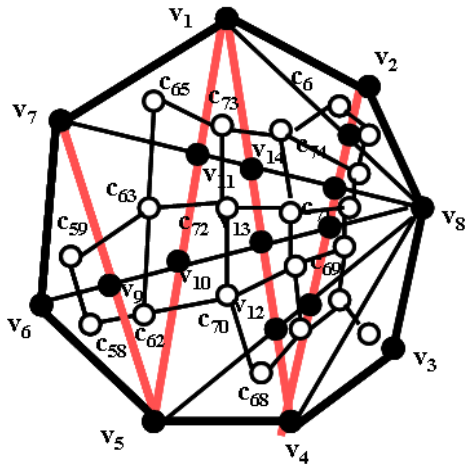

Рис. 18.14. Проведение соединения $(v_4,v_2)$.

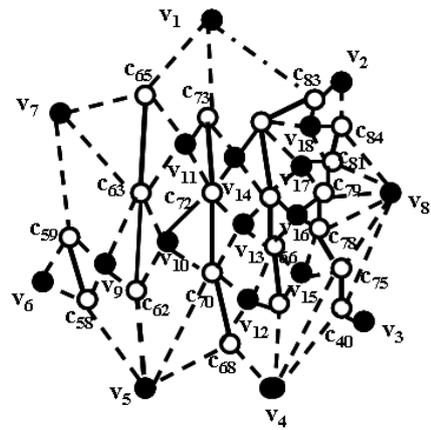

Рис. 18.15. Смешанный граф циклов.

В результате проведения соединений $\{e_{24},e_4,e_3,e_9\}$ образуется множество циклов, характеризующих топологический рисунок:

$c_{40} = (v_8,v_4)+(v_4,v_3)+(v_3,v_8)$;
$c_{58} = (v_5,v_9)+(v_9,v_6)+(v_6,v_5)$;
$c_{59} = (v_7,v_6)+(v_6,v_9)+(v_9,v_7)$;
$c_{62} = (v_5,v_{10})+(v_{10},v_9)+(v_9,v_5)$;
$c_{64} = (v_{10},v_{11})+(v_{11},v_7)+(v_7,v_9)+(v_9,v_{10})$;
$c_{65} = (v_1,v_7)+(v_7,v_{11})+(v_{11},v_1)$;
$c_{68} = (v_4,v_{12})+(v_{12},v_5)+(v_5,v_4)$;
$c_{70} = (v_{12},v_{13})+(v_{13},v_{10})+(v_{10},v_5)+(v_5,v_{12})$;
$c_{72} = (v_{13},v_{14})+(v_{14},v_{11})+(v_{11},v_{10})+(v_{10},v_{13})$;
$c_{73} = (v_1,v_{11})+(v_{11},v_{14})+(v_{14},v_1)$;
$c_{75} = (v_4,v_8)+(v_8,v_{15})+(v_{15},v_4)$;
$c_{76} = (v_4,v_{15})+(v_{15},v_{12})+(v_{12},v_4)$;
$c_{77} = (v_{15},v_8)+(v_8,v_{16})+(v_{16},v_{15})$;
$c_{78} = (v_{15},v_{16})+(v_{16},v_{13})+(v_{13},v_{12})+(v_{12},v_{15})$;
$c_{79} = (v_{16},v_8)+(v_8,v_{17})+(v_{17},v_{16})$;
$c_{80} = (v_{16},v_{17})+(v_{17},v_{14})+(v_{14},v_{13})+(v_{13},v_{16})$;
$c_{81} = (v_{17},v_8)+(v_8,v_{18})+(v_{18},v_{17})$;
$c_{82} = (v_{17},v_{18})+(v_{18},v_1)+(v_1,v_{14})+(v_{14},v_{17})$;
$c_{83} = (v_2,v_1)+(v_1,v_{18})+(v_{18},v_2)$;

$c_{84} = (v_2,v_{18})+(v_{18},v_8)+(v_8,v_2)$;

обод $= (v_1,v_2)+(v_2,v_8)+(v_8,v_3)+(v_3,v_4)+(v_4,v_5)+(v_5,v_6)+(v_6,v_7)+(v_7,v_1)$.

Формирование смешанного графа циклов после проведения соединения $(v_4,v_2)$ показано на рис. 18.15. Все четыре непересекающихся соединения внутри гамильтонова цикла проведены. Выделим следующую часть максимально плоского суграфа между гамильтоновым циклом и ободом. С этой целью из топологического рисунка внутри гамильтонова цикла удалим циклы содержащие мнимые вершины (рис. 18.16).

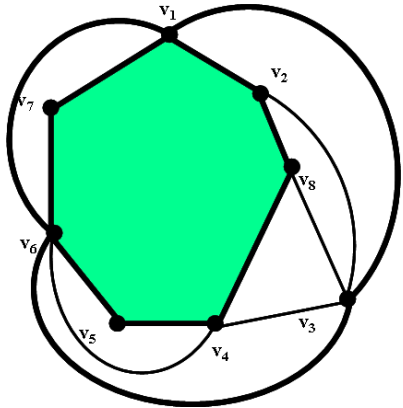 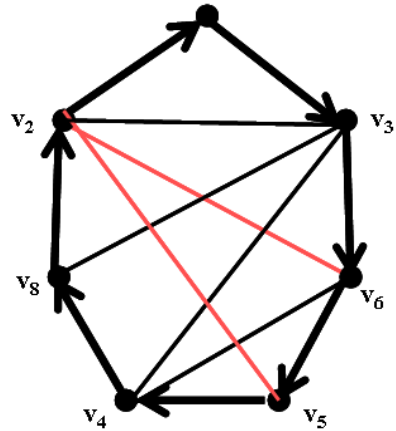

Рис. 18.16. Удаление циклов с мнимыми вершинами.

Рис. 18.17. Проекция непересекающихся соединений на КБС.

### 18.6. Соединение $e_{11} = (v_2,v_6)$

Выделим систему циклов, описывающую выделенную часть максимально плоского суграфа как внутренность обода, за исключением содержащей циклы с мнимыми вершинами части гамильтонова цикла:

$c_1 = (v_1,v_2)+(v_2,v_3)+(v_3,v_1)$;
$c_{47} = (v_4,v_5)+(v_5,v_6)+(v_6,v_4)$;
$c_{40} = (v_8,v_4)+(v_4,v_3)+(v_3,v_8)$;
$c_{26} = (v_8,v_3)+(v_3,v_2)+(v_2,v_8)$;
$c_{38} = (v_4,v_6)+(v_6,v_3)+(v_3,v_4)$;
обод $= (v_2,v_1)+(v_1,v_3)+(v_3,v_6)+(v_6,v_5)+(v_5,v_4)+(v_4,v_8)+(v_8,v_2)$.

Построим координатно-базисную систему векторов, состоящую из ребер контура выделенного подпространства (рис. 18.17). Определим проекцию соединений на КБС. Установим, что соединения $(v_2,v_5)$ и $(v_2,v_6)$ не пересекаются и построим смешанный граф циклов (рис. 18.18).

Выделим минимальный маршрут для проведения соединения $(v_2,v_6)$ и запишем его в виде кортежа $<v_2,c_{26},c_{40},c_{38},v_6>$ с учетом введения мнимых вершин $\{v_{19},v_{20}\}$.

Определим сопряженное неориентированное ребро как пересечение циклов $c_{26} \cap c_{40}$ с включением мнимой вершины $v_{19}$:

$c_{26} \cap c_{40} = [(v_8,v_3)+(v_3,v_2)+(v_2,v_8)] \cap [(v_8,v_4)+(v_4,v_3)+(v_3,v_8)] = (v_8,v_3)$.

Определим сопряженное неориентированное ребро как пересечение циклов $c_{40} \cap c_{38}$ с включением мнимой вершины $v_{20}$:

$c_{40} \cap c_{38} = [(v_8,v_4)+(v_4,v_3)+(v_3,v_8)] \cap [(v_4,v_6)+(v_6,v_3)+(v_3,v_4)] = (v_3,v_4)$.

Рис. 18.18. Граф циклов.  Рис. 18.19. Смешанный граф циклов.

Рассмотрим следующую последовательность циклов:

$c_{26} = (v_8,v_3)+(v_3,v_2)+(v_2,v_8)$;
$c_{40} = (v_8,v_4)+(v_4,v_3)+(v_3,v_8)$;
$c_{38} = (v_4,v_6)+(v_6,v_3)+(v_3,v_4)$.

Введём мнимые вершины в соответствующие рёбра:

$c_{26} = (v_8,v_{19})+(v_{19},v_3)+(v_3,v_2)+(v_2,v_8)$;
$c_{40} = (v_8,v_4)+(v_4,v_{20})+(v_{20},v_3)+(v_3,v_{19})+(v_{19},v_8)$;
$c_{38} = (v_4,v_6)+(v_6,v_3)+(v_3,v_{20})+(v_{20},v_4)$.

Упорядочим запись циклов, оставляя на первом месте ориентированные рёбра с вершинами из кортежа вершин минимального маршрута и отделяя рёбра с мнимыми вершинами:

$c_{26} = [(v_2,v_8)+(v_8,v_{19})]+[(v_{19},v_3)+(v_3,v_2)]$;
$c_{40} = [(v_{19},v_8)+(v_8,v_4)+(v_4,v_{20})]+[(v_{20},v_3)+(v_3,v_{19})]$;
$c_{38} = [(v_6,v_3)+(v_3,v_{20})]+[(v_{20},v_4)+(v_4,v_6)]$.

Вставим отрезки соединения $(v_2,v_6)$:

$c_{26} = [(v_2,v_8)+(v_8,v_{19})]+[(v_{19},v_2)+(v_2,v_{19})]+[(v_{19},v_3)+(v_3,v_2)]$;
$c_{40} = [(v_{19},v_8)+(v_8,v_4)+(v_4,v_{20})]+[(v_{20},v_{19})+(v_{19},v_{20})]+[(v_{20},v_3)+(v_3,v_{19})]$;
$c_{38} = [(v_6,v_3)+(v_3,v_{20})]+[(v_{20},v_6)+(v_6,v_{20})]+[(v_{20},v_4)+(v_4,v_6)]$.

Выделим замкнутую последовательность рёбер и сформируем новую систему циклов вместо циклов $c_{26},c_{40},c_{38}$:

$c_{85} = (v_2,v_8)+(v_8,v_{19})+(v_{19},v_2)$;
$c_{86} = (v_2,v_{19})+(v_{19},v_3)+(v_3,v_2)$;
$c_{87} = (v_{19},v_8)+(v_8,v_4)+(v_4,v_{20})+(v_{20},v_{19})$;
$c_{88} = (v_{19},v_{20})+(v_{20},v_3)+(v_3,v_{19})$;
$c_{89} = (v_6,v_3)+(v_3,v_{20})+(v_{20},v_6)$;
$c_{90} = (v_6,v_{20})+(v_{20},v_4)+(v_4,v_6)$.

Новая система циклов определяет смешанный граф циклов (рис. 18.21):

$c_1 = (v_1,v_2)+(v_2,v_3)+(v_3,v_1);$
$c_{47} = (v_4,v_5)+(v_5,v_6)+(v_6,v_4);$
$c_{85} = (v_2,v_8)+(v_8,v_{19})+(v_{19},v_2);$
$c_{86} = (v_2,v_{19})+(v_{19},v_3)+(v_3,v_2);$
$c_{87} = (v_{19},v_8)+(v_8,v_4)+(v_4,v_{20})+(v_{20},v_{19});$
$c_{88} = (v_{19},v_{20})+(v_{20},v_3)+(v_3,v_{19});$
$c_{89} = (v_6,v_3)+(v_3,v_{20})+(v_{20},v_6);$
$c_{90} = (v_6,v_{20})+(v_{20},v_4)+(v_4,v_6).$

### 18.7. Соединение $e_{10} = (v_2,v_5)$

Выделим минимальный маршрут для проведения соединения $(v_2,v_5)$ и запишем его в виде кортежа $<v_2,c_{85},c_{87},c_{90},c_{47},v_5>$ с учетом введения мнимых вершин $\{v_{21},v_{22},v_{23}\}$.

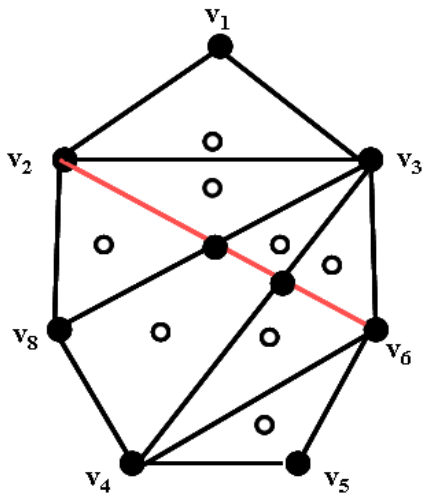 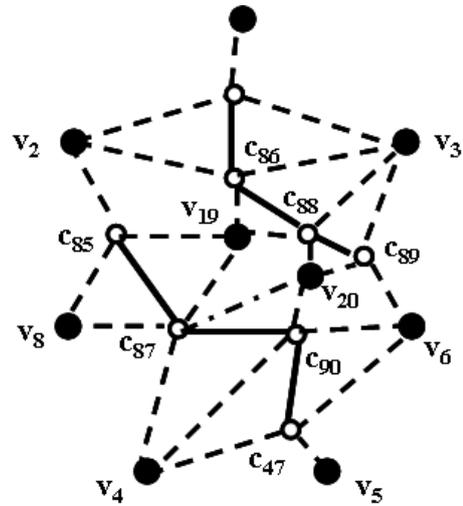

Рис. 18.20. Проведение соединения $(v_2,v_6)$.    Рис. 18.21. Смешанный граф циклов.

Определим сопряженное неориентированное ребро как пересечение циклов $c_{85} \cap c_{87}$ с включением мнимой вершины $v_{21}$:

$c_{85} \cap c_{87} = [(v_2,v_8)+(v_8,v_{19})+(v_{19},v_2)] \cap [(v_{19},v_8)+(v_8,v_4)+(v_4,v_{20})+(v_{20},v_{19})] = (v_8,v_{19}).$

Определим сопряженное неориентированное ребро как пересечение циклов $c_{87} \cap c_{90}$ с включением мнимой вершины $v_{22}$:

$c_{87} \cap c_{90} = [(v_{19},v_8)+(v_8,v_4)+(v_4,v_{20})+(v_{20},v_{19})] \cap [(v_6,v_{20})+(v_{20},v_4)+(v_4,v_6)] = (v_{20},v_4).$

Определим сопряженное неориентированное ребро как пересечение циклов $c_{90} \cap c_{47}$ с включением мнимой вершины $v_{23}$:

$c_{90} \cap c_{47} = [(v_6,v_{20})+(v_{20},v_4)+(v_4,v_6)] \cap [(v_4,v_5)+(v_5,v_6)+(v_6,v_4)] = (v_4,v_6).$

Рассмотрим следующую последовательность циклов:

$c_{85} = (v_2,v_8)+(v_8,v_{19})+(v_{19},v_2);$
$c_{87} = (v_{19},v_8)+(v_8,v_4)+(v_4,v_{20})+(v_{20},v_{19});$
$c_{90} = (v_6,v_{20})+(v_{20},v_4)+(v_4,v_6);$
$c_{47} = (v_4,v_5)+(v_5,v_6)+(v_6,v_4).$

Введём мнимые вершины в соответствующие рёбра:

$c_{85} = (v_2,v_8)+(v_8,v_{21})+(v_{21},v_{19})+(v_{19},v_2);$

$c_{87} = (v_{19},v_{21})+(v_{21},v_8)+(v_8,v_4)+(v_4,v_{22})+(v_{22},v_{20})+(v_{20},v_{19})$;
$c_{90} = (v_6,v_{20})+(v_{20},v_{22})+(v_{22},v_4)+(v_4,v_{23})+(v_{23},v_6)$;
$c_{47} = (v_4,v_5)+(v_5,v_6)+(v_6,v_{23})+(v_{23},v_4)$.

Упорядочим запись циклов, оставляя на первом месте ориентированные ребра с вершинами из кортежа вершин минимального маршрута и отделяя ребра с мнимыми вершинами:

$c_{85} = [(v_2,v_8)+(v_8,v_{21})]+[(v_{21},v_{19})+(v_{19},v_2)]$;
$c_{87} = [(v_{21},v_8)+(v_8,v_4)+(v_4,v_{22})]+[(v_{22},v_{20})+(v_{20},v_{19})+(v_{19},v_{21})]$;
$c_{90} = [(v_{22},v_4)+(v_4,v_{23})]+[(v_{23},v_6)+(v_6,v_{20})+(v_{20},v_{22})]$;
$c_{47} = [(v_5,v_6)+(v_6,v_{23})]+[(v_{23},v_4)+(v_4,v_5)]$.

Вставим отрезки соединения $(v_2,v_6)$:

$c_{85} = [(v_2,v_8)+(v_8,v_{21})]+[(v_{21},v_2)+(v_2,v_{21})]+[(v_{21},v_{19})+(v_{19},v_2)]$;
$c_{87} = [(v_{21},v_8)+(v_8,v_4)+(v_4,v_{22})]+[(v_{22},v_{21})+(v_{21},v_{22})]+[(v_{22},v_{20})+(v_{20},v_{19})+(v_{19},v_{21})]$;
$c_{90} = [(v_{22},v_4)+(v_4,v_{23})]+[(v_{23},v_{22})+(v_{22},v_{23})]+[(v_{23},v_6)+(v_6,v_{20})+(v_{20},v_{22})]$;
$c_{47} = [(v_5,v_6)+(v_6,v_{23})]+[(v_{23},v_5)+(v_5,v_{23})]+[(v_{23},v_4)+(v_4,v_5)]$.

Выделим замкнутую последовательность ребер и сформируем новую систему циклов вместо циклов $c_{85}, c_{87}, c_{90}, c_{47}$:

$c_{91} = (v_2,v_8)+(v_8,v_{21})+(v_{21},v_2)$;
$c_{92} = (v_2,v_{21})+(v_{21},v_{19})+(v_{19},v_2)$;
$c_{93} = (v_{21},v_8)+(v_8,v_4)+(v_4,v_{22})+(v_{22},v_{21})$;
$c_{94} = (v_{21},v_{22})+(v_{22},v_{20})+(v_{20},v_{19})+(v_{19},v_{21})$;
$c_{95} = (v_{22},v_4)+(v_4,v_{23})+(v_{23},v_{22})$;
$c_{96} = (v_{22},v_{23})+(v_{23},v_6)+(v_6,v_{20})+(v_{20},v_{22})$;
$c_{97} = (v_5,v_6)+(v_6,v_{23})+(v_{23},v_5)$;
$c_{98} = (v_5,v_{23})+(v_{23},v_4)+(v_4,v_5)$.

Результат проведения соединений представлен на рис. 18.22. Смешанный граф циклов представлен на рис. 18.23.

Новая система циклов:

$c_1 = (v_1,v_2)+(v_2,v_3)+(v_3,v_1)$;
$c_{86} = (v_2,v_{19})+(v_{19},v_3)+(v_3,v_2)$;
$c_{88} = (v_{19},v_{20})+(v_{20},v_3)+(v_3,v_{19})$;
$c_{89} = (v_6,v_3)+(v_3,v_{20})+(v_{20},v_6)$;
$c_{91} = (v_2,v_8)+(v_8,v_{21})+(v_{21},v_2)$;
$c_{92} = (v_2,v_{21})+(v_{21},v_{19})+(v_{19},v_2)$;
$c_{93} = (v_{21},v_8)+(v_8,v_4)+(v_4,v_{22})+(v_{22},v_{21})$;
$c_{94} = (v_{21},v_{22})+(v_{22},v_{20})+(v_{20},v_{19})+(v_{19},v_{21})$;
$c_{95} = (v_{22},v_4)+(v_4,v_{23})+(v_{23},v_{22})$;
$c_{96} = (v_{22},v_{23})+(v_{23},v_6)+(v_6,v_{20})+(v_{20},v_{22})$;
$c_{97} = (v_5,v_6)+(v_6,v_{23})+(v_{23},v_5)$;
$c_{98} = (v_5,v_{23})+(v_{23},v_4)+(v_4,v_5)$.

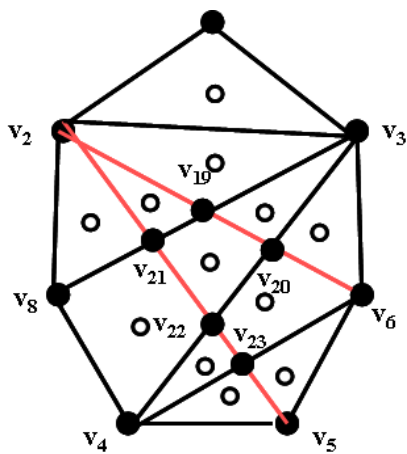 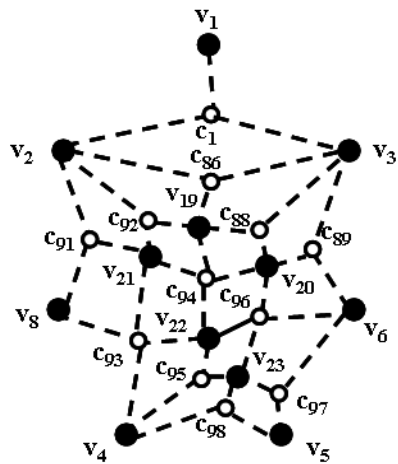

Рис. 18.22. Проведение соединений.      Рис. 18.23. Смешанный граф циклов.

Удалим из гамильтонова цикла цикл $c_{40}$ и присоединим полученную систему. В результате получим следующую систему циклов и обод:

$c_{58} = (v_5,v_9)+(v_9,v_6)+(v_6,v_5)$;
$c_{59} = (v_7,v_6)+(v_6,v_9)+(v_9,v_7)$;
$c_{62} = (v_5,v_{10})+(v_{10},v_9)+(v_9,v_5)$;
$c_{64} = (v_{10},v_{11})+(v_{11},v_7)+(v_7,v_9)+(v_9,v_{10})$;
$c_{65} = (v_1,v_7)+(v_7,v_{11})+(v_{11},v_1)$;
$c_{68} = (v_4,v_{12})+(v_{12},v_5)+(v_5,v_4)$;
$c_{70} = (v_{12},v_{13})+(v_{13},v_{10})+(v_{10},v_5)+(v_5,v_{12})$;
$c_{72} = (v_{13},v_{14})+(v_{14},v_{11})+(v_{11},v_{10})+(v_{10},v_{13})$;
$c_{73} = (v_1,v_{11})+(v_{11},v_{14})+(v_{14},v_1)$;
$c_{75} = (v_4,v_8)+(v_8,v_{15})+(v_{15},v_4)$;
$c_{76} = (v_4,v_{15})+(v_{15},v_{12})+(v_{12},v_4)$;
$c_{77} = (v_{15},v_8)+(v_8,v_{16})+(v_{16},v_{15})$;
$c_{78} = (v_{15},v_{16})+(v_{16},v_{13})+(v_{13},v_{12})+(v_{12},v_{15})$;
$c_{79} = (v_{16},v_8)+(v_8,v_{17})+(v_{17},v_{16})$;
$c_{80} = (v_{16},v_{17})+(v_{17},v_{14})+(v_{14},v_{13})+(v_{13},v_{16})$;
$c_{81} = (v_{17},v_8)+(v_8,v_{18})+(v_{18},v_{17})$;
$c_{82} = (v_{17},v_{18})+(v_{18},v_1)+(v_1,v_{14})+(v_{14},v_{17})$;
$c_{83} = (v_2,v_1)+(v_1,v_{18})+(v_{18},v_2)$;
$c_{84} = (v_2,v_{18})+(v_{18},v_8)+(v_8,v_2)$:
$c_1 = (v_1,v_2)+(v_2,v_3)+(v_3,v_1)$;
$c_{86} = (v_2,v_{19})+(v_{19},v_3)+(v_3,v_2)$;
$c_{88} = (v_{19},v_{20})+(v_{20},v_3)+(v_3,v_{19})$;
$c_{89} = (v_6,v_3)+(v_3,v_{20})+(v_{20},v_6)$;
$c_{91} = (v_2,v_8)+(v_8,v_{21})+(v_{21},v_2)$;
$c_{92} = (v_2,v_{21})+(v_{21},v_{19})+(v_{19},v_2)$;
$c_{93} = (v_{21},v_8)+(v_8,v_4)+(v_4,v_{22})+(v_{22},v_{21})$;
$c_{94} = (v_{21},v_{22})+(v_{22},v_{20})+(v_{20},v_{19})+(v_{19},v_{21})$;
$c_{95} = (v_{22},v_4)+(v_4,v_{23})+(v_{23},v_{22})$;
$c_{96} = (v_{22},v_{23})+(v_{23},v_6)+(v_6,v_{20})+(v_{20},v_{22})$;
$c_{97} = (v_5,v_6)+(v_6,v_{23})+(v_{23},v_5)$;
$c_{98} = (v_5,v_{23})+(v_{23},v_4)+(v_4,v_5)$.
обод = $(v_6,v_7)+(v_7,v_1)+(v_1,v_3)+(v_3,v_6)$.

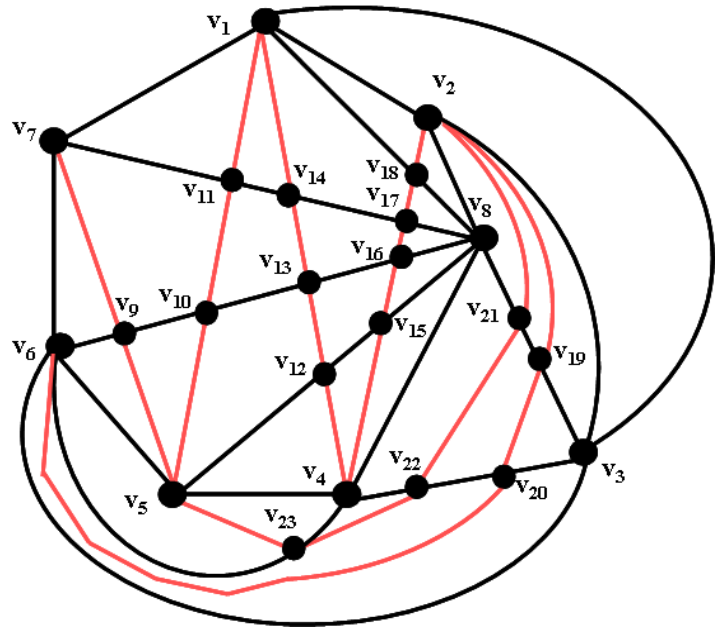

Рис. 18.24. Топологический рисунок с ободом $<v_6,v_7,v_1,v_3>$.

Рассмотрим часть пространства вне цикла $<v_6,v_7,v_1,v_3>$. Поставим в соответствие выделенному ободу координатно-базисную систему векторов и определим количество пересечений соединений исключенных в процессе планаризации (рис. 18.26). Соединения (красный цвет) не пересекаются (рис. 18.25).

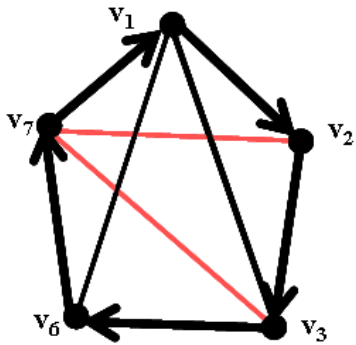

Рис. 18.25. Соединения вне гамильтонова цикла.

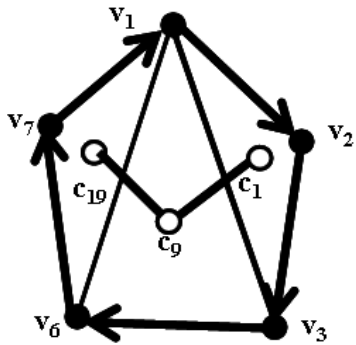

Рис. 18.26. Граф циклов.

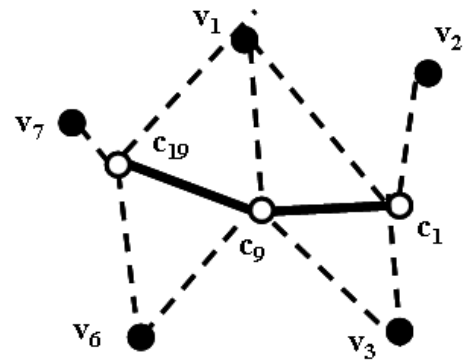

Рис. 18.27. Смешанный граф циклов.

### 18.8. Соединение $e_{17} = (v_7, v_3)$

Построим множество циклов для выделенной части планарного графа:

$c_1 = (v_1,v_2)+(v_2,v_3)+(v_3,v_1)$;
$c_9 = (v_1,v_3)+(v_3,v_6)+(v_6,v_1)$;
$c_{19} = (v_7,v_1)+(v_1,v_6)+(v_6,v_7)$;
обод = $(v_3,v_2)+(v_2,v_1)+(v_1,v_7)+(v_7,v_6)+(v_6,v_3)$.

Выделим минимальный маршрут для проведения соединения $(v_7,v_3)$ и запишем его в виде кортежа $<v_7,c_{19},c_9,v_3>$ с учетом введения мнимой вершины $v_{24}$.

Определим сопряженное неориентированное ребро как пересечение циклов $c_{19} \cap c_9$ с включением мнимой вершины $v_{24}$:

$c_{19} \cap c_9 = [(v_7,v_1)+(v_1,v_6)+(v_6,v_7)] \cap [(v_6,v_1)+(v_1,v_3)+(v_3,v_6)] = (v_6,v_1)$.

Рассмотрим следующую последовательность циклов:

$c_{19} = (v_7,v_1)+(v_1,v_6)+(v_6,v_7)$;
$c_9 = (v_6,v_1)+(v_1,v_3)+(v_3,v_6)$.

Введём мнимые вершины в соответствующие рёбра:

$c_{19} = (v_7,v_1)+(v_1,v_{24})+(v_{24},v_6)+(v_6,v_7)$;
$c_9 = (v_6,v_{24})+(v_{24},v_1)+(v_1,v_3)+(v_3,v_6)$.

Упорядочим запись циклов, оставляя на первом месте ориентированные рёбра с вершинами из кортежа вершин минимального маршрута и отделяя рёбра с мнимыми вершинами:

$c_{19} = [(v_7,v_1)+(v_1,v_{24})]+[(v_{24},v_6)+(v_6,v_7)]$;
$c_9 = [(v_3,v_6)\ (v_6,v_{24})]+[(v_{24},v_1)+(v_1,v_3)]$.

Вставим отрезки соединения $(v_7,v_3)$:

$c_{19} = [(v_7,v_1)+(v_1,v_{24})]+[(v_{24},v_7)+(v_7,v_{24})]+[(v_{24},v_6)+(v_6,v_7)]$;
$c_9 = [(v_3,v_6)\ (v_6,v_{24})]+[(v_{24},v_3)+(v_3,v_{24})]+[(v_{24},v_1)+(v_1,v_3)]$.

Выделим замкнутую последовательность рёбер и сформируем новую систему циклов вместо циклов $c_9, c_{19}$:

$c_{99} = (v_7,v_1)+(v_1,v_{24})+(v_{24},v_7)$;
$c_{100} = (v_7,v_{24})+(v_{24},v_6)+(v_6,v_7)$;
$c_{101} = (v_3,v_6)\ (v_6,v_{24})+(v_{24},v_3)$;
$c_{102} = (v_3,v_{24})+(v_{24},v_1)+(v_1,v_3)$.

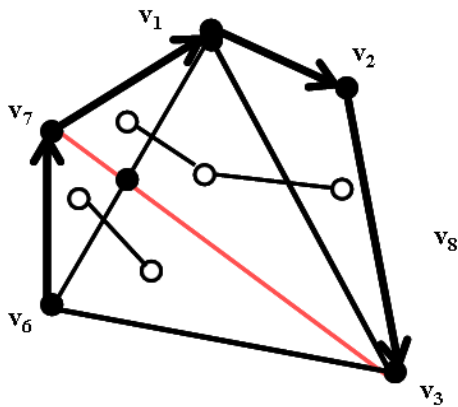
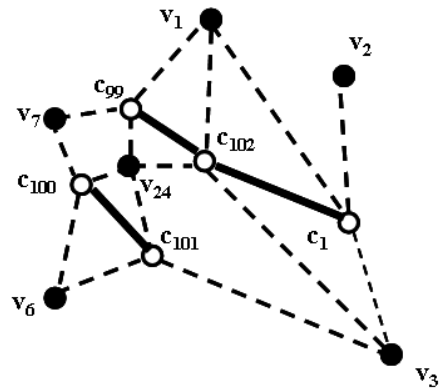

Рис. 18.28. Проведение соединения $(v_7,v_3)$.    Рис. 18.29. Смешанный граф циклов.

В результате получим следующую систему циклов:

$c_1 = (v_1,v_2)+(v_2,v_3)+(v_3,v_1)$;
$c_{99} = (v_7,v_1)+(v_1,v_{24})+(v_{24},v_7)$;
$c_{100} = (v_7,v_{24})+(v_{24},v_6)+(v_6,v_7)$;
$c_{101} = (v_3,v_6)\ (v_6,v_{24})+(v_{24},v_3)$;
$c_{102} = (v_3,v_{24})+(v_{24},v_1)+(v_1,v_3)$;
обод $= (v_3,v_2)+(v_2,v_1)+(v_1,v_7)+(v_7,v_6)+(v_6,v_3)$.

### 18.9. Соединение $e_{12} = (v_7, v_2)$

Выделим минимальный маршрут для проведения соединения $(v_7, v_2)$, запишем минимальный маршрут в виде кортежа $<v_7, c_{99}, c_{102}, c_1, v_2>$ с учетом введения мнимых вершин $\{v_{25}, v_{26}\}$.

Определим сопряженное неориентированное ребро как пересечение циклов $c_{99} \cap c_{102}$ с включением мнимой вершины $v_{25}$:

$c_{99} \cap c_{102} = [(v_7,v_1)+(v_1,v_{24})+(v_{24},v_7)] \cap [(v_3,v_{24})+(v_{24},v_1)+(v_1,v_3)] = (v_{24},v_1)$.

Определим сопряженное неориентированное ребро как пересечение циклов $c_{102} \cap c_1$ с включением мнимой вершины $v_{26}$:

$c_{102} \cap c_1 = [(v_3,v_{24})+(v_{24},v_1)+(v_1,v_3)] \cap [(v_1,v_2)+(v_2,v_3)+(v_3,v_1)] = (v_3,v_1)$.

Рассмотрим следующую последовательность циклов:

$c_{99} = (v_7,v_1)+(v_1,v_{24})+(v_{24},v_7)$;
$c_{102} = (v_3,v_{24})+(v_{24},v_1)+(v_1,v_3)$:
$c_1 = (v_1,v_2)+(v_2,v_3)+(v_3,v_1)$.

Введём мнимые вершины в соответствующие рёбра:

$c_{99} = (v_7,v_1)+(v_1,v_{25})+(v_{25},v_{24})+(v_{24},v_7)$;
$c_{102} = (v_3,v_{24})+(v_{24},v_{25})+(v_{25},v_1)+(v_1,v_{26})+(v_{26},v_3)$;
$c_1 = (v_1,v_2)+(v_2,v_3)+(v_3,v_{26})+(v_{26},v_1)$.

Упорядочим запись циклов, оставляя на первом месте ориентированные ребра с вершинами из кортежа вершин минимального маршрута и отделяя ребра с мнимыми вершинами:

$c_{99} = [(v_7,v_1)+(v_1,v_{25})]+[(v_{25},v_{24})+(v_{24},v_7)]$;
$c_{102} = [(v_{25},v_1)+(v_1,v_{26})]+[(v_{26},v_3)+(v_3,v_{24})+(v_{24},v_{25})]$:
$c_1 = [(v_2,v_3)+(v_3,v_{26})]+[(v_{26},v_1)+(v_1,v_2)]$.

Вставим отрезки соединения $(v_7, v_2)$:

$c_{99} = [(v_7,v_1)+(v_1,v_{25})]+[(v_{25},v_7)+(v_7,v_{25})]+[(v_{25},v_{24})+(v_{24},v_7)]$;
$c_{102} = [(v_{25},v_1)+(v_1,v_{26})]+[(v_{26},v_{25})+(v_{25},v_{26})]+[(v_{26},v_3)+(v_3,v_{24})+(v_{24},v_{25})]$:
$c_1 = [(v_2,v_3)+(v_3,v_{26})]+[(v_{26},v_2)+(v_2,v_{26})]+[(v_{26},v_1)+(v_1,v_2)]$.

Выделим замкнутую последовательность рёбер и сформируем новую систему циклов вместо циклов $\{c_{99}, c_{102}, c_1\}$:

$c_{103} = (v_7,v_1)+(v_1,v_{25})+(v_{25},v_7)$;
$c_{104} = (v_7,v_{25})+(v_{25},v_{24})+(v_{24},v_7)$;
$c_{105} = (v_{25},v_1)+(v_1,v_{26})+(v_{26},v_{25})$;
$c_{106} = (v_{25},v_{26})+(v_{26},v_3)+(v_3,v_{24})+(v_{24},v_{25})$:
$c_{107} = (v_2,v_3)+(v_3,v_{26})+(v_{26},v_2)$;
$c_{108} = (v_2,v_{26})+(v_{26},v_1)+(v_1,v_2)$.

Образуем систему циклов для топологического рисунка:

$c_{100} = (v_7,v_{24})+(v_{24},v_6)+(v_6,v_7)$;
$c_{101} = (v_3,v_6)(v_6,v_{24})+(v_{24},v_3)$;
$c_{103} = (v_7,v_1)+(v_1,v_{25})+(v_{25},v_7)$;
$c_{104} = (v_7,v_{25})+(v_{25},v_{24})+(v_{24},v_7)$;

$c_{105} = (v_{25},v_1)+(v_1,v_{26})+(v_{26},v_{25});$
$c_{106} = (v_{25},v_{26})+(v_{26},v_3)+(v_3,v_{24})+(v_{24},v_{25}):$
$c_{107} = (v_2,v_3)+(v_3,v_{26})+(v_{26},v_2);$
$c_{108} = (v_2,v_{26})+(v_{26},v_1)+(v_1,v_2).$
обод = $(v_3,v_2)+(v_2,v_1)+(v_1,v_7)+(v_7,v_6)+(v_6,v_3).$

Общая система циклов после проведения соединений $\{e_{24},e_4,e_3,e_9,e_{11},e_{10},e_{17},e_{12}\}$ для топологического рисунка (рис. 18.30).

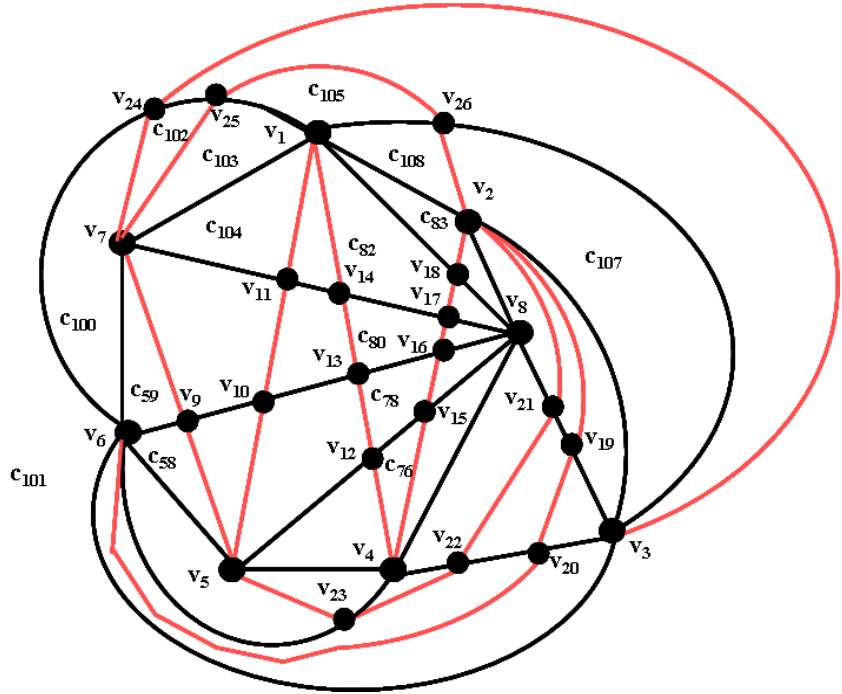

Рис. 18.30. Топологический рисунок с мнимыми вершинами

$c_{58} = (v_5,v_9)+(v_9,v_6)+(v_6,v_5);$
$c_{59} = (v_7,v_6)+(v_6,v_9)+(v_9,v_7);$
$c_{62} = (v_5,v_{10})+(v_{10},v_9)+(v_9,v_5);$
$c_{64} = (v_{10},v_{11})+(v_{11},v_7)+(v_7,v_9)+(v_9,v_{10});$
$c_{65} = (v_1,v_7)+(v_7,v_{11})+(v_{11},v_1);$
$c_{68} = (v_4,v_{12})+(v_{12},v_5)+(v_5,v_4);$
$c_{70} = (v_{12},v_{13})+(v_{13},v_{10})+(v_{10},v_5)+(v_5,v_{12});$
$c_{72} = (v_{13},v_{14})+(v_{14},v_{11})+(v_{11},v_{10})+(v_{10},v_{13});$
$c_{73} = (v_1,v_{11})+(v_{11},v_{14})+(v_{14},v_1);$
$c_{75} = (v_4,v_8)+(v_8,v_{15})+(v_{15},v_4);$
$c_{76} = (v_4,v_{15})+(v_{15},v_{12})+(v_{12},v_4);$
$c_{77} = (v_{15},v_8)+(v_8,v_{16})+(v_{16},v_{15});$
$c_{78} = (v_{15},v_{16})+(v_{16},v_{13})+(v_{13},v_{12})+(v_{12},v_{15});$
$c_{79} = (v_{16},v_8)+(v_8,v_{17})+(v_{17},v_{16});$
$c_{80} = (v_{16},v_{17})+(v_{17},v_{14})+(v_{14},v_{13})+(v_{13},v_{16});$
$c_{81} = (v_{17},v_8)+(v_8,v_{18})+(v_{18},v_{17});$
$c_{82} = (v_{17},v_{18})+(v_{18},v_1)+(v_1,v_{14})+(v_{14},v_{17});$
$c_{83} = (v_2,v_1)+(v_1,v_{18})+(v_{18},v_2);$
$c_{84} = (v_2,v_{18})+(v_{18},v_8)+(v_8,v_2):$
$c_{86} = (v_2,v_{19})+(v_{19},v_3)+(v_3,v_2);$
$c_{88} = (v_{19},v_{20})+(v_{20},v_3)+(v_3,v_{19});$
$c_{89} = (v_6,v_3)+(v_3,v_{20})+(v_{20},v_6);$
$c_{91} = (v_2,v_8)+(v_8,v_{21})+(v_{21},v_2);$

$c_{92} = (v_2,v_{21})+(v_{21},v_{19})+(v_{19},v_2)$;
$c_{93} = (v_{21},v_8)+(v_8,v_4)+(v_4,v_{22})+(v_{22},v_{21})$;
$c_{94} = (v_{21},v_{22})+(v_{22},v_{20})+(v_{20},v_{19})+(v_{19},v_{21})$;
$c_{95} = (v_{22},v_4)+(v_4,v_{23})+(v_{23},v_{22})$;
$c_{96} = (v_{22},v_{23})+(v_{23},v_6)+(v_6,v_{20})+(v_{20},v_{22})$;
$c_{97} = (v_5,v_6)+(v_6,v_{23})+(v_{23},v_5)$;
$c_{98} = (v_5,v_{23})+(v_{23},v_4)+(v_4,v_5)$.
$c_{100} = (v_7,v_{24})+(v_{24},v_6)+(v_6,v_7)$;
$c_{101} = (v_3,v_6)\ (v_6,v_{24})+(v_{24},v_3)$;
$c_{103} = (v_7,v_1)+(v_1,v_{25})+(v_{25},v_7)$;
$c_{104} = (v_7,v_{25})+(v_{25},v_{24})+(v_{24},v_7)$;
$c_{105} = (v_{25},v_1)+(v_1,v_{26})+(v_{26},v_{25})$;
$c_{106} = (v_{25},v_{26})+(v_{26},v_3)+(v_3,v_{24})+(v_{24},v_{25})$:
$c_{107} = (v_2,v_3)+(v_3,v_{26})+(v_{26},v_2)$;
$c_{108} = (v_2,v_{26})+(v_{26},v_1)+(v_1,v_2)$;
обод $= \varnothing$.

### 18.10. Соединение $e_{15} = (v_3,v_5)$

Рассмотрим следующий маршрут $<v_3,c_{101},c_{100},c_{59},c_{58},v_5>$ с введением мнимых вершин $\{v_{27},v_{28},v_{29}\}$:

$c_{101} = (v_3,v_6)\ (v_6,v_{24})+(v_{24},v_3)$;
$c_{100} = (v_7,v_{24})+(v_{24},v_6)+(v_6,v_7)$;
$c_{59} = (v_7,v_6)+(v_6,v_9)+(v_9,v_7)$;
$c_{58} = (v_5,v_9)+(v_9,v_6)+(v_6,v_5)$.

Определим сопряженное неориентированное ребро как пересечение циклов $c_{101} \cap c_{100}$ с включением мнимой вершины $v_{27}$:

$c_{101} \cap c_{100} = [(v_3,v_6)\ (v_6,v_{24})+(v_{24},v_3)] \cap [\ (v_7,v_{24})+(v_{24},v_6)+(v_6,v_7)] = (v_6,v_{24})$.

Определим сопряженное неориентированное ребро как пересечение циклов $c_{19} \cap c_9$ с включением мнимой вершины $v_{28}$:

$c_{100} \cap c_{59} = [(v_7,v_{24})+(v_{24},v_6)+(v_6,v_7)] \cap [(v_7,v_6)+(v_6,v_9)+(v_9,v_7)] = (v_6,v_7)$.

Определим сопряженное неориентированное ребро как пересечение циклов $c_{19} \cap c_9$ с включением мнимой вершины $v_{29}$:

$c_{59} \cap c_{58} = [(v_7,v_6)+(v_6,v_9)+(v_9,v_7)] \cap [\ (v_5,v_9)+(v_9,v_6)+(v_6,v_5)] = (v_6,v_9)$.

Рассмотрим следующую последовательность циклов:

$c_{101} = (v_3,v_6)\ (v_6,v_{24})+(v_{24},v_3)$;
$c_{100} = (v_7,v_{24})+(v_{24},v_6)+(v_6,v_7)$;
$c_{59} = (v_7,v_6)+(v_6,v_9)+(v_9,v_7)$;
$c_{58} = (v_5,v_9)+(v_9,v_6)+(v_6,v_5)$.

Введём мнимые вершины в соответствующие рёбра:

$c_{101} = (v_3,v_6)\ (v_6,v_{27})+(v_{27},v_{24})+(v_{24},v_3)$;
$c_{100} = (v_7,v_{24})+(v_{24},v_{27})+(v_{27},v_6)+(v_6,v_{28})+(v_{28},v_7)$;
$c_{59} = (v_7,v_{28})+(v_{28},v_6)+(v_6,v_{29})+(v_{29},v_9)+(v_9,v_7)$;
$c_{58} = (v_5,v_9)+(v_9,v_{29})+(v_{29},v_6)+(v_6,v_5)$.

Упорядочим запись циклов, оставляя на первом месте ориентированные рёбра с вершинами из кортежа вершин минимального маршрута и отделяя рёбра с мнимыми вершинами:

$c_{101} = [(v_3,v_6)+(v_6,v_{27})]+[(v_{27},v_{24})+(v_{24},v_3)];$
$c_{100} = [(v_{27},v_6)+(v_6,v_{28})]+[(v_{28},v_7)+(v_7,v_{24})+[(v_{24},v_{27})];$
$c_{59} = [(v_{28},v_6)+(v_6,v_{29})]+[(v_{29},v_9)+(v_9,v_7)+(v_7,v_{28})];$
$c_{58} = [(v_5,v_9)+(v_9,v_{29})]+[(v_{29},v_6)+(v_6,v_5)].$

Вставим отрезки соединения $(v_7,v_3)$:

$c_{101} = [(v_3,v_6)+(v_6,v_{27})]+[(v_{27},v_3)+(v_3,v_{27})]+[(v_{27},v_{24})+(v_{24},v_3)];$
$c_{100} = [(v_{27},v_6)+(v_6,v_{28})]+[(v_{28},v_{29})+(v_{29},v_{28})]+[(v_{28},v_7)+(v_7,v_{24})+[(v_{24},v_{27})];$
$c_{59} = [(v_{28},v_6)+(v_6,v_{29})]+[(v_{29},v_{28})+(v_{28},v_{29})]+[(v_{29},v_9)+(v_9,v_7)+(v_7,v_{28})];$
$c_{58} = [(v_5,v_9)+(v_9,v_{29})]+[(v_{29},v_5)+(v_5,v_{29})]+[(v_{29},v_6)+(v_6,v_5)].$

Выделим замкнутую последовательность рёбер и сформируем новую систему циклов вместо циклов $\{c_{101},c_{100},c_{59},c_{58}\}$:

$c_{109} = (v_3,v_6)+(v_6,v_{27})+(v_{27},v_3);$
$c_{110} = (v_3,v_{27})+(v_{27},v_{24})+(v_{24},v_3);$
$c_{111} = (v_{27},v_6)+(v_6,v_{28})+(v_{28},v_{29});$
$c_{112} = (v_{29},v_{28})+(v_{28},v_7)+(v_7,v_{24})+(v_{24},v_{27});$
$c_{113} = (v_{28},v_6)+(v_6,v_{29})+(v_{29},v_{28});$
$c_{114} = (v_{28},v_{29})+(v_{29},v_9)+(v_9,v_7)+(v_7,v_{28});$
$c_{115} = (v_5,v_9)+(v_9,v_{29})+(v_{29},v_5);$
$c_{116} = (v_5,v_{29})+(v_{29},v_6)+(v_6,v_5).$

### 18.11. Соединение $e_{21} = (v_7,v_4)$

Рассмотрим следующий маршрут $<v_7,c_{103},c_{104},c_{107},c_{83},c_{82},c_{80},c_{78},c_{76},v_4>$ с введением мнимых вершин $\{v_{30},v_{31},v_{32},v_{33},v_{34},v_{35},v_{36}\}$:

$c_{103} = (v_7,v_1)+(v_1,v_{25})+(v_{25},v_7);$
$c_{105} = (v_{25},v_1)+(v_1,v_{26})+(v_{26},v_{25});$
$c_{108} = (v_2,v_{26})+(v_{26},v_1)+(v_1,v_2);$
$c_{83} = (v_2,v_1)+(v_1,v_{18})+(v_{18},v_2);$
$c_{82} = (v_{17},v_{18})+(v_{18},v_1)+(v_1,v_{14})+(v_{14},v_{17});$
$c_{80} = (v_{16},v_{17})+(v_{17},v_{14})+(v_{14},v_{13})+(v_{13},v_{16});$
$c_{78} = (v_{15},v_{16})+(v_{16},v_{13})+(v_{13},v_{12})+(v_{12},v_{15});$
$c_{76} = (v_4,v_{15})+(v_{15},v_{12})+(v_{12},v_4).$

Определим сопряжённое неориентированное ребро как пересечение циклов $c_{102} \cap c_{105}$ с включением мнимой вершины $v_{30}$:

$c_{102} \cap c_{105} = [(v_7,v_1)+(v_1,v_{25})+(v_{25},v_7)] \cap [(v_{25},v_1)+(v_1,v_{26})+(v_{26},v_{25})] = (v_1,v_{25}).$

Определим сопряжённое неориентированное ребро как пересечение циклов $c_{105} \cap c_{108}$ с включением мнимой вершины $v_{31}$:

$c_{105} \cap c_{108} = [(v_{25},v_1)+(v_1,v_{26})+(v_{26},v_{25})] \cap [(v_2,v_{26})+(v_{26},v_1)+(v_1,v_2)] = (v_1,v_{26}).$

Определим сопряжённое неориентированное ребро как пересечение циклов $c_{108} \cap c_{83}$ с включением мнимой вершины $v_{32}$:

$c_{108} \cap c_{83} = [(v_2,v_{26})+(v_{26},v_1)+(v_1,v_2))] \cap [(v_2,v_1)+(v_1,v_{18})+(v_{18},v_2)] = (v_1,v_2).$

Определим сопряженное неориентированное ребро как пересечение циклов $c_{83} \cap c_{82}$ с включением мнимой вершины $v_{33}$:

$c_{83} \cap c_{82} = [(v_2,v_1)+(v_1,v_{18})+(v_{18},v_2)] \cap [(v_{17},v_{18})+(v_{18},v_1)+(v_1,v_{14})+(v_{14},v_{17})] = (v_1,v_{18})$.

Определим сопряженное неориентированное ребро как пересечение циклов $c_{82} \cap c_{80}$ с включением мнимой вершины $v_{34}$:

$c_{82} \cap c_{80} = [(v_{17},v_{18})+(v_{18},v_1)+(v_1,v_{14})+(v_{14},v_{17})] \cap [(v_{16},v_{17})+(v_{17},v_{14})+(v_{14},v_{13})+(v_{13},v_{16})] = (v_{14},v_{17})$.

Определим сопряженное неориентированное ребро как пересечение циклов $c_{80} \cap c_{78}$ с включением мнимой вершины $v_{35}$:

$c_{80} \cap c_{78} = [(v_{16},v_{17})+(v_{17},v_{14})+(v_{14},v_{13})+(v_{13},v_{16})] \cap [(v_{15},v_{16})+(v_{16},v_{13})+(v_{13},v_{12})+(v_{12},v_{15})] = (v_{13},v_{16})$.

Определим сопряженное неориентированное ребро как пересечение циклов $c_{78} \cap c_{76}$ с включением мнимой вершины $v_{36}$:

$c_{78} \cap c_{76} = [(v_{15},v_{16})+(v_{16},v_{13})+(v_{13},v_{12})+(v_{12},v_{15})] \cap [(v_4,v_{15})+(v_{15},v_{12})+(v_{12},v_4)] = (v_{12},v_{15})$.

Рассмотрим следующую последовательность циклов:

$c_{103} = (v_7,v_1)+(v_1,v_{25})+(v_{25},v_7)$;
$c_{105} = (v_{25},v_1)+(v_1,v_{26})+(v_{26},v_{25})$;
$c_{108} = (v_2,v_{26})+(v_{26},v_1)+(v_1,v_2)$;
$c_{83} = (v_2,v_1)+(v_1,v_{18})+(v_{18},v_2)$;
$c_{82} = (v_{17},v_{18})+(v_{18},v_1)+(v_1,v_{14})+(v_{14},v_{17})$;
$c_{80} = (v_{16},v_{17})+(v_{17},v_{14})+(v_{14},v_{13})+(v_{13},v_{16})$;
$c_{78} = (v_{15},v_{16})+(v_{16},v_{13})+(v_{13},v_{12})+(v_{12},v_{15})$;
$c_{76} = (v_4,v_{15})+(v_{15},v_{12})+(v_{12},v_4)$.

Введём мнимые вершины в соответствующие рёбра:

$c_{103} = (v_7,v_1)+(v_1,v_{30})+(v_{30},v_{25})+(v_{25},v_7)$;
$c_{105} = (v_{25},v_{30})+(v_{30},v_1)+(v_1,v_{31})+(v_{31},v_{26})+(v_{26},v_{25})$;
$c_{108} = (v_2,v_{26})+(v_{26},v_{31})+(v_{31},v_1)+(v_1,v_{32})+(v_{32},v_2)$;
$c_{83} = (v_2,v_{32})+(v_{32},v_1)+(v_1,v_{33})+(v_{33},v_{18})+(v_{18},v_2)$;
$c_{82} = (v_{17},v_{18})+(v_{18},v_{33})+(v_{33},v_1)+(v_1,v_{14})+(v_{14},v_{34})+(v_{34},v_{17})$;
$c_{80} = (v_{16},v_{17})+(v_{17},v_{34})+(v_{34},v_{14})+(v_{14},v_{13})+(v_{13},v_{35})+(v_{35},v_{16})$;
$c_{78} = (v_{15},v_{16})+(v_{16},v_{35})+(v_{35},v_{13})+(v_{13},v_{12})+(v_{12},v_{36})+(v_{36},v_{15})$;
$c_{76} = (v_4,v_{15})+(v_{15},v_{36})+(v_{36},v_{12})+(v_{12},v_4)$.

Упорядочим запись циклов, оставляя на первом месте ориентированные рёбра с вершинами из кортежа вершин минимального маршрута и отделяя рёбра с мнимыми вершинами:

$c_{103} = [(v_7,v_1)+(v_1,v_{30})]+[(v_{30},v_{25})+(v_{25},v_7)]$;
$c_{105} = [(v_{30},v_1)+(v_1,v_{31})]+[(v_{31},v_{26})+(v_{26},v_{25})+(v_{25},v_{30})[$;
$c_{108} = [(v_{31},v_1)+(v_1,v_{32})]+[(v_{32},v_2)+(v_2,v_{26})+(v_{26},v_{31})]$;
$c_{83} = [(v_{32},v_1)+(v_1,v_{33})]+[(v_{33},v_{18})+(v_{18},v_2)+(v_2,v_{32})[$;
$c_{82} = [(v_{33},v_1)+(v_1,v_{14})+(v_{14},v_{34})]+[(v_{34},v_{17})+(v_{17},v_{18})+(v_{18},v_{33})]$;
$c_{80} = [(v_{34},v_{14})+(v_{14},v_{13})+(v_{13},v_{35})]+[(v_{35},v_{16})+(v_{16},v_{17})+(v_{17},v_{34})]$;
$c_{78} = [(v_{35},v_{13})+(v_{13},v_{12})+(v_{12},v_{36})]+[(v_{36},v_{15})+(v_{15},v_{16})+(v_{16},v_{35})]$;
$c_{76} = [(v_4,v_{15})+(v_{15},v_{36})]+[(v_{36},v_{12})+(v_{12},v_4)]$.

Вставим отрезки соединения $(v_7,v_4)$:

$c_{103} = [(v_7,v_1)+(v_1,v_{30})]+[(v_{30},v_7)+(v_7,v_{30})]+[(v_{30},v_{25})+(v_{25},v_7)];$
$c_{105} = [(v_{30},v_1)+(v_1,v_{31})]+[(v_{31},v_{30})+(v_{30},v_{31})]+[(v_{31},v_{26})+(v_{26},v_{25})+(v_{25},v_{30})[;$
$c_{108} = [(v_{31},v_1)+(v_1,v_{32})]+[(v_{32},v_{31})+(v_{31},v_{32})]+[(v_{32},v_2)+(v_2,v_{26})+(v_{26},v_{31})];$
$c_{83} = [(v_{32},v_1)+(v_1,v_{33})]+[(v_{33},v_{32})+(v_{32},v_{33})]+[(v_{33},v_{18})+(v_{18},v_2)+(v_2,v_{32})[;$
$c_{82} = [(v_{33},v_1)+(v_1,v_{14})+(v_{14},v_{34})]+[(v_{34},v_{35})+(v_{35},v_{34})]+[(v_{34},v_{17})+(v_{17},v_{18})+(v_{18},v_{33})];$
$c_{80} = [(v_{34},v_{14})+(v_{14},v_{13})+(v_{13},v_{35})]+[(v_{35},v_{36})+(v_{36},v_{35})]+[(v_{35},v_{16})+(v_{16},v_{17})+(v_{17},v_{34})];$
$c_{78} = [(v_{35},v_{13})+(v_{13},v_{12})+(v_{12},v_{36})]+[(v_{36},v_{35})+(v_{35},v_{36})]+[(v_{36},v_{15})+(v_{15},v_{16})+(v_{16},v_{35})];$
$c_{76} = [(v_4,v_{15})+(v_{15},v_{36})]+[(v_{36},v_4)+(v_4,v_{36})]+[(v_{36},v_{12})+(v_{12},v_4)].$

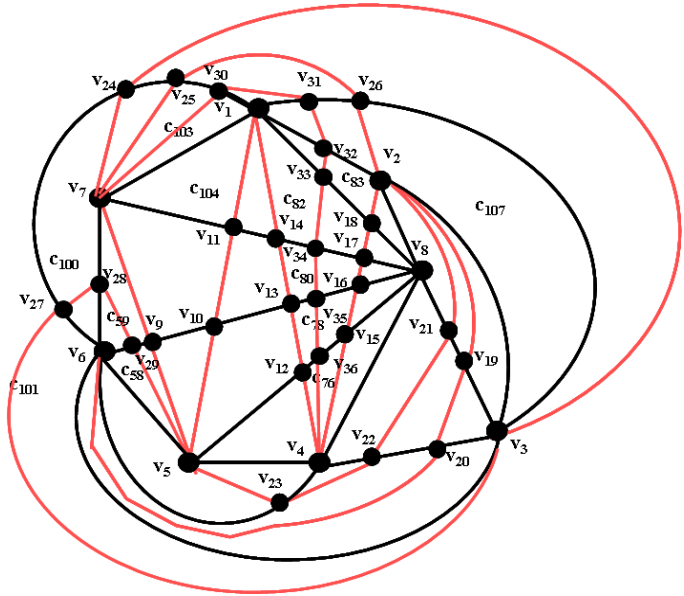

Рис. 18.31. Топологический рисунок непланарного графа
после проведения всех соединений.

Выделим замкнутую последовательность ребер и сформируем новую систему циклов вместо циклов $\{c_{103}, c_{104}, c_{108}, c_{83}, c_{82}, c_{80}, c_{78}, c_{76}\}$:

$c_{117} = (v_7,v_1)+(v_1,v_{30})+(v_{30},v_7);$
$c_{118} = (v_7,v_{30})+(v_{30},v_{25})+(v_{25},v_7);$
$c_{119} = (v_{30},v_1)+(v_1,v_{31})+(v_{31},v_{30});$
$c_{120} = (v_{30},v_{31})+(v_{31},v_{26})+(v_{26},v_{25})+(v_{25},v_{30});$
$c_{121} = (v_{31},v_1)+(v_1,v_{32})+(v_{32},v_{31});$
$c_{122} = (v_{31},v_{32})+(v_{32},v_2)+(v_2,v_{26})+(v_{26},v_{31});$
$c_{123} = (v_{32},v_1)+(v_1,v_{33})+(v_{33},v_{32});$
$c_{124} = (v_{32},v_{33})+(v_{33},v_{18})+(v_{18},v_2)+(v_2,v_{32});$
$c_{125} = (v_{33},v_1)+(v_1,v_{14})+(v_{14},v_{34})+(v_{34},v_{35});$
$c_{126} = (v_{35},v_{34})+(v_{34},v_{17})+(v_{17},v_{18})+(v_{18},v_{33});$
$c_{127} = (v_{34},v_{14})+(v_{14},v_{13})+(v_{13},v_{35})+(v_{35},v_{36});$
$c_{128} = (v_{36},v_{35})+(v_{35},v_{16})+(v_{16},v_{17})+(v_{17},v_{34});$
$c_{129} = (v_{35},v_{13})+(v_{13},v_{12})+(v_{12},v_{36})+(v_{36},v_{35});$
$c_{130} = (v_{35},v_{36})+(v_{36},v_{15})+(v_{15},v_{16})+(v_{16},v_{35});$
$c_{131} = (v_4,v_{15})+(v_{15},v_{36})+(v_{36},v_4);$
$c_{132} = (v_4,v_{36})+(v_{36},v_{12})+(v_{12},v_4).$

Система циклов, описывающая топологический рисунок непланарного графа после проведения соединений (рис. 18.31):

$c_{62} = (v_5,v_{10})+(v_{10},v_9)+(v_9,v_5);$
$c_{64} = (v_{10},v_{11})+(v_{11},v_7)+(v_7,v_9)+(v_9,v_{10});$

$c_{65} = (v_1,v_7)+(v_7,v_{11})+(v_{11},v_1);$
$c_{68} = (v_4,v_{12})+(v_{12},v_5)+(v_5,v_4);$
$c_{70} = (v_{12},v_{13})+(v_{13},v_{10})+(v_{10},v_5)+(v_5,v_{12});$
$c_{72} = (v_{13},v_{14})+(v_{14},v_{11})+(v_{11},v_{10})+(v_{10},v_{13});$
$c_{73} = (v_1,v_{11})+(v_{11},v_{14})+(v_{14},v_1);$
$c_{75} = (v_4,v_8)+(v_8,v_{15})+(v_{15},v_4);$
$c_{77} = (v_{15},v_8)+(v_8,v_{16})+(v_{16},v_{15});$
$c_{79} = (v_{16},v_8)+(v_8,v_{17})+(v_{17},v_{16});$
$c_{81} = (v_{17},v_8)+(v_8,v_{18})+(v_{18},v_{17});$
$c_{84} = (v_2,v_{18})+(v_{18},v_8)+(v_8,v_2):$
$c_{86} = (v_2,v_{19})+(v_{19},v_3)+(v_3,v_2);$
$c_{88} = (v_{19},v_{20})+(v_{20},v_3)+(v_3,v_{19});$
$c_{89} = (v_6,v_3)+(v_3,v_{20})+(v_{20},v_6);$
$c_{91} = (v_2,v_8)+(v_8,v_{21})+(v_{21},v_2);$
$c_{92} = (v_2,v_{21})+(v_{21},v_{19})+(v_{19},v_2);$
$c_{93} = (v_{21},v_8)+(v_8,v_4)+(v_4,v_{22})+(v_{22},v_{21});$
$c_{94} = (v_{21},v_{22})+(v_{22},v_{20})+(v_{20},v_{19})+(v_{19},v_{21});$
$c_{95} = (v_{22},v_4)+(v_4,v_{23})+(v_{23},v_{22});$
$c_{96} = (v_{22},v_{23})+(v_{23},v_6)+(v_6,v_{20})+(v_{20},v_{22});$
$c_{97} = (v_5,v_6)+(v_6,v_{23})+(v_{23},v_5);$
$c_{98} = (v_5,v_{23})+(v_{23},v_4)+(v_4,v_5).$
$c_{102} = (v_7,v_{25})+(v_{25},v_{24})+(v_{24},v_7);$
$c_{106} = (v_{25},v_{26})+(v_{26},v_3)+(v_3,v_{24})+(v_{24},v_{25}):$
$c_{107} = (v_2,v_3)+(v_3,v_{26})+(v_{26},v_2);$
$c_{109} = (v_3,v_6)+(v_6,v_{27})+(v_{27},v_3);$
$c_{110} = (v_3,v_{27})+(v_{27},v_{24})+(v_{24},v_3);$
$c_{111} = (v_{27},v_6)+(v_6,v_{28})+(v_{28},v_{29});$
$c_{112} = (v_{29},v_{28})+(v_{28},v_7)+(v_7,v_{24})+(v_{24},v_{27});$
$c_{113} = (v_{28},v_6)+(v_6,v_{29})+(v_{29},v_{28});$
$c_{114} = (v_{28},v_{29})+(v_{29},v_9)+(v_9,v_7)+(v_7,v_{28});$
$c_{115} = (v_5,v_9)+(v_9,v_{29})+(v_{29},v_5);$
$c_{116} = (v_5,v_{29})+(v_{29},v_6)+(v_6,v_5).$
$c_{117} = (v_7,v_1)+(v_1,v_{30})+(v_{30},v_7);$
$c_{118} = (v_7,v_{30})+(v_{30},v_{25})+(v_{25},v_7);$
$c_{119} = (v_{30},v_1)+(v_1,v_{31})+(v_{31},v_{30});$
$c_{120} = (v_{30},v_{31})+(v_{31},v_{26})+(v_{26},v_{25})+(v_{25},v_{30});$
$c_{121} = (v_{31},v_1)+(v_1,v_{32})+(v_{32},v_{31});$
$c_{122} = (v_{31},v_{32})+(v_{32},v_2)+(v_2,v_{26})+(v_{26},v_{31});$
$c_{123} = (v_{32},v_1)+(v_1,v_{33})+(v_{33},v_{32});$
$c_{124} = (v_{32},v_{33})+(v_{33},v_{18})+(v_{18},v_2)+(v_2,v_{32});$
$c_{125} = (v_{33},v_1)+(v_1,v_{14})+(v_{14},v_{34})+(v_{34},v_{35});$
$c_{126} = (v_{35},v_{34})+(v_{34},v_{17})+(v_{17},v_{18})+(v_{18},v_{33});$
$c_{127} = (v_{34},v_{14})+(v_{14},v_{13})+(v_{13},v_{35})+(v_{35},v_{36});$
$c_{128} = (v_{36},v_{35})+(v_{35},v_{16})+(v_{16},v_{17})+(v_{17},v_{34});$
$c_{129} = (v_{35},v_{13})+(v_{13},v_{12})+(v_{12},v_{36})+(v_{36},v_{35});$
$c_{130} = (v_{35},v_{36})+(v_{36},v_{15})+(v_{15},v_{16})+(v_{16},v_{35});$
$c_{131} = (v_4,v_{15})+(v_{15},v_{36})+(v_{36},v_4);$
$c_{132} = (v_4,v_{36})+(v_{36},v_{12})+(v_{12},v_4).$
обод $= \varnothing$ .

## 18.12. Построение каналов для проведения маршрутов

**Определение 18.1.** Цепочка вершин в графе циклов, соединенных сопряженными ребрами не принадлежащими соединениям удаленным в процессе планаризации, называется *каналом*.

Построим множество циклов для выделенной части планарного графа:

$c_1 = (v_1,v_2)+(v_2,v_3)+(v_3,v_1)$;
$c_9 = (v_1,v_3)+(v_3,v_6)+(v_6,v_1)$;
$c_{19} = (v_7,v_1)+(v_1,v_6)+(v_6,v_7)$;
обод = $(v_3,v_2)+(v_2,v_1)+(v_1,v_7)+(v_7,v_6)+(v_6,v_3)$.

Введение мнимой вершины $v_{24}$ разбивает соединение $(v_7,v_3)$ на части $(v_7,v_{24})$ и $(v_{24},v_3)$. После проведения соединения $(v_7,v_3)$, образуются циклы:

$c_{99} = (v_7,v_1)+(v_1,v_{24})+(v_{24},v_7)$;
$c_{100} = (v_7,v_{24})+(v_{24},v_6)+(v_6,v_7)$;
$c_{101} = (v_3,v_6)\,(v_6,v_{24})+(v_{24},v_3)$;
$c_{102} = (v_3,v_{24})+(v_{24},v_1)+(v_1,v_3)$.

Циклы $c_{99}$ и $c_{102}$ имеют сопряженное ребро $(v_7,v_{24})$, а циклы $c_{100}$ и $c_{101}$ – сопряженное ребро $(v_{24},v_3)$. Для построения каналов соответствующие ребра $(c_{99},c_{102})$ и $(c_{100},c_{101})$ должны быть удалены из смешанного графа циклов (рис. 18.33).

В результате получим следующую систему циклов:

$c_1 = (v_1,v_2)+(v_2,v_3)+(v_3,v_1)$;
$c_{99} = (v_7,v_1)+(v_1,v_{24})+(v_{24},v_7)$;
$c_{100} = (v_7,v_{24})+(v_{24},v_6)+(v_6,v_7)$;
$c_{101} = (v_3,v_6)\,(v_6,v_{24})+(v_{24},v_3)$;
$c_{102} = (v_3,v_{24})+(v_{24},v_1)+(v_1,v_3)$.
обод = $(v_3,v_2)+(v_2,v_1)+(v_1,v_7)+(v_7,v_6)+(v_6,v_3)$.

Введение мнимых вершин $v_{25}$ и $v_{26}$ разбивает соединение $(v_7,v_2)$ на части $(v_7,v_{25})$, $(v_{25},v_{26})$ и $(v_{26},v_2)$. После проведения соединения $(v_7,v_3)$, образуются циклы:

$c_{103} = (v_7,v_1)+(v_1,v_{25})+(v_{25},v_7)$;
$c_{104} = (v_7,v_{25})+(v_{25},v_{24})+(v_{24},v_7)$;
$c_{105} = (v_{25},v_1)+(v_1,v_{26})+(v_{26},v_{25})$;
$c_{106} = (v_{25},v_{26})+(v_{26},v_3)+(v_3,v_{24})+(v_{24},v_{25})$:
$c_{107} = (v_2,v_3)+(v_3,v_{26})+(v_{26},v_2)$;
$c_{108} = (v_2,v_{26})+(v_{26},v_1)+(v_1,v_2)$.

Циклы $c_{103}$ и $c_{104}$ имеют сопряженное ребро $(v_7,v_{24})$, циклы $c_{105}$ и $c_{106}$ – сопряженное ребро $(v_{25},v_{26})$, а циклы $c_{107}$ и $c_{108}$ – $(v_2,v_{26})$. Для построения каналов соответствующие ребра $(c_{103},c_{104})$, $(c_{105},c_{106})$ и $(c_{107},c_{108})$ должны быть удалены из смешанного графа циклов (рис. 18.34).

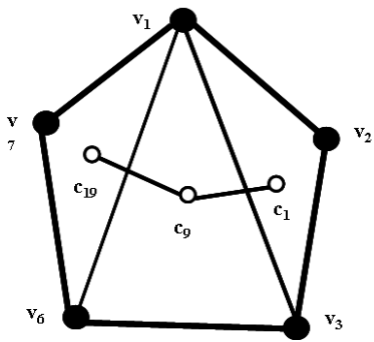 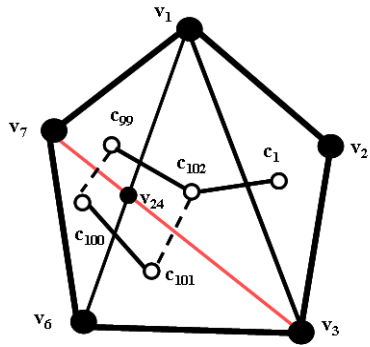 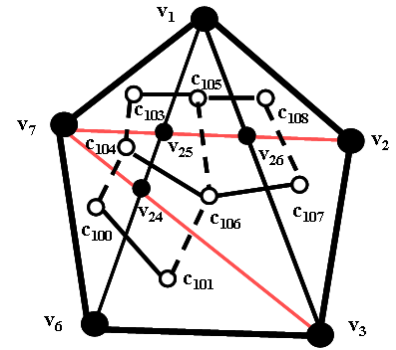

| Рис. 18.32. Граф циклов. | Рис. 18.33. Удаление ребер в графе циклов вне обода. | Рис. 18.34. Граф циклов после проведения соединений $e_{15}, e_{21}$. |

Образуем систему циклов для топологического рисунка:

$c_{100} = (v_7,v_{24})+(v_{24},v_6)+(v_6,v_7)$;
$c_{101} = (v_3,v_6)\,(v_6,v_{24})+(v_{24},v_3)$;
$c_{103} = (v_7,v_1)+(v_1,v_{25})+(v_{25},v_7)$;
$c_{104} = (v_7,v_{25})+(v_{25},v_{24})+(v_{24},v_7)$;
$c_{105} = (v_{25},v_1)+(v_1,v_{26})+(v_{26},v_{25})$;
$c_{106} = (v_{25},v_{26})+(v_{26},v_3)+(v_3,v_{24})+(v_{24},v_{25})$:
$c_{107} = (v_2,v_3)+(v_3,v_{26})+(v_{26},v_2)$;
$c_{108} = (v_2,v_{26})+(v_{26},v_1)+(v_1,v_2)$;
обод = $(v_3,v_2)+(v_2,v_1)(v_1,v_7)+(v_7,v_6)+(v_6,v_3)$.

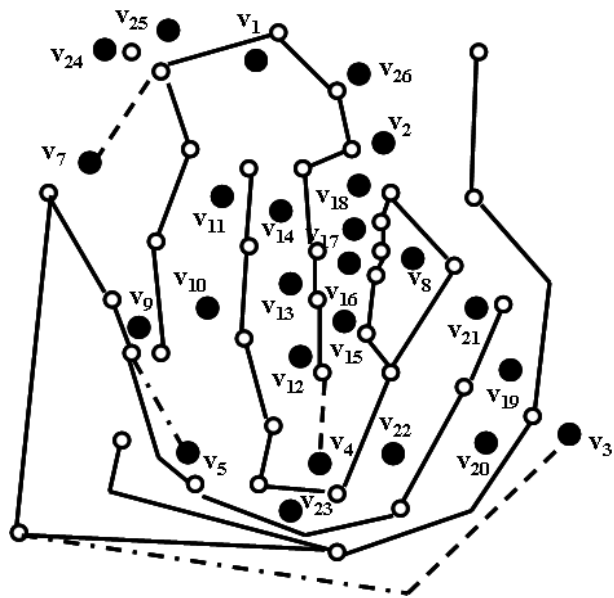

Рис. 18.35. Построение каналов для проведения маршрутов.

На рис. 18.35 показано распределение каналов для проведения соединений $e_{15}$ и $e_{21}$.

**Комментарии**

В главе подробно рассмотрен метод формирования системы циклов для процедуры проведения соединений. На примере графа К$_8$ показан процесс построения каналов для проведения соединений, исключённых в процессе планаризации.

## Глава 19. ТОПОЛОГИЧЕСКИЙ РИСУНОК СЛОЯ

### 19.1. Процесс образования слоя.

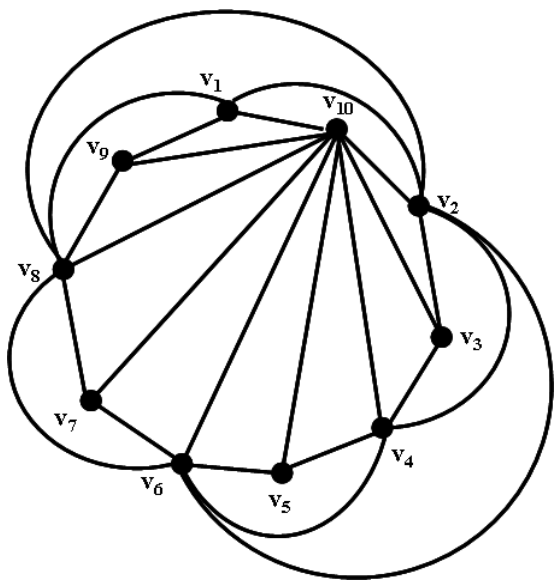

Рис. 19.1. Максимально плоский суграф графа К$_{10}$.

Рассмотрим полный граф К$_{10}$. Выделим максимально плоский суграф графа К$_{10}$ (рис. 19.1). Гамильтонов цикл разделяет топологическое пространство на две части: внешнюю и внутреннюю (рис. 19.2). На ребрах гамильтонова цикла графа К$_{10}$ построим координатно-базисную систему векторов и определим проекции соединений удаленных в процессе планаризации (рис. 19.3).

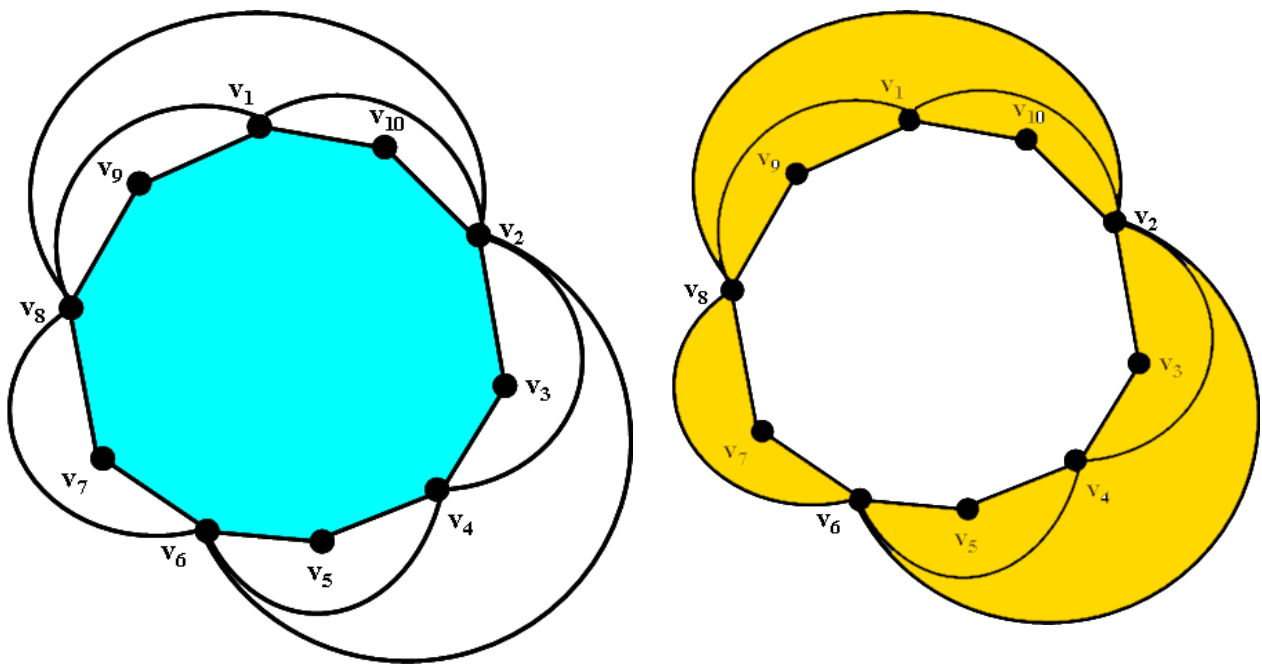

Рис. 19.2. Внутренняя и внешняя области гамильтонова цикла.

Среди множества удаленных в процессе планаризации соединений выделим подмножество непересекающихся соединений (на рис. 19.4 непересекающиеся соединения обозначены красным цветом).

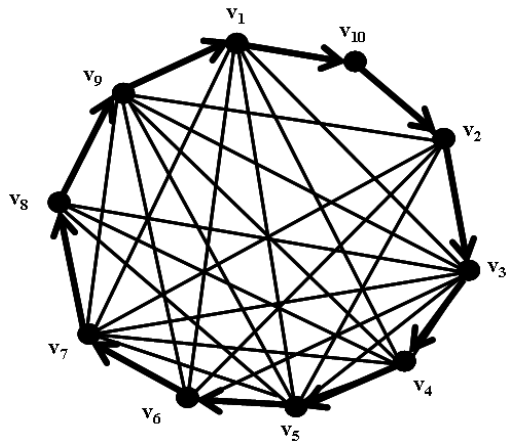

Рис. 19.3. Координатно-базисная система и удаленные соединения.

Относительно выбранных непересекающихся соединений построим топологический рисунок внутренней части гамильтонова цикла с мнимыми вершинами (рис. 19.4).

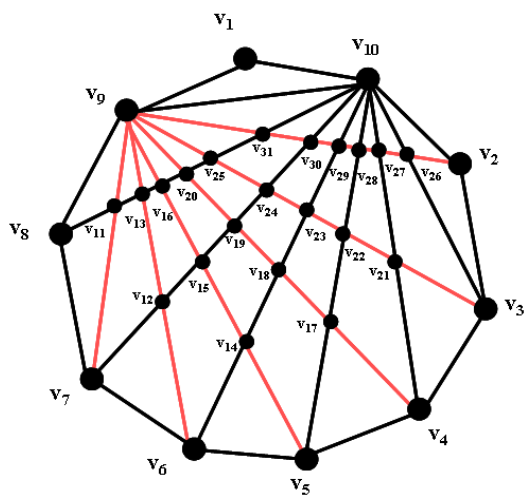

Рис. 19.4. Топологический рисунок внутри гамильтонова цикла.

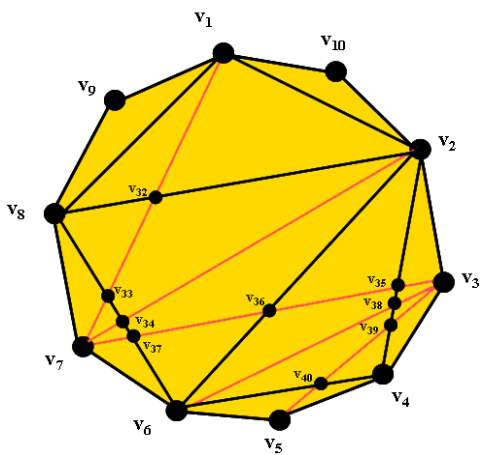

Рис. 19.5. Построение соединений с мнимыми вершинами для внешней части гамильтонова цикла.

Выберем непересекающееся подмножество соединений для внешней части гамильтонова цикла (рис. 19.5).

Построим топологический рисунок для внешней части гамильтонова цикла путем зеркального отображения (рис. 19.6).

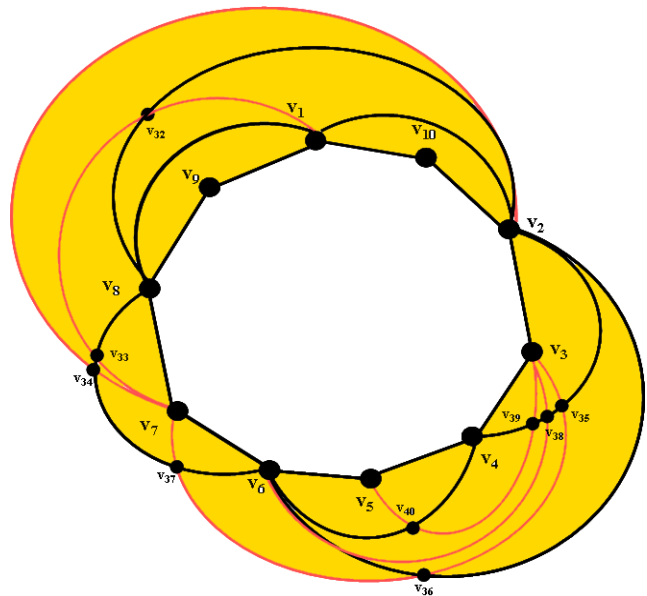

Рис. 19.6. Топологический рисунок внешний части гамильтонова цикла.

Общий топологический рисунок первого слоя, состоящий из внутренней и внешней части топологического пространства, представлен на рис. 19.7.

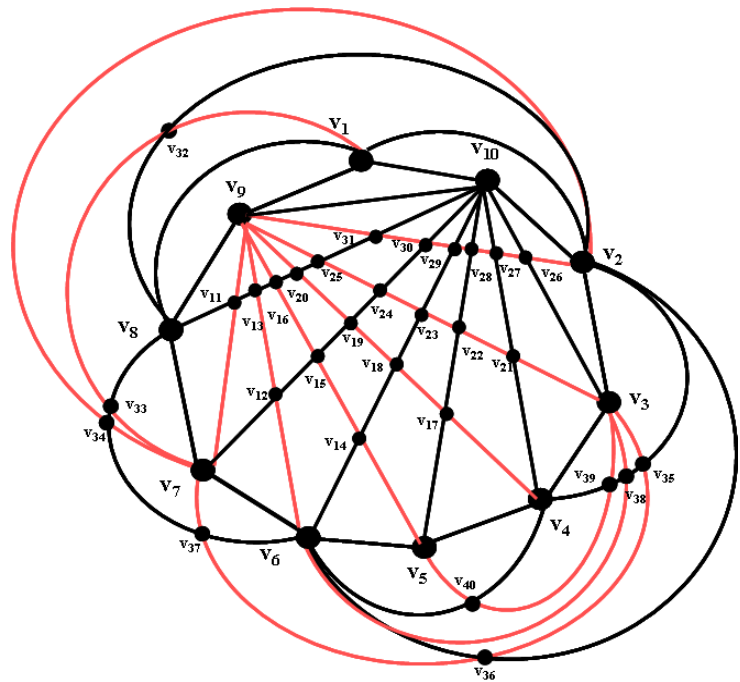

Рис. 19.7. Топологический рисунок первого слоя.

Построим топологический рисунок второго слоя для внутренней части гамильтонова цикла (рис. 19.8).

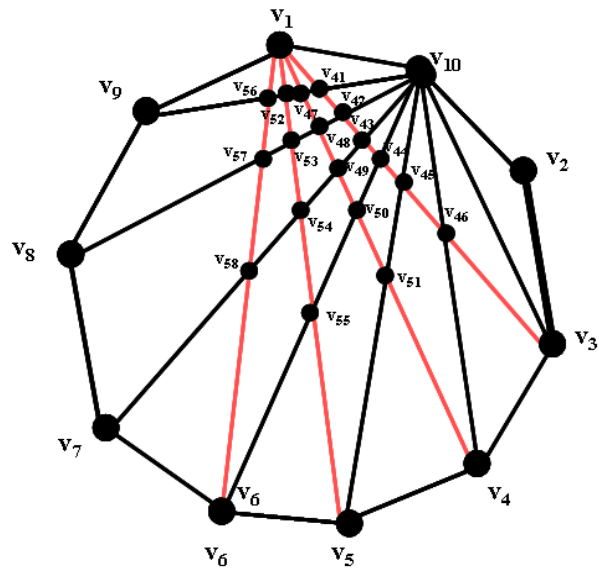

Рис. 19.8. Топологический рисунок второго слоя для
внутренней части гамильтонова цикла.

Построим топологический рисунок второго слоя для внешней части гамильтонова цикла (рис. 19.9).

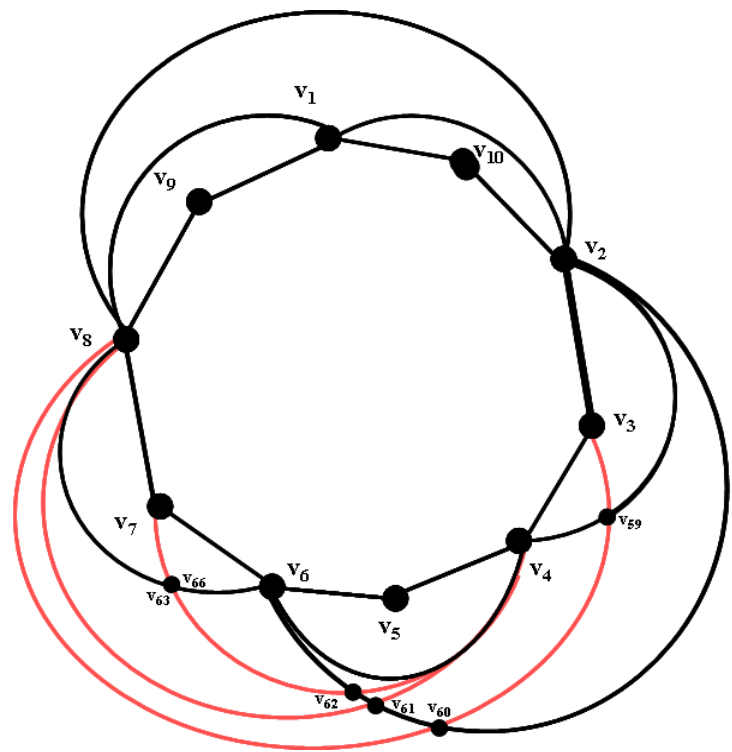

Рис. 19.9. Топологический рисунок второго слоя для внешней части гамильтонова цикла.

Топологический рисунок графа с толщиной для полного графа К$_{10}$ представляется тремя слоями. Первый слой – топологический рисунок максимально плоского суграфа (рис. 19.1). Топологический рисунок второго слоя представлен на рис. 19.7. Топологический рисунок третьего слоя представлен на рис. 19.11.

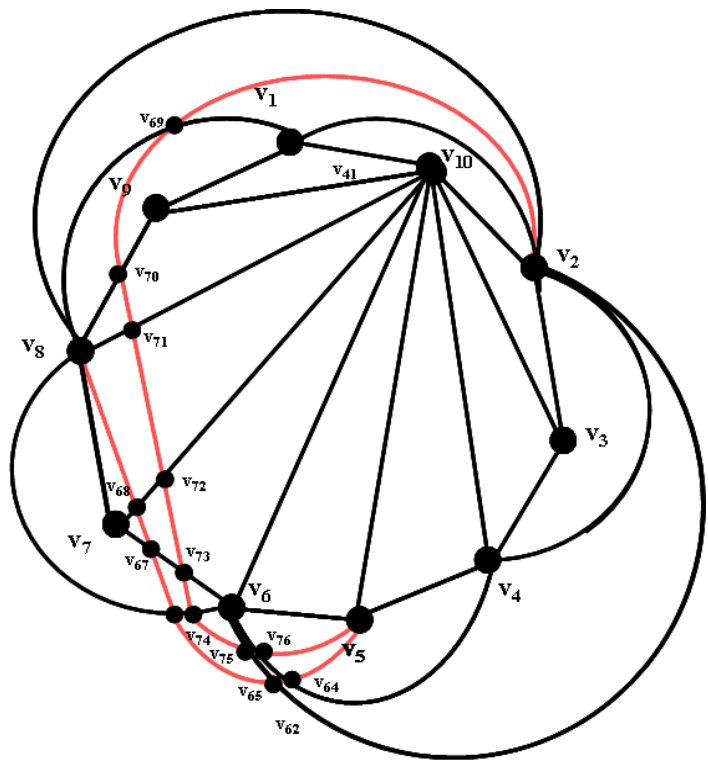

Рис. 19.10. Топологический рисунок канальных соединений.

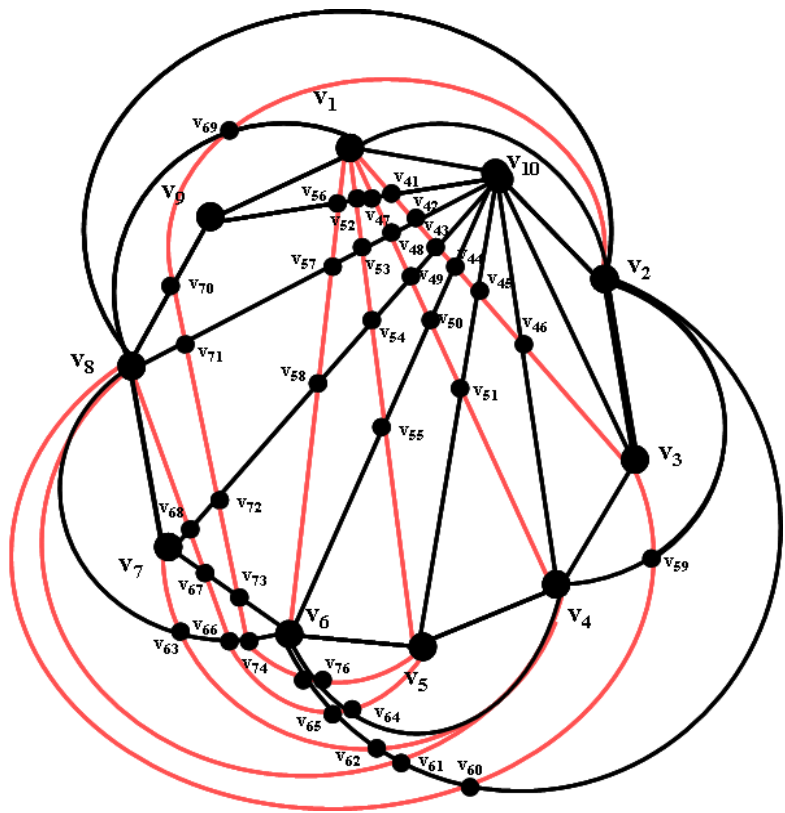

Рис. 19.11. Общий топологический рисунок второго слоя.

Топологический рисунок графа с толщиной для отображения третьего слоя на второй представлен на рис. 19.12.

Топологический рисунок каждого слоя можно представить без введения мнимых вершин приняв за основу гамильтонов цикл графа.

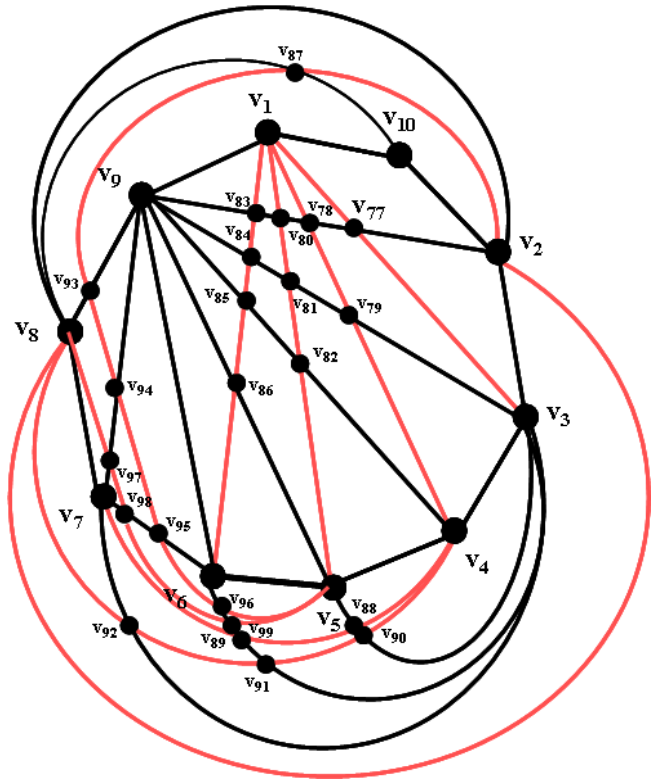

Рис. 19.12. Топологический рисунок проекции третьего слоя на второй слой.

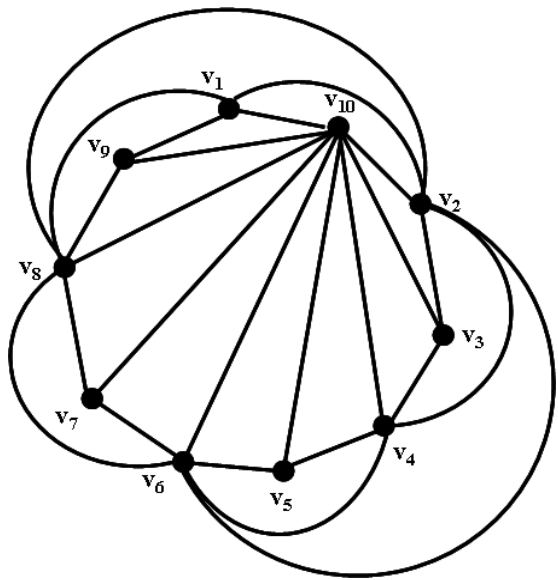

Рис. 19.13. Топологический рисунок первого слоя.

Первый слой топологического рисунка графа с толщиной характеризуется следующей системой циклов (рис. 19.13):

$c_1 = (v_8,v_2)+(v_2,v_1)+(v_1,v_8);$
$c_2 = (v_8,v_1)+(v_1,v_9)+(v_9,v_8);$
$c_3 = (v_1,v_2)+(v_2,v_{10})+(v_{10},v_1);$
$c_4 = (v_9,v_1)+(v_1,v_{10})+(v_{10},v_9);$
$c_5 = (v_9,v_{10})+(v_{10},v_8)+(v_8,v_9);$
$c_6 = (v_8,v_{10})+(v_{10},v_7)+(v_7,v_8);$
$c_7 = (v_7,v_{10})+(v_{10},v_6)+(v_6,v_7);$

$c_8 = (v_6,v_{10})+(v_{10},v_5)+(v_5,v_6)$;
$c_9 = (v_5,v_{10})+(v_{10},v_4)+(v_4,v_5)$;
$c_{10} = (v_4,v_{10})+(v_{10},v_3)+(v_3,v_4)$;
$c_{11} = (v_3,v_{10})+(v_{10},v_2)+(v_2,v_3)$;
$c_{12} = (v_8,v_7)+(v_7,v_6)+(v_6,v_8)$;
$c_{13} = (v_6,v_5)+(v_5,v_4)+(v_4,v_6)$;
$c_{14} = (v_4,v_3)+(v_3,v_2)+(v_2,v_4)$;
$c_{15} = (v_6,v_4)+(v_4,v_2)+(v_2,v_6)$;
Обод $c_0 = (v_6,v_8)+(v_8,v_2)+(v_2,v_6)$.

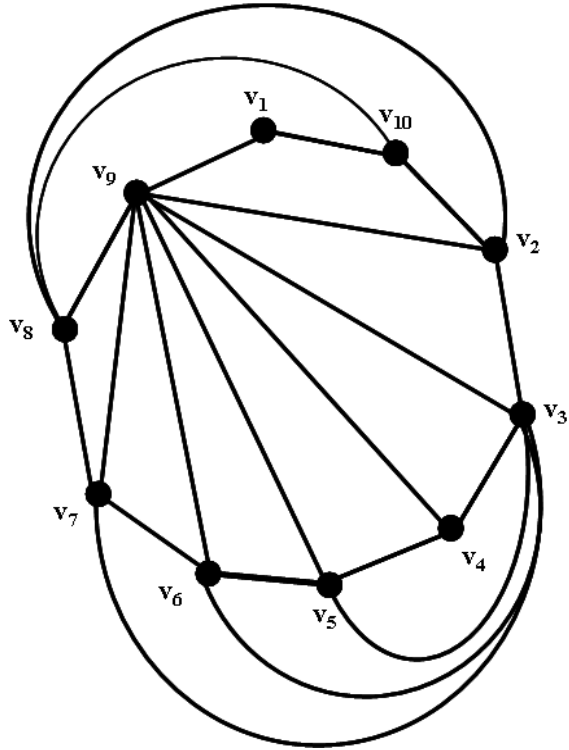

Рис. 19.14. Топологический рисунок второго слоя.

Второй слой топологического рисунка графа с толщиной характеризуется следующей системой циклов (рис. 19.14):

$c_{16} = (v_8,v_2)+(v_2,v_{10})+(v_{10},v_8)$;
$c_{17} = (v_8,v_{10})+(v_{10},v_1)+(v_1,v_9)+(v_9,v_8)$;
$c_{18} = (v_9,v_1)+(v_1,v_{10})+(v_{10},v_2)+(v_2,v_9)$;
$c_{19} = (v_9,v_2)+(v_2,v_3)+(v_3,v_9)$;
$c_{20} = (v_9,v_3)+(v_3,v_4)+(v_4,v_9)$;
$c_{21} = (v_9,v_4)+(v_4,v_5)+(v_5,v_9)$;
$c_{22} = (v_9,v_5)+(v_5,v_6)+(v_6,v_9)$;
$c_{23} = (v_9,v_6)+(v_6,v_7)+(v_7,v_9)$;
$c_{24} = (v_9,v_7)+(v_7,v_8)+(v_8,v_9)$;
$c_{25} = (v_7,v_6)+(v_6,v_3)+(v_3,v_7)$;
$c_{26} = (v_6,v_5)+(v_5,v_3)+(v_3,v_6)$;
$c_{28} = (v_5,v_4)+(v_4,v_3)+(v_3,v_5)$.
Обод $c_0 = (v_2,v_8)+(v_8,v_7)+(v_7,v_3)+(v_3,v_2)$.

Третий слой топологического рисунка графа с толщиной характеризуется следующей системой циклов (рис. 19.15):

$c_{29} = (v_{93},v_2)+(v_2,v_{10})+(v_{10},v_1)+(v_1,v_9)+(v_9,v_{93})$;

$c_{30} = (v_8,v_{98})+(v_{98},v_7)+(v_7,v_8)$;
$c_{31} = (v_{98},v_8)+(v_8,v_{93})+(v_{93},v_{95})+(v_{95},v_{98})$;
$c_{32} = (v_{93},v_9)+(v_9,v_1)+(v_1,v_6)+(v_6,v_{95})+(v_{95},v_{93})$;
$c_{33} = (v_6,v_1)+(v_1,v_5)+(v_5,v_6)$;
$c_{34} = (v_5,v_1)+(v_1,v_4)+(v_4,v_5)$;
$c_{35} = (v_4,v_1)+(v_1,v_3)+(v_3,v_4)$;
$c_{36} = (v_3,v_1)+(v_1,v_{10})+(v_{10},v_2)+(v_2,v_3)$;
$c_{37} = (v_8,v_4)+(v_4,v_3)+(v_3,v_2)+(v_2,v_8)$;
$c_{38} = (v_8,v_7)+(v_7,v_4)+(v_4,v_8)$;
$c_{39} = (v_7,v_{98})+(v_{98},v_5)+(v_5,v_4)+(v_4,v_7$;
$c_{40} = (v_{98},v_{95})+(v_{95},v_5)+(v_5,v_{98})$.
$c_{41} = (v_{95},v_6)+(v_6,v_5)+(v_5,v_{95})$;
Обод $c_0 = (v_2,v_{93})+(v_{93},v_8)+(v_8,v_2)$.

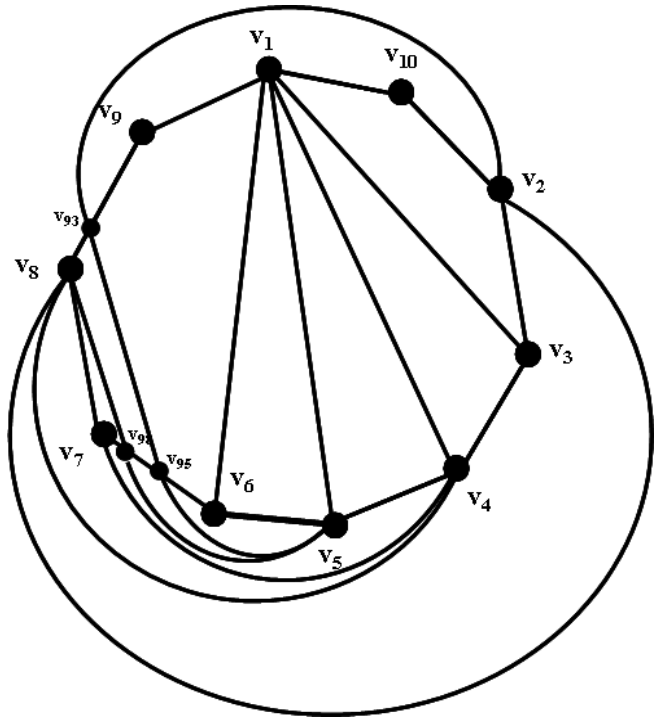

Рис. 19.15. Топологический рисунок третьего слоя.

**Комментарии**

В данной главе рассматриваются вопросы разбиения топологического пространства на части с целью определения непересекаюшихся соединений. Гамильтонов цикл максимального плоского суграфа разделяет топологическое пространство на внешнюю и внутреннею части. Для определения подмножества непересекающихся соединений используются методы векторной алгебры пересечений. Выделенные непересекающиеся соединения каждой части образуют пересечения с ребрами максимально плоского суграфа и характеризуются введением мнимых вершин и построением новой системы независимых циклов. Объединение систем циклов характеризующих внешнюю и внутреннюю части топологического пространства порождает топологический рисунок слоя и позволяет выделить каналы для проведения соединений (рис. 19.10).

**Заключение**

В представленной работе предложен метод построения топологического рисунка непланарного несепарабельного графа. Метод основан на разделении общего топологического рисунка непланарного графа на совокупность плоских топологических рисунков отдельных частей (слоев) графа. Выделение максимально плоского суграфа определяет первый плоский слой, а остальные плоские топологические рисунки строятся как проекции на первый слой. Таким образом, топологический рисунок максимально плоского суграфа является координатно-базисной системой для построения топологических рисунков других слоев.

В работе рассмотрено построение топологических рисунков на примере графов $K_7$, $K_8$ и $K_{10}$. Предполагается возможность применения метода к произвольному несепарабельному непланарному графу.